\documentclass[11pt]{amsart}
\usepackage{amsthm,amssymb,amsmath}
\usepackage{euscript}
\usepackage[matrix,arrow,curve]{xy}
\usepackage{graphics}
\usepackage{color}

\numberwithin{equation}{section}

\theoremstyle{plain}

\newtheorem{thm}{Theorem}[section]
\newtheorem{theorem}[thm]{Theorem}
\newtheorem*{theorem*}{Theorem}
\newtheorem{cor}{Corollary}[section]
\newtheorem{corollary}[cor]{Corollary}
\newtheorem{lemma}{Lemma}[section]

\newtheorem{prop}{Proposition}[section]
\newtheorem{proposition}[prop]{Proposition}

\newtheorem{definition}{Definition}[section]
\newtheorem{remark}{Remark}[section]

\newtheorem{example}{Example}[section]

\newcommand{\R}{\mathbb{R}}
\newcommand{\Z}{\mathbb{Z}}
\newcommand{\C}{\mathbb{C}}

\newcommand{\half}{{\textstyle\frac{1}{2}}}

\newcommand{\iso}{\cong}           
\newcommand{\htp}{\simeq}          
\newcommand{\smooth}{C^\infty}

\newcommand{\RP}[1]{\R {\mathrm P}^{#1}}
\newcommand{\leftsc}{\langle}
\newcommand{\rightsc}{\rangle}

\newcommand{\Cone}{\mathrm{Cone}}

\newcommand{\A}{\EuScript{A}}
\newcommand{\B}{\EuScript{B}}

\newcommand{\E}{\EuScript{E}}
\newcommand{\F}{\EuScript{F}}

\newcommand{\JJ}{\EuScript{J}}
\newcommand{\PP}{\EuScript{P}}
\newcommand{\TT}{\EuScript{T}}
\newcommand{\UU}{\EuScript{U}}
\newcommand{\VV}{\EuScript{V}}

\newcommand{\K}{\mathbb{K}}
\newcommand{\scrM}{\EuScript{M}}
\newcommand{\scrR}{\EuScript{R}}
\newcommand{\scrS}{\EuScript{S}}
\renewcommand{\top}{\mathrm{top}}
\newcommand{\geom}{\mathrm{geom}}
\newcommand{\frakC}{\mathfrak{C}}

\begin{document}
\title[Lefschetz fibrations]{Fukaya $A_\infty$-structures associated to Lefschetz fibrations. I}
\author{Paul Seidel}
\maketitle

\section{Introduction}

This paper discusses Lagrangian Floer cohomology in the context of Lefschetz fibrations, with emphasis on the algebraic structures encountered there. Part of the intent is expository, which means we'll review structures that have appeared in the literature before, but take a point of view which differs a bit from the usual one. Going beyond that, we obtain a new comparison theorem, which corresponds to a weak version of \cite[Conjecture 4]{seidel06}. This is intended to be the start of a wider project concerning Fukaya categories of Lefschetz fibrations. A sneak peek at the possible further developments is afforded by the papers \cite{maydanskiy09,maydanskiy-seidel09}, which share the basic outlook with this one but take a shortcut through the theory, aiming directly for specific applications. Finally, there is a significant overlap in the eventual goals, even though not necessarily in the intermediate steps, with \cite{bee} (obviously, \cite{bee} encompasses much more than what we're aiming for).

Consider a symplectic Lefschetz fibration $\pi: E \rightarrow \C$, with fibre $M$. We fix a basis of vanishing paths in the base. This gives rise to a collection of Lefschetz thimbles, which are noncompact Lagrangian submanifolds in $E$, as well as the associated vanishing cycles, which are Lagrangian spheres in $M$. Using either one as the starting point, one can associate to this basis a directed $A_\infty$-algebra $\A$. In terms of $M$, one first constructs the full Fukaya $A_\infty$-algebra, and then extracts the directed subalgebra by forgetting part of the structure in a purely algebraic process, which may seem arbitrary.  However, when one thinks in terms of $E$, the directedness of $\A$ comes from geometric conventions which govern the treatment of the intersection points at infinity. In spite of this, it is not hard to see that the two constructions yield the same result. This has been long implicit in the literature, starting with Kontsevich's unpublished work, and appears more recently in \cite[Section 18e]{seidel04}.

As noted in \cite{seidel06}, $\A$ comes with an important additional datum, namely a map of $A_\infty$-bimodules over $\A$ from the dual diagonal bimodule $\A^\vee$ to the diagonal bimodule $\A$. From the point of view of $M$, this is again constructed algebraically from the full $A_\infty$-algebra, and represents the first level of information beyond that contained in $\A$. In terms of $E$, the same kind of datum can be obtained geometrically, by comparing two different ways of treating intersection points at infinity. Our main result is the following:

\begin{theorem*}
The two maps $\A^\vee \rightarrow \A$ obtained by working in terms of $M$ and $E$, respectively, agree up to homotopy and composition with a quasi-isomorphism from $\A^\vee$ to itself.
\end{theorem*}

Obviously, the precise statement depends on first defining all elements involved. We refer to Corollary \ref{th:main} later in the paper for this, and only briefly outline the general assumptions and conventions used. We will impose strong exactness assumptions on $M$ and $E$. These are probably mostly unnecessary, but any attempt to remove them would complicate the technical side of the story considerably. Next, all Floer cohomology groups are taken with coefficients in $\K = \Z/2$, and are ungraded. This is again for simplicity, but the issues are far less complicated than the previously mentioned one. The general assumptions and constructions needed to introduce signs and gradings are standard, and one could easily integrate them into the developments here, at the cost of lengthening the exposition somewhat.

{\em Acknowledgments}. I thank Claude Viterbo and Kai Cieliebak for inviting me to give lectures in which I inflicted early versions of this material on unsuspecting audiences. The work described here was partially supported by NSF grant DMS-0652620.

\section{Geometry of the fibre\label{sec:fibre}}

\subsection{Geometric data\label{subsec:gromov}}
We begin by summarizing the basic setup for the pseudo-holomorphic map equation and perturbations thereof. The terminology here follows \cite[Chapter 2]{seidel04}, but the same material is present in many places in the literature, going back to Floer's and Fukaya's work \cite{floer88c,floer88,fukaya93}.

\begin{definition} \label{def:exact-symplectic}
An {\em exact symplectic manifold} $(M,\omega_M,\theta_M,I_M)$ is a compact $2n$-dimensional manifold with boundary, equipped with an exact symplectic form $\omega_M = d\theta_M$, as well as with an $\omega_M$-compatible almost complex structure $I_M$ which is weakly convex (meaning that $I_M$-holomorphic maps from an open disc to $M$ can touch $\partial M$ only if they remain entirely inside it).
\end{definition}

From now on, we will work with a fixed such $M$. An exact Lagrangian submanifold is an $n$-dimensional submanifold $V \subset M$, disjoint from $\partial M$, such that $\theta_M|V$ is an exact one-form. All Lagrangian submanifolds occurring here are implicitly assumed to be connected.

\begin{definition} \label{def:riemann-1}
A {\em punctured-boundary Riemann surface} is a pair $(S,V)$ of the following kind. $S$ is a Riemann surface of the form $S = \bar{S} \setminus \Sigma$, where $\bar{S}$ is a compact connected Riemann surface with nonempty boundary, and $\Sigma$ is a nonempty finite set of boundary points. The boundary conditions $V = (V_z)$ are a locally constant family of exact Lagrangian submanifolds parametrized by points $z \in \partial S$ (or equivalently, the choice of one such submanifold for each connected component of the boundary).
\end{definition}

While this structure is conceptually satisfactory, it is technically convenient to rigidify it a little. First of all, we make an (arbitrary) distinction of the punctures into inputs and outputs, $\Sigma = \Sigma_{\mathrm{in}} \cup \Sigma_{\mathrm{out}}$. Secondly, we choose strip-like ends $\epsilon = (\epsilon_\zeta)$, which are proper holomorphic embeddings
\begin{equation}\left\{
\begin{aligned}
& \epsilon_\zeta: \R^+ \times [0,1] \longrightarrow S, \\
& \epsilon_\zeta^{-1}(\partial S) = \R^+ \times \{0,1\}, \\
& \textstyle \lim_{s \rightarrow +\infty} \epsilon_\zeta(s,t) = \zeta
\end{aligned}\right.
\end{equation}
for $\zeta \in \Sigma_{\mathrm{in}}$, respectively with $\R^- = (-\infty,0]$ instead of $\R^+ = [0,\infty)$ for $\zeta \in \Sigma_{\mathrm{out}}$. We also ask that the images of the $\epsilon_\zeta$ for different $\zeta$ should be pairwise disjoint. The resulting $(S,V,\epsilon)$ is called a {\em Riemann surface with strip-like ends}. Strip-like ends provide distinguished local coordinates near the punctures, which is useful for studying the asymptotic behaviour of pseudo-holomorphic maps, as well as in gluing processes.

\subsection{Pseudo-holomorphic maps\label{subsec:pseudo}}
We will use almost complex structures on $M$ which are $\omega_M$-compatible and agree with the given $I_M$ in some neighbourhood of $\partial M$. Denote the set of such almost complex structures by $\JJ(M)$.

Take a Riemann surface with strip-like ends. In addition, choose a family $J = (J_z)$, $z \in S$, of almost complex structures in $\JJ(M)$. On the strip-like ends the family should be translation-invariant, which means that $J_{\epsilon_\zeta(s,t)} = J_{\zeta,t}$ for some $J_\zeta = (J_{\zeta,t})$ depending only on $t \in [0,1]$. The associated pseudo-holomorphic map equation is
\begin{equation} \label{eq:gromov}
\left\{
\begin{aligned}
& u: S \longrightarrow M, \\
& u(z) \in V_z \text{ for all $z \in \partial S$,} \\
& (du)^{0,1} = \half(du(z) + J_z \circ du(z) \circ j_S) = 0,
\end{aligned}
\right.
\end{equation}
where $j_S$ is the complex structure on $S$. At each end $\zeta$ we have two Lagrangian submanifolds $V_{\zeta,k}$ ($k = 0,1$), which are those associated to the points $\epsilon_\zeta(s,k)$ for any $s$. The natural limit condition for \eqref{eq:gromov} is that as $s \rightarrow \pm\infty$, $u(\epsilon_\zeta(s,\cdot))$ converges to a constant $x_\zeta \in V_{\zeta,0} \cap V_{\zeta,1}$.

A familiar variation on the story above, which is particularly useful when the boundary conditions fail to intersect transversally, is to introduce a Hamiltonian inhomogeneous term. We will consider the simplest form of such terms. Namely, choose $H \in \smooth(M,\R)$ which vanishes near the boundary. Denote by $X_H$ its Hamiltonian vector field, and by $\phi_H$ the time-one map of its flow (the time $t$ map for general $t$ is of course just $\phi_{tH}$). Additionally, choose $\gamma \in \Omega^1(S)$ which is zero on the boundary (meaning that $\gamma|\partial S \in \Omega^1(\partial S)$ is zero), and which over the ends has the form $\epsilon_\zeta^*(\gamma) = g_\zeta(t)\, \mathit{dt}$ for some function $g_\zeta \in \smooth([0,1],\R)$. The perturbation of \eqref{eq:gromov} associated to $(H,\gamma)$ is
\begin{equation} \label{eq:gromov-2}
(du - X_H \otimes \gamma)^{0,1} = 0,
\end{equation}
where $X_H \otimes \gamma$ is thought of as a section of $\mathit{TM} \otimes \mathit{TS}^* = \mathit{Hom}(\mathit{TS},\mathit{TM}) \rightarrow S \times M$, and evaluated at $(z,u(z))$. This time, the natural limits are of the form $\lim_{s \rightarrow \pm\infty} u(\epsilon_\zeta(s,t)) = u_\zeta(t)$ for
\begin{equation} \label{eq:gromov-2-limits}
\left\{
\begin{aligned}
& u_\zeta: [0,1] \longrightarrow M, \\
& du_\zeta/dt = g_\zeta(t) X_H, \\
& u_\zeta(k) \in V_{\zeta,k} \text{ for $k = 0,1$.}
\end{aligned}
\right.
\end{equation}
Clearly, such flow lines correspond bijectively to points $x_\zeta = u_\zeta(1) \in \phi_{c_\zeta H}(V_{\zeta,0}) \cap V_{\zeta,1}$, where $c_\zeta = \int_0^1 g_\zeta(t) \mathit{dt}$.

\begin{remark} \label{th:r1}
The class of perturbations used in \eqref{eq:gromov-2} is by no means the most general one, compare \cite[Chapter 2]{seidel04} or the treatment of pseudo-holomorphic maps in Lefschetz fibrations later on \eqref{eq:u-equation}. However, it is particularly simple from an expository viewpoint, and sufficient for the purpose at hand; the only drawback is that it requires a little more care when discussing transversality, see Remark \ref{th:constant-maps} below.
\end{remark}

\subsection{An $A_\infty$-algebra\label{subsec:directed-algebra}}
Recall first that for any pair $(V_0,V_1)$ of exact Lagrangian submanifolds we have a Floer cohomology group $\mathit{HF}(V_0,V_1)$, which in this very stripped-down framework is a single vector space over the coefficient field $\K = \Z/2$. It is invariant under isotopies of either $V_0$ or $V_1$ (within the same class of submanifolds). Moreover, for $V_0 = V_1 = V$ we have $\mathit{HF}(V,V) \iso H^*(V;\K)$ canonically. Our interest is in the cochain level algebraic structures underlying Floer cohomology. For the moment, we restrict ourselves to the most elementary instance of such structures, namely directed Fukaya $A_\infty$-algebras. This means that we consider an ordered collection of Lagrangian submanifolds $(V_1,\dots,V_m)$ in $M$, and define an $A_\infty$-algebra $\bar\A$ whose cohomology is the direct sum of Floer cohomology groups
\begin{equation} \label{eq:directed-cohomology}
H(\bar{\A}) = \bigoplus_{i<j} \mathit{HF}(V_i,V_j).
\end{equation}
The directedness property is technically useful in excluding self-transversality problems, which means that we can use the unperturbed equation \eqref{eq:gromov}. More importantly, directedness has an inherent meaning in terms of Lefschetz fibrations, but that will only become clear later.

From now on, we assume that the $V_i$ are in general position in the following sense:
\begin{equation} \label{eq:generic-1}
\parbox{28em}{The intersections $V_i \cap V_j$, $i < j$, are transverse, and the triple intersections $V_i \cap V_j \cap V_k$, $i<j<k$, are empty.}
\end{equation}
For each $i<j$, let $\mathit{CF}(V_i,V_j)$ be the Floer cochain group in its simplest version, which is the vector space over $\K$ freely generated by the points of $V_i \cap V_j$. Choose a generic one-parameter family of almost complex structures $J_{i,j}^1 = (J_{i,j,t}^1)_{0 \leq t \leq 1}$ in $\JJ(M)$. Consider the pseudo-holomorphic strip equation, which is the special case of \eqref{eq:gromov} given by $S = \R \times [0,1]$, $V_{s,0} = V_i$, $V_{s,1} = V_j$, and $J_{s,t} = J_{i,j,t}^1$. Given intersection points $x_0,x_1 \in V_i \cap V_j$, let $n^1(x_0,x_1) \in \K$ be the (modulo $2$) number of solutions of this equation which are non-stationary and isolated mod translation, with limits $x_0$ as $s \rightarrow -\infty$ and $x_1$ as $s \rightarrow +\infty$. The Floer differential on $\mathit{CF}(V_i,V_j)$ is defined in the standard way by
\begin{equation} \label{eq:mu1}
\mu^1(x_1) = \sum_{x_0} n^1(x_0,x_1) \, x_0,
\end{equation}
and its cohomology is $\mathit{HF}(V_i,V_j)$. In particular, if we take the direct sum
\begin{equation} \label{eq:fibre-a}
\bar{\A} = \bigoplus_{i<j} \mathit{CF}(V_i,V_j),
\end{equation}
then its cohomology is \eqref{eq:directed-cohomology}. Consider the semisimple ring $R = \K^m = \K e_1 \oplus \cdots \oplus \K e_m$, where $e_i^2 = e_i$, and $e_ie_j = 0$ for $i \neq j$. One can view $\bar{\A}$ as an $R$-bimodule, where $e_k$ acts on the left by projecting to the summand with $j = k$, and on the right by projecting to the summand with $i = k$. Fukaya's construction turns $\bar{\A}$ into a (non-unital and ungraded) $A_\infty$-algebra over $R$. The $A_\infty$-structure is given by a sequence of maps
\begin{equation}
\mu^d: \bar{\A}^{\otimes d} \longrightarrow \bar{\A}, \quad d \geq 1.
\end{equation}
Here, the tensor product is taken over $R$, and the map itself respects the $R$-bimodule structure (from now on, this will be the standing convention for all tensor products and maps involving $R$-bimodules). When we unpack the notation, the components of $\mu^d$ are maps
\begin{equation} \label{eq:a}
\begin{aligned}
& \mathit{CF}(V_{i_{d-1}},V_{i_d}) \otimes \cdots \otimes \mathit{CF}(V_{i_0},V_{i_1}) \longrightarrow \mathit{CF}(V_{i_0},V_{i_d}), \\
& \mu^d(x_d,\dots,x_1) = \sum_{x_0} n^d(x_0,\dots,x_d)\, x_0
\end{aligned}
\end{equation}
where $i_0 < \cdots < i_d$, $x_0 \in V_{i_0} \cap V_{i_d}$, and $x_k \in V_{i_{k-1}} \cap V_{i_k}$ for $k>0$. For $d = 1$ this is the previously defined differential.

The coefficients in \eqref{eq:a} for $d \geq 2$ are obtained by counting pseudo-holomorphic polygons. Concretely, we consider surfaces $S$ which are $(d+1)$-punctured discs, which means that $\bar{S}$ is a closed disc, and $|\Sigma_{\mathrm{out}}| = 1$, $|\Sigma_{\mathrm{in}}| = d$. Write $\Sigma = \{\zeta_0,\dots,\zeta_d\}$, where $\zeta_0$ is the unique point of $\Sigma_{\mathrm{out}}$, and the other points are ordered in accordance with the boundary orientation. Boundary conditions are $V = (V_{i_0},\dots,V_{i_d})$, where $V_{i_0}$ is associated to the boundary component between $\zeta_0$ and $\zeta_1$, and so on. For each such surface we choose a set of strip-like ends, as well as a family $J$ of almost complex structures, which over the end reduces to the previously chosen $J_{i_0,i_d}^1$ and $J_{i_{k-1},i_k}^1$, $k \geq 1$. The last-mentioned piece of data can be more conveniently thought of as a family of almost complex structures $J_{i_0,\dots,i_d}^d$ parametrized by points of the universal family $\scrS^d \rightarrow \scrR^d$ of $(d+1)$-pointed discs. Besides the previously mentioned condition, there are more complicated consistency requirements expressed in terms of the Deligne-Mumford compactification $\bar\scrR^d$. We refer to \cite[Section 9]{seidel04} for the details, and assume that a generic choice of almost complex structures, subject to these requirements, has been made. Given $(i_0,\dots,i_d)$ and $(x_0,\dots,x_d)$ as in \eqref{eq:a}, one considers the moduli space $\scrM^d(x_0,\dots,x_d)$ of pairs $(S,u)$, where $S$ is a $(d+1)$-punctured disc, equipped with the previously specified boundary conditions and almost complex structures, and $u$ is a solution of \eqref{eq:gromov} with limits $x_i$. Genericity of the almost complex structures, even within the overall allowed class, means that this moduli space is smooth and its zero-dimensional part is a finite set. One defines $n^d(x_0,\dots,x_d)$ by counting points in that set.

Later on in this section, we carry out several more constructions of the same kind. To free up some notation, we now introduce the convention that all structures relating to $\bar{\A}$ will have that added as a subscript. Here is a list for future reference: \vskip1em

\begin{center}
\begin{tabular}{|l|l|}
\hline
$\mathit{CF}_{\bar{\A}}(V_i,V_j)$ & unperturbed Floer cochain spaces, defined for $i<j$ \\
$\mu^d_{\bar{\A}}$ & $A_\infty$-operations on $\bar{\A}$ \\
$\scrS^d_{\bar{\A}} \rightarrow \scrR^d_{\bar{\A}}$ & moduli space of $(d+1)$-punctured discs \\
$J^d_{\bar{\A},i_0,\dots,i_d}$ & families of almost complex structures \\
$\scrM_{\bar{\A}}^d(x_0,\dots,x_d)$ & moduli spaces of pseudo-holomorphic maps \\
$n_{\bar{\A}}^d(x_0,\dots,x_d)$ & numbers counting points in the spaces above \\
\hline
\end{tabular}
\end{center} \vskip1em

Finally, we'd like to consider a minor variant. By construction, $\bar\A$ is non-unital. The standard way of addressing this shortcoming is to adjoin strict units, meaning to pass to $\A = R \oplus \bar{\A}$ with the $A_\infty$-structure extended as follows:
\begin{equation} \label{eq:adjoin-units}
\begin{aligned}
& \mu^1_\A(e_i) = 0, \\
& \mu^2_\A(e_i,a) = e_i a, \;\; \mu^2_\A(a,e_i) = a e_i, \\
& \mu^d_\A(\dots,e_i,\dots) = 0 \text{ for all $d > 2$}. \\
\end{aligned}
\end{equation}

\subsection{A bimodule\label{subsec:b-bimodule}}
Think for a moment of general $(d+1)$-punctured discs with boundary components (in their natural order) labeled by $V_{i_0},\dots,V_{i_d}$. By definition, the $A_\infty$-structure on $\bar{\A}$ involves only a subset of possible boundary conditions, namely ones where $i_0 < \cdots < i_d$. We will now consider the next step after that, which is when $i_k < i_{k+1}$ for all $k$ except one. After introducing a suitable Hamiltonian perturbation, the resulting information can be encoded algebraically into a bimodule $\B$ over $\bar{\A}$, whose cohomology is
\begin{equation} \label{eq:b-cohomology}
H(\B) = \bigoplus_{i,j} \mathit{HF}(V_i,V_j).
\end{equation}
For generalities about $A_\infty$-bimodules, see for instance \cite{tradler01,seidel08}.

Choose $H$ as in Section \ref{subsec:pseudo}. By making a generic choice, we can achieve general position in the following sense:
\begin{equation} \label{eq:b-general-position}
\parbox{28em}{$\phi_{H}(V_i) \cap V_j$ is transverse for any $i,j$. Moreover, the triple intersections $\phi_{H}(V_i) \cap \phi_{H}(V_j) \cap V_k$ and $V_i \cap V_j \cap \phi_{H}(V_k)$, for $i<j$ and arbitrary $k$, are empty.
}
\end{equation}
This time, write $\mathit{CF}(V_i,V_j)$ for the $\K$-vector space freely generated by points of $\phi_{H}(V_i) \cap V_j$. The direct sum of these groups, without any ordering condition, forms the $R$-bimodule
\begin{equation} \label{eq:b}
\B = \bigoplus_{i,j} \mathit{CF}(V_i,V_j).
\end{equation}
Fix a cutoff function $g \in \smooth([0,1],\R)$ vanishing near $t = 0$ and $t = 1$, and such that $\int_0^1 g(t)\, \mathit{dt} = 1$. Additionally, for each $(i,j)$ in $\{1,\dots,m\}$ choose a family $J^{0|1|0}_{i,j} = (J_{i,j,t}^{0|1|0})_{0 \leq t \leq 1}$ in $\JJ(M)$. Consider the perturbed pseudo-holomorphic strip equation, which is the special case of \eqref{eq:gromov-2} with $S = \R \times [0,1]$, $V_{s,0} = V_i$, $V_{s,1} = V_j$, $J_{s,t} = J_{i,j,t}^{0|1|0}$, and $\gamma = g(t)\,\mathit{dt}$. The limits \eqref{eq:gromov-2-limits} over the ends correspond bijectively to points of $\phi_H(V_i) \cap V_j$. In the standard way \eqref{eq:mu1}, we use the numbers $n^{0|1|0}(x_0,x_1)$ obtained by counting solutions of this equation (mod translation) to define a differential on $\B$, denoted by $\mu^{0|1|0}$. Its cohomology is again Floer cohomology, explaining \eqref{eq:b-cohomology}.

More generally, the $\bar{\A}$-bimodule structure of $\B$ consists of maps
\begin{equation} \label{eq:b-structure}
\mu^{q|1|p}: \bar{\A}^{\otimes q} \otimes \B \otimes \bar{\A}^{\otimes p} \longrightarrow \B, \quad p,q \geq 0.
\end{equation}
To define these maps for $p+q>0$, one proceeds as follows. Take a $(p+q+2)$-punctured disc $S$, equipped with boundary conditions $V = (V_{i_0},\dots,V_{i_{p+q+1}})$ where $i_0 < \cdots < i_p$ and $i_{p+1} < \cdots < i_{p+q+1}$. It will be convenient for us to think of $S$ as a punctured strip, which means to identify
\begin{equation} \label{eq:punctured-strip}
S = (\R \times [0,1]) \setminus \{\zeta_1,\dots,\zeta_p,\zeta_{p+2},\dots,\zeta_{p+q+1}\}
\end{equation}
where $\zeta_1,\dots,\zeta_p$ lie on $\R \times \{0\}$, and $\zeta_{p+2},\dots,\zeta_{p+q+1}$ on $\R \times \{1\}$ (thus, $\zeta_0$ corresponds to the end $s = -\infty$, and $\zeta_{p+1}$ to $s = +\infty$). Choose strip-like ends $\epsilon$ for our surface, so that $\epsilon_0$ and $\epsilon_{p+1}$ are (up to a translation in $s$-direction) the obvious local coordinates near the respective ends of $\R \times [0,1]$, and the other ends are disjoint from the support of the function $(s,t) \mapsto g(t)$. Suppose that we have a family of almost complex structures $J = (J_{s,t})$ which for $\pm s \ll 0$ reduces to $J_{i_0,i_{p+q+1}}^{0|1|0}$ and $J_{i_p,i_{p+1}}^{0|1|0}$, respectively, while being given over all the other $\epsilon_k$ by $J_{\bar{\A},i_k,i_{k+1}}^{1}$. We then consider the same inhomogeneous $\bar\partial$-equation \eqref{eq:gromov-2} as before, for maps $u: S \rightarrow M$ with boundary conditions as in \eqref{eq:gromov}. Limits as $s \rightarrow \mp \infty$ correspond to points $x_0 \in \phi_H(V_{i_0}) \cap V_{i_{p+q+1}}$ and $x_{p+1} \in \phi_H(V_{i_p}) \cap V_{i_{p+1}}$, while over the other ends they are just constant intersection points $x_k \in V_{i_k} \cap V_{i_{k+1}}$.

One has to consider the universal families of Riemann surfaces $S$ up to isomorphism, which means dividing by translation in \eqref{eq:punctured-strip}. Denote these families by $\scrS^{q|1|p} \rightarrow \scrR^{q|1|p}$, and choose families of almost complex structures $J_{i_0,\dots,i_{p+q+1}}^{q|1|p}$ parametrized by their total spaces. These are related to each other and to the previous choices made for $\bar{\A}$ by appropriate consistency conditions. More precisely, one can first assume that the choices for $\bar{\A}$ as well as the $J^{0|1|0}$ have been fixed. Then the consistency conditions for the higher $J^{q|1|p}$ are inductive with respect to $p+q$ (alternatively, one could solve the consistency conditions for both $\bar{\A}$ and the bimodule structure together by common induction on the number of marked points, but that would be less in tune with the way we are gradually building up the theory here). Generically within the class of consistent choices, the associated moduli spaces $\scrM^{q|1|p}(x_0,\dots,x_{p+q+1})$ yield numbers $n^{q|1|p}(x_0,\dots,x_{p+q+1})$, which as in \eqref{eq:a} appear as coefficients of the structure maps \eqref{eq:b-structure}. Figure \ref{fig:bimodule} shows the three basic ways in which the Riemann surfaces split, corresponding to the terms in the bimodule structure formula
\begin{equation}
\begin{aligned}
& \sum_{r,s} \mu^{q-s|1|p-r}(a_{p+q},\dots,a_{p+s+1},\mu^{s|1|r}(a_{p+s},\dots,a_{p+1},b,a_p,\dots,a_{p-r+1}), \\[-1em] & \hspace{25em} a_{p-r},
\dots,a_1) \\ = & \sum_{i,r} \mu^{q|1|p-r+1}(a_{p+q},\dots,a_{p+1},b,a_p,\dots,a_{i+r+1},
\\[-1em] & \hspace{16em}
\mu_{\bar\A}^r(a_{i+r},\dots,a_{i+1}),a_i,\dots,a_1) \\
+ & \sum_{j,s} \mu^{q-s+1|1|p}(a_{p+q},\dots,\mu_{\bar\A}^s(a_{p+j+s},\dots,a_{p+j+1}),
\\[-1em] & \hspace{18em} a_{p+j},\dots,
a_{p+1},b,a_p,\dots,a_1).
\end{aligned}
\end{equation}
We again conclude the construction by adding $\B$ to the relevant notation, as follows: \vskip1em
\begin{center}
\begin{tabular}{|l|l|}
\hline
$H_\B$ & Hamiltonian \\
$g_\B$ & cut-off function \\
$\mathit{CF}_\B(V_i,V_j)$ & $H_\B$-perturbed Floer cochain spaces, \\
& defined for all $i,j$ \\
$\mu^{q|1|p}_\B$ & $A_\infty$-bimodule operations on $\B$ \\
$\scrS^{q|1|p}_\B \rightarrow \scrR^{q|1|p}_\B$ & moduli space of $(p+q+2)$-punctured discs
\\ & \hspace{2em} (these are isomorphic to $\scrR_{\bar\A}^{p+1+q}$, but \\ & \hspace{2em} used differently) \\
$J^{q|1|p}_{\B,i_0,\dots,i_{p+q+1}}$ & families of almost complex structures \\
$\scrM_\B^{q|1|p}(x_0,\dots,x_{p+q+1})$ & moduli spaces of perturbed \\ & \hspace{2em} pseudo-holomorphic maps \\
$n_\B^{q|1|p}(x_0,\dots,x_{p+q+1})$ & numbers counting points in the spaces above \\
\hline
\end{tabular}
\end{center} \vskip1em
\begin{figure}
\begin{center}
\begin{picture}(0,0)%
\includegraphics{bimodule.pstex}%
\end{picture}%
\setlength{\unitlength}{3158sp}%
\begingroup\makeatletter\ifx\SetFigFont\undefined%
\gdef\SetFigFont#1#2#3#4#5{%
  \reset@font\fontsize{#1}{#2pt}%
  \fontfamily{#3}\fontseries{#4}\fontshape{#5}%
  \selectfont}%
\fi\endgroup%
\begin{picture}(6027,4694)(-2489,-3408)
\put(-100,-2391){\makebox(0,0)[lb]{\smash{{\SetFigFont{10}{12}{\rmdefault}{\mddefault}{\updefault}{\color[rgb]{0,0,0}$a_{p+j+s}$}%
}}}}
\put(1801,-1561){\makebox(0,0)[lb]{\smash{{\SetFigFont{10}{12}{\rmdefault}{\mddefault}{\updefault}{\color[rgb]{0,0,0}$a_{i+r}$}%
}}}}
\put(1351,-2086){\makebox(0,0)[lb]{\smash{{\SetFigFont{10}{12}{\rmdefault}{\mddefault}{\updefault}{\color[rgb]{0,0,0}$a_{p+j+1}$}%
}}}}
\put(-2024,-1936){\makebox(0,0)[lb]{\smash{{\SetFigFont{10}{12}{\rmdefault}{\mddefault}{\updefault}{\color[rgb]{0,0,0}$a_1$}%
}}}}
\put(-1274,-1936){\makebox(0,0)[lb]{\smash{{\SetFigFont{10}{12}{\rmdefault}{\mddefault}{\updefault}{\color[rgb]{0,0,0}$a_p$}%
}}}}
\put(-1124,-661){\makebox(0,0)[lb]{\smash{{\SetFigFont{10}{12}{\rmdefault}{\mddefault}{\updefault}{\color[rgb]{0,0,0}$a_{p+1}$}%
}}}}
\put(-2444,-661){\makebox(0,0)[lb]{\smash{{\SetFigFont{10}{12}{\rmdefault}{\mddefault}{\updefault}{\color[rgb]{0,0,0}$a_{p+q}$}%
}}}}
\put(-824,-1336){\makebox(0,0)[lb]{\smash{{\SetFigFont{10}{12}{\rmdefault}{\mddefault}{\updefault}{\color[rgb]{0,0,0}$b$}%
}}}}
\put(1426,164){\makebox(0,0)[lb]{\smash{{\SetFigFont{10}{12}{\rmdefault}{\mddefault}{\updefault}{\color[rgb]{0,0,0}$a_{p-r}$}%
}}}}
\put(2551,164){\makebox(0,0)[lb]{\smash{{\SetFigFont{10}{12}{\rmdefault}{\mddefault}{\updefault}{\color[rgb]{0,0,0}$a_{p-r+1}$}%
}}}}
\put(2226,1139){\makebox(0,0)[lb]{\smash{{\SetFigFont{10}{12}{\rmdefault}{\mddefault}{\updefault}{\color[rgb]{0,0,0}$a_{p+s}$}%
}}}}
\put(1426,1139){\makebox(0,0)[lb]{\smash{{\SetFigFont{10}{12}{\rmdefault}{\mddefault}{\updefault}{\color[rgb]{0,0,0}$a_{p+s+1}$}%
}}}}
\put(751,-1561){\makebox(0,0)[lb]{\smash{{\SetFigFont{10}{12}{\rmdefault}{\mddefault}{\updefault}{\color[rgb]{0,0,0}$a_{i+1}$}%
}}}}
\end{picture}%
\caption{\label{fig:bimodule}}
\end{center}
\end{figure}

Because of the directed nature of $\bar\A$, every bimodule over it admits certain natural submodules. We apply this idea to $\B$, as follows. For each $i$ take a cocycle $u_i \in \mathit{CF}_\B(V_i,V_i)$
which represents the unit class $1 \in H^*(V_i;\K) \iso \mathit{HF}(V_i,V_i)$. Then,
\begin{equation} \label{eq:bplus}
\B^+ \stackrel{\mathrm{def}}{=} \bigoplus_i \K u_i \oplus \bigoplus_{i < j} \mathit{CF}_\B(V_i,V_j)
\subset \B
\end{equation}
is automatically a sub-bimodule. Note that $\B^+$ is independent of the choice of representatives $u_i$ up to quasi-isomorphism. Namely, suppose that we have two different choices $u_i$ and $u_i'$, and the resulting bimodules $\B^+,\, (\B')^+ \subset \B$. By including both choices as well as bounding cochains $w_i$, satisfying $u_i' - u_i = \mu^{0|1|0}_\B(w_i)$, one gets another submodule of $\B$, which contains both $\B^+$ and $(\B')^+$ as quasi-isomorphic submodules.

Write $\B^- = \B/\B^+$. To the short exact sequence of bimodules
\begin{equation} \label{eq:b-plus-minus}
0 \rightarrow \B^+ \longrightarrow \B \longrightarrow \B^- \rightarrow 0
\end{equation}
one can associate a boundary map, which is a bimodule homomorphism
\begin{equation} \label{eq:boundary}
\Delta: \B^- \longrightarrow \B^+,
\end{equation}
unique up to homotopy. To express this in more abstract terms, take the differential category $\frakC$ of all bimodules over $\bar\A$. The cohomological category $H(\frakC)$, which is the homotopy category of bimodules, carries a natural triangulated structure, and the homotopy class $[\Delta]$ completes the maps from \eqref{eq:b-plus-minus} to an exact triangle.

\begin{remark} \label{th:ext}
From the definition \eqref{eq:bplus}, it is obvious that $H(\B^+) \rightarrow H(\B)$ is injective, which implies that the cohomology level map $H(\Delta)$, induced by the linear term $\Delta^{0|1|0}$, is zero. In spite of that, the whole of $\Delta$ still has an effect which is visible on the cohomology level. Namely, suppose that we choose our chain level representative so that $\Delta^{0|1|0} = 0$, which is always possible (by starting with a splitting of \eqref{eq:b-plus-minus} which is compatible with the differential). Then the next order terms $\Delta^{1|1|0}: \bar{\A} \otimes \B^- \rightarrow \B^+$ and $\Delta^{0|1|1}: \B^- \otimes \bar{\A} \rightarrow \B^+$ are chain maps, and their induced cohomology level morphisms define a class
\begin{equation} \label{eq:ext}
\big[[\Delta^{1|1|0}],[\Delta^{0|1|1}]\big] \in \mathit{Ext}_{H(\A) \otimes H(\A)^{\mathrm{opp}}}^1(H(\B^-),H(\B^+)).
\end{equation}
This is independent of all choices, and describes the $H(\A)$-bimodule $H(\B)$ as an extension of $H(\B^-)$ by $H(\B^+)$. Here is an example showing the difference (suggested by discussions with Maydanskiy). Let $V_1,V_2$ be two simple closed curves on a surface $M$, intersecting transversally in exactly two points, which are essential (can't be removed by a Hamiltonian isotopy). Then $\A$ is the algebra corresponding to the Kronecker quiver
\begin{equation}
\xymatrix{
\bullet \ar@/^{.5pc}/[rr] \ar@/_{.5pc}/[rr] && \bullet
}
\end{equation}
A simple computation shows that the $\mathit{Ext}$-group from \eqref{eq:ext} is $\K^3$. If the two curves are Hamiltonian isotopic, the bimodule extension is nontrivial, but otherwise it is trivial.
\end{remark}

\begin{remark} \label{th:factoring}
$\Delta$ can be factored into a sequence of smaller steps. Namely, assign to our Lagrangian submanifolds $V_i$ real numbers $o_i$, which are strictly decreasing and otherwise arbitrary. Fix some $c \in \R$ which satisfies $c \neq o_j - o_i$ for all $i,j$. Consider the sub-bimodule
\begin{equation}
{\F}^c = \bigoplus_{o_j - o_i < c} \mathit{CF}_\B(V_i,V_j) \subset \B.
\end{equation}
Define $\B^c$ to be the mapping cone $\Cone(\B^+ \rightarrow \B/\F^c)$, where the map is inclusion followed by projection. For $c \ll 0$, we have $\F^c = 0$ and $\B^c$ is quasi-isomorphic to $\B^-$. For $c \gg 0$, we have $\F^c = \B$ and $\B^c = \B^+$. In general,
\begin{align}
& \text{if $c>0$:} \quad \label{eq:bc-one}
H(\B^c) = \bigoplus_{i<j} \mathit{HF}(V_i,V_j) \oplus \bigoplus_i H^0(V_i;\K) \oplus
\bigoplus_{\substack{i > j \\ o_j - o_i > c}} \mathit{HF}(V_i,V_j), \\
& \text{if $c<0$:} \quad \label{eq:bc-two}
H(\B^c) =
\bigoplus_{\substack{i < j \\ o_j - o_i < c}} \mathit{HF}(V_i,V_j) \oplus \bigoplus_i \tilde{H}^*(V_i;\K) \oplus \bigoplus_{i>j} \mathit{HF}(V_i,V_j).
\end{align}
There are natural projection maps $\B^{c^-} \rightarrow \B^{c^+}$ for all $c^- \leq c^+$, which are compatible with composition, and which agree with $\Delta$ up to homotopy when $\pm c^\pm \gg 0$.
\end{remark}

\begin{remark}
Even though we are working with ungraded Floer groups throughout, it is worthwhile mentioning what would happen if one puts oneself in a situation where $\Z$-gradings are well-defined (or even if one just keeps track of $\Z/2$-gradings). It then makes sense to define $\B^-$ as $(\B/\B^+)[-1]$, where the upwards shift in the grading ensures that $\Delta$ has degree zero. The same holds for the definition of $\B^c$. As a consequence, in \eqref{eq:bc-one} the summands in the last term have their degrees shifted up by one. In \eqref{eq:bc-two}, the same applies to the last two terms.
\end{remark}

\subsection{A Hochschild cocycle\label{subsec:hh}}
There are connections between $\bar{\A}$ and $\B$ that go beyond the bimodule structure. This is clear on the level of Floer cohomology, where
\begin{equation} \label{eq:ab-inclusion}
H(\A) = R \oplus H(\bar\A) = R \oplus \bigoplus_{i<j} \mathit{HF}(V_i,V_j) = H(\B^+) \subset H(\B).
\end{equation}
To implement a cochain level version of this observation, we have to use an analytic framework in which the Hamiltonian term is being gradually turned on. In the simplest instance, our Riemann surface is the closed upper half-plane $S = \R \times \R^+$, with boundary condition given by some $V_j$. We consider this as a Riemann surface with one output, whose associated strip-like end is
\begin{equation} \label{eq:half-plane-end}
\epsilon(s,t) = \exp(-\pi (s+it-i)).
\end{equation}
Choose a family $J^0_j$ of almost complex structures parametrized by this surface, and which over the strip-like end reduces to $J^{0|1|0}_{\B,j,j}$. Take $H = H_\B$ to be the same Hamiltonian as before, except that we now make one more minor assumption:
\begin{equation} \label{eq:generic-3}
\parbox{28em}{
No critical point of $H$ lies on any $V_j$.
}
\end{equation}
Additionally choose $\gamma \in \Omega^1(S)$ which vanishes near the boundary, and whose pullback by the strip-like end is $g_\B(t)\, \mathit{dt}$. By counting solutions of the resulting equation \eqref{eq:gromov-2}, one gets numbers $n^0(x_0) \in \K$ for each $x_0 \in \phi_H(V_j) \cap V_j$, which together form a distinguished cocycle $\bar{\Phi}^0 \in \mathit{CF}_\B(V_j,V_j)$  representing $1 \in H^*(V_j;\K) \iso \mathit{HF}(V_j,V_j)$ (this is obvious if the isomorphism between Floer cohomology and classical cohomology is constructed as a Piunikhin-Salamon-Schwarz type map; see \cite{piunikhin-salamon-schwarz94} for the original Hamiltonian version, and \cite{albers,biran-cornea09c} for the Lagrangian counterpart relevant here).

\begin{remark} \label{th:constant-maps}
At this point, it is maybe appropriate to discuss some transversality issues which have only been mentioned implicitly so far. Using families of almost complex structures to achieve transversality, as we are doing, has one potential failure point. Namely, maps which satisfy
\begin{equation} \label{eq:zero-energy}
du = X_H \otimes \gamma
\end{equation}
remain solutions of \eqref{eq:gromov-2} for any choice of almost complex structure, hence cannot be perturbed away. If $\gamma$ is not closed, any $u$ as in \eqref{eq:zero-energy} is necessarily a constant map taking value at a stationary point of $X_H$. In the construction of $\bar\Phi^0$, this was precisely avoided by imposing \eqref{eq:generic-3} (and the same will apply to the generalizations below).
For closed $\gamma$, solutions of \eqref{eq:zero-energy} are less constrained, but still must lie on a single orbit of $X_H$. In the definition of $\bar\A$ where $H = 0$, these would occur whenever several of the $V_i$ intersect, but that was ruled out by \eqref{eq:generic-1}; except of course for the case $d = 1$, where they are not part of the moduli space, since they violate the stability (finite automorphisms) condition. In the definition of $\B$, the same role is played by \eqref{eq:b-general-position}.

To be more technically accurate, those conditions imply that \eqref{eq:zero-energy} must be violated at some point of $S$ which lies outside the strip-like ends. Locally near such a point, the choice of family of almost complex structures $(J_z)$ is essentially unconstrained, and that is sufficient to make standard transversality arguments go through. The exception to this is the definition of the differentials, where the translation-invariance requires additional arguments, whose structure we recall briefly. Following \cite{floer-hofer-salamon94,oh96c} (see also the discussion in Section \ref{subsec:bimodule-comparison} below), the main step is to establish a version of ``somewhere injectivity'', meaning that for any non-stationary solution $u: \R \times [0,1] \rightarrow M$, there is a point $(s_0,t_0)$ where $\partial_s u$ is nonzero and $u(s_0,t_0) \notin u(\R \setminus \{s_0\},t_0)$. Additionally, one can require that $u(s_0,t_0)$ is also distinct from $u_\zeta(t_0)$ for both limits $\zeta$. Roughly speaking, this means that a perturbation of the family of almost complex structures, near the parameter value $t_0$ and the point $u(s_0,t_0)$, affects the Cauchy-Riemann equation only locally near $(s_0,t_0)$.
\end{remark}

Before continuing, we need to recall some algebraic notions. For any $\bar{\A}$-bimodule $\PP$, the Hochschild complex of $\bar{\A}$ with coefficients in $\PP$ is
\begin{equation}
\mathit{CC}(\bar{\A},\PP) = \prod_{d \geq 0} \mathit{hom}_{R \otimes R}(\bar{\A}^{\otimes d},\PP).
\end{equation}
The $\mathit{hom}$'s are maps of $R$-bimodules; in particular, the $d = 0$ term is $\mathit{hom}_{R \otimes R}(R,\PP) = \bigoplus_k e_k\PP e_k$. The Hochschild differential is
\begin{equation} \label{eq:hochschild-d}
\begin{aligned}
& (\delta \bar{\Phi})^d(a_d,\dots,a_1) \\ & = \sum_{i,j} \mu_\PP^{d-i-j|1|i}(a_d,\dots,a_{i+j+1},\bar{\Phi}^j(a_{i+j},\dots,a_{i+1}),a_i,\dots,a_1) \\ & +
\sum_{i,j} \bar{\Phi}^{d-j+1}(a_d,\dots,a_{i+j+1},\mu_{\bar{\A}}^j(a_{i+j},\dots,a_{i+1}),a_i,\dots,a_1).
\end{aligned}
\end{equation}
It is a basic fact that Hoch\-schild cohomology has an alternative interpretation in terms of bimodule homomorphisms. Namely, let $\frakC$ be the differential category of bimodules over $\bar\A$. In particular, this contains the diagonal bimodule $\A = \bar{\A} \oplus R$ (the name comes from the fact that $\frakC$ can also be thought of as the category of strictly unital $\A$-bimodules). Then, there is a canonical chain map
\begin{equation} \label{eq:functorial-hh}
X: \mathit{CC}(\bar{\A},\PP) \longrightarrow \mathit{hom}_{\frakC}(\A,\PP).
\end{equation}
Explicitly, $\Phi = X(\bar{\Phi})$ is given by
\begin{multline} \label{eq:r-formula}
\Phi^{q|1|p}(a_{p+q+1},\dots,a_1) = \\ \sum_{i+j<p+1} \mu_\PP^{p+q-i-j+2|1|i}(a_{p+q+1},\dots,a_{p+1},\dots,\bar{\Phi}^j(a_{i+j},\dots,a_{i+1}),a_i,\dots,a_1).
\end{multline}
The restriction on the indices means that the element $a_{p+1} \in \A$ always appears to the left of $\bar{\Phi}$; and to plug that into $\mu_\PP$, we are implicitly using its extension to a strictly unital $\A$-bimodule. In the simplest case,
\begin{equation} \label{eq:phi1}
\Phi^{0|1|0}(a) = \mu_\PP^{1|1|0}(a,\bar{\Phi}^{0}).
\end{equation}

\begin{remark}
Concerning the asymmetry in \eqref{eq:r-formula}: the alternative formula with $a_{p+1}$ appearing on the right of $\bar{\Phi}$ also leads to a chain map, which turns out to be chain homotopic to \eqref{eq:functorial-hh}.
\end{remark}

\begin{lemma}
\eqref{eq:functorial-hh} is a quasi-isomorphism.
\end{lemma}

\proof This is again a well-known fact, so we only sketch the argument. One can filter both sides of \eqref{eq:functorial-hh} by length, meaning the number of inputs. Denoting these (complete decreasing) filtrations by $F^l$, one has $X(F^l \mathit{CC}(\bar{\A},\B)) \subset F^{l+1} \mathit{hom}_{\frakC}(\A,\B)$. Consider the induced maps between spectral sequences. At the $E^2$ level one gets two versions of the classical Hochschild cohomology of $H(\A)$ with coefficients in $H(\B)$ based of different projective resolutions of $H(\A)$, of the form $H(\A) \otimes T(H(\bar{\A})) \otimes H(\A)$ and $H(\A) \otimes T(H(\bar{\A})) \otimes H(\A) \otimes T(H(\bar{\A})) \otimes H(\A)$, respectively. \qed

We will use geometry to define a distinguished cocycle $\bar{\Phi} \in \mathit{CC}(\bar{\A},\B^+)$. Here, the choice made in defining $\B^+ \subset \B$ as in \eqref{eq:bplus} is to take $u_i = \bar\Phi^0$, which is also the leading order term of our Hochschild cocycle. Take a surface $S$ which is a $(d+1)$-punctured disc, $d \geq 0$, with an additional interior marked point. There is a unique identification
\begin{equation} \label{eq:punctured-plane}
S = (\R \times \R^+) \setminus \{\zeta_1,\dots,\zeta_d\}
\end{equation}
where $\zeta_1 < \cdots < \zeta_d$ are in $\R$, and the additional marked point becomes $i = (0,1) \in \R \times \R^+$. Equip this surface with boundary conditions $V_{i_0},\dots,V_{i_d}$ for some $i_0 < \cdots < i_d$, as well as with a family of almost complex structures which reduces to $J_{\B,i_0,i_d}^{0|1|0}$ on the end $|z| \gg 0$, and to $J_{\bar{\A},i_k,i_{k+1}}^1$ near the other $\zeta_k$. Using the same one-form $\gamma$ as before, one then has all the necessary ingredients for \eqref{eq:gromov-2}. As usual, we actually have to consider the universal family $\scrS^d \rightarrow \scrR^d$ of surfaces (of $(d+1)$-marked discs with an additional interior marked point), and choose almost complex structures $J^d_{i_0,\dots,i_d}$ on the total space of that family. Counting points in the resulting moduli spaces $\scrM^d(x_0,\dots,x_d)$ (where $x_0 \in \phi_H(V_{i_0}) \cap V_{i_d}$, and $x_k \in V_{i_{k-1}} \cap V_{i_k}$ for $k>0$) yields numbers $n^d(x_0,\dots,x_d)$, which assemble into the desired $\bar{\Phi}^d$. Figure \ref{fig:hochschild} shows the degenerations responsible for the two terms in the equation $\delta \bar\Phi = 0$.
\begin{figure}
\begin{center}
\begin{picture}(0,0)%
\includegraphics{hochschild.pstex}%
\end{picture}%
\setlength{\unitlength}{3947sp}%
\begingroup\makeatletter\ifx\SetFigFont\undefined%
\gdef\SetFigFont#1#2#3#4#5{%
  \reset@font\fontsize{#1}{#2pt}%
  \fontfamily{#3}\fontseries{#4}\fontshape{#5}%
  \selectfont}%
\fi\endgroup%
\begin{picture}(4371,2444)(-1661,-2358)
\put(886,-1186){\makebox(0,0)[lb]{\smash{{\SetFigFont{10}{12.0}{\rmdefault}{\mddefault}{\updefault}{\color[rgb]{0,0,0}$a_{i+j}$}%
}}}}
\put(1951,-1636){\makebox(0,0)[lb]{\smash{{\SetFigFont{10}{12.0}{\rmdefault}{\mddefault}{\updefault}{\color[rgb]{0,0,0}$a_{i+1}$}%
}}}}
\put(-1249,-1936){\makebox(0,0)[lb]{\smash{{\SetFigFont{10}{12.0}{\rmdefault}{\mddefault}{\updefault}{\color[rgb]{0,0,0}$a_1$}%
}}}}
\put(-1274,-661){\makebox(0,0)[lb]{\smash{{\SetFigFont{10}{12.0}{\rmdefault}{\mddefault}{\updefault}{\color[rgb]{0,0,0}$a_d$}%
}}}}
\put(2251,-1036){\makebox(0,0)[lb]{\smash{{\SetFigFont{10}{12.0}{\rmdefault}{\mddefault}{\updefault}{\color[rgb]{0,0,0}$a_{i+1}$}%
}}}}
\put(2476,-61){\makebox(0,0)[lb]{\smash{{\SetFigFont{10}{12.0}{\rmdefault}{\mddefault}{\updefault}{\color[rgb]{0,0,0}$a_{i+j}$}%
}}}}
\end{picture}%
\caption{\label{fig:hochschild}}
\end{center}
\end{figure}

\begin{remark} \label{th:moduli1}
To clarify, the moduli spaces $\scrR^d$ have two equivalent descriptions, both of which have advantages. One is as configuration spaces of ordered $d$-tuples of points $(\zeta_1,\dots,\zeta_d)$ on the real line, as in \eqref{eq:punctured-plane} (and $\scrR^0$ is a point). In that case one does not divide by any symmetries, which is convenient when one thinks of equipping $(\R \times \R^+) \setminus \{\zeta_1,\dots,\zeta_d\}$ with an inhomogeneous term. From this point of view, the correct compactification $\bar{\scrR}^d$ is constructed by starting with \cite{fulton-macpherson} and then adding more boundary strata which account for the possibility of some marked points going to $\pm\infty$. The second way to describe $\scrR^d$ is as moduli space of $(d+1)$-pointed discs with one interior marked point $z$, divided by the group of automorphisms $\mathit{Aut}(\bar{S}) \iso \mathit{PSL}(2,\R)$. The advantage of this viewpoint is that the construction of $\bar\scrR^d$ is standard, being a special case of the (real) Deligne-Mumford compactification. The transition between the two points of view, as already explained above, occurs by using $z$ and the point $\zeta_0 \in \partial\bar{S}$ to specify an identification between $\bar{S}$ and a fixed closed disc, hence between $S \setminus \{\zeta_0\}$ and the upper half plane.
\end{remark}

As usual, we retrospectively adjust the notation:
\begin{center}
\begin{tabular}{|l|l|}
\hline
$\gamma_{\bar{\Phi}}$ & one-form on the upper half plane \\
$\scrS^d_{\bar{\Phi}} \rightarrow \scrR^d_{\bar{\Phi}}$ & moduli space of $d+1$-punctured discs \\
& \hspace{1em} with an additional interior marked point \\
$J^d_{\bar{\Phi},i_0,\dots,i_{p+q+1}}$ & family of almost complex structures \\
$\scrM_{\bar{\Phi}}^d(x_0,\dots,x_d)$ & moduli spaces of pseudo-holomorphic maps \\
$n_{\bar{\Phi}}^d(x_0,\dots,x_d)$ & numbers counting points in the spaces above \\
\hline
\end{tabular}
\end{center}

From \eqref{eq:phi1} one sees that the resulting bimodule homomorphism $\Phi = X(\bar\Phi)$ implements the isomorphism \eqref{eq:ab-inclusion}, which shows that:

\begin{prop} \label{th:plus-quasi-iso}
$\Phi: \A \rightarrow \B^+$ is an $\bar\A$-bimodule quasi-isomorphism.
\end{prop}

\section{Poincar{\'e} duality\label{sec:fibre2}}

\subsection{The dual diagonal bimodule\label{subsec:dual-diagonal}}
The material in this section, while a natural continuation of the previous discussion, is not as important for the main results of the paper, and we'll therefore cover it in less detail. Assume that the $V_i$ are $\K$-homology spheres (see Remark \ref{th:generalized-directed} for further discussion of this assumption). Then, by Poincar{\'e} duality in Floer cohomology, $H(\B^-)$ is canonically dual to $H(\B^+) \iso H(\A)$.

To formulate the corresponding chain level results, we need a little more algebraic generalities. Given any finite-dimensional $\bar{\A}$-bimodule $\PP$, consider the dual vector space $\PP^\vee = \mathit{hom}_{\K}(\PP,\K)$ with the induced $\bar{\A}$-bimodule structure. This can be used as coefficient bimodule for the Hochschild complex, which turns into
\begin{equation} \label{eq:cc-dual}
\mathit{CC}(\bar{\A},\PP^\vee) = \prod_d \mathit{hom}_{R \otimes R}(\bar{\A}^{\otimes d} \otimes \PP,R).
\end{equation}
The cohomology of this is better known as the dual of Hochschild homology. Let $\A^\vee$ be the dual diagonal bimodule. Then there is a canonical quasi-isomorphism
\begin{equation} \label{eq:dual-iso}
Y: \mathit{CC}(\bar{\A},\PP^\vee) \longrightarrow \mathit{hom}_{\frakC}(\PP,\A^\vee).
\end{equation}
In terms of the right hand side of \eqref{eq:cc-dual} (whose components are written as $\bar{\Psi}^{d|1}$ for $d \geq 1$), $\Psi = Y(\bar{\Psi})$ is given by
\begin{multline}
\qquad \leftsc \Psi^{q|1|p}(a_{p+q+1},\dots,b,\dots,a_1),a' \rightsc =
\sum_{i, j} \bar{\Psi}^{i|1}\big(a_{p+j+i},\dots,a_{p+j+1}, \\ \mu_{\PP}^{j|1|p-i}(a_{p+j},\dots,a_{p+1},b,a_p,\dots,a_1,a',a_{p+q+1},\dots,
a_{p+j+i+1})\big). \qquad
\end{multline}

Now return to the specific bimodule $\B$ from \eqref{eq:b}. In parallel with our previous discussion, geometry yields a distinguished cocycle $\bar{\Psi} \in \mathit{CC}(\bar{\A},\B^\vee)$. This uses the same Riemann surfaces as in \eqref{eq:punctured-plane}, but where the end $|z| \rightarrow \infty$ is now considered as an input. The leading order term $\bar{\Psi}^0: \mathit{CF}_\B(V_j,V_j) \rightarrow \K$ is a cochain level realization of the integration map $\mathit{HF}(V_j,V_j) \iso H^*(V_j;\K) \rightarrow \K$. If we restrict $\bar{\Psi}$ to the sub-bimodule $\B^+$, then (because of directedness and the $R$-bimodule structure) its only possible nonvanishing term is
\begin{equation}
\bar{\Psi}^0|\B^+: \bigoplus_j e_j \B^+ e_j = \bigoplus_j R u_j \longrightarrow \K.
\end{equation}
By construction $u_j$ represents the identity in $\mathit{HF}(V_j,V_j)$, hence its image under $\bar\Psi^0$ must be zero by the previous observation. This implies that $\bar\Psi$ actually lies in the subspace $\mathit{CC}(\bar{\A},(\B^-)^\vee)$. Moreover, just by looking at $\bar\Psi^0$, it follows that $\Psi = Y(\bar\Psi)$ implements the isomorphism $H(\B^-) \iso H(\A)^\vee$. Hence:

\begin{prop} \label{th:minus-quasi-iso}
$\Psi: \B^- \rightarrow \A^\vee$ is an $\bar\A$-bimodule quasi-isomorphism.
\end{prop}

Hence, using $\Psi$ and $\Phi$, we can rewrite the short exact sequence \eqref{eq:b-plus-minus} as an exact triangle in $H(\frakC)$ of the form
\begin{equation} \label{eq:dual-triangle}
\xymatrix{
\A \ar[r] & \B \ar[r] & \A^\vee \ar@/^1pc/[ll]
}
\end{equation}
which in particular means that its boundary map \eqref{eq:boundary} turns into an element
\begin{equation} \label{eq:boundary-2}
[\Phi^{-1} \circ \Delta \circ \Psi^{-1}] \in H(\mathit{hom}_{\frakC}(\A^\vee,\A)),
\end{equation}
where $\Phi^{-1}$, $\Psi^{-1}$ are quasi-inverses (unique up to homotopy).

\begin{remark} \label{th:generalized-directed}
It is instructive to think of what the appropriate generalization for Lagrangian submanifolds $(V_1,\dots,V_m)$ with different topology might be. On the level of chain complexes, we leave $\B$ as before \eqref{eq:b}. However, our original definition of $\A$ is clearly not appropriate any more, since it is no longer ``half of $\B$''. Instead, suppose that for each $i$ we are given an abstract $(n+1)$-dimensional compact manifold $L_i$ with $\partial L_i = V_i$. Choose Morse functions $h_{V_i}$ on $V_i$ and $h_{L_i}$ on $L_i$, in such a way that there is a surjective map of the associated Morse complexes
\begin{equation} \label{eq:morse-projection}
C^*(h_{L_i}) \longrightarrow C^*(h_{V_i})
\end{equation}
realizing the restriction map on cohomology (this can be achieved by taking $h_{L_i}(r,y) = \psi_i(r) + h_{V_i}(y)$ in a collar neighbourhood $(r,y) \in (-\epsilon,0] \times V_i \subset L_i$, where the radial component $\psi_i$ has a single local minimum). Now choose additional machinery (families of Morse functions parametrized by metrized ribbon graphs) which makes both Morse complexes into $A_\infty$-algebras, and such that \eqref{eq:morse-projection} is an $A_\infty$-homomorphism. We then re-define
\begin{equation} \label{eq:generalized-a}
\A = \bigoplus_i C^*(h_{L_i}) \oplus \bar{\A}
\end{equation}
and equip that with the structure of a cohomologically unital $A_\infty$-algebra (built from the $A_\infty$-structures of each summand, together with the map \eqref{eq:morse-projection} and auxiliary moduli spaces mixing Morse flow trees and holomorphic discs). Similarly, $\B$ can be made into a cohomologically unital $\A$-bimodule. This comes with an $\A$-bimodule homomorphism $\Phi: \A \rightarrow \B$ whose cohomology level effect is the direct sum of the identity maps on $\mathit{HF}(V_i,V_j)$, $i<j$, and the restriction maps $H^*(L_i;\K) \rightarrow H^*(V_i;\K)$. The cohomology of its mapping cone is
\begin{equation}
\bigoplus_i H^*(L_i,V_i;\K) \oplus \bigoplus_{i>j} \mathit{HF}(V_i,V_j),
\end{equation}
which is dual to $H(\A)$. It is therefore reasonable to expect a generalization of the picture from \eqref{eq:dual-triangle}. However, setting up the relevant argument requires additional work, which we have not carried out. Such generalizations may be useful in the context of Landau-Ginzburg models with non-isolated singularities (this was pointed out to the author by Abouzaid and Auroux).

We'll underline two points about this generalization. First, even though we are initially considering $n$-dimensional manifolds $V_i$, $(n+1)$-dimensional bounding manifolds $L_i$ play a natural role. Secondly, if we review our original construction from this viewpoint, it appears to be based on a technically convenient but potentially misleading shortcut, which was to assume that the $L_i$ are $\K$-homology balls, and then to interpret $\A$ as the result of the purely algebraic process of adding units to $\bar{\A}$. No such strategy can work in general, since \eqref{eq:generalized-a} contains additional information not available in $\bar{\A}$.
\end{remark}

\subsection{Full Fukaya categories\label{subsec:strict}}
There is an alternative approach, which will play no further role in this paper, but deserves to be mentioned because it leads to slightly stronger results, and also relates to material that appears elsewhere in the literature.

Fix $(V_1,\dots,V_m)$ as before. The associated Fukaya $A_\infty$-structure, in the non-directed sense of the word, is an $A_\infty$-algebra $\B$ over $R$ with $H(\B) = \bigoplus_{i,j} \mathit{HF}(V_i,V_j)$. Several definitions are available, which are all equivalent up to quasi-isomorphism. For instance, \cite{seidel04} uses perturbations of the Cauchy-Riemann equation similar to \eqref{eq:gromov-2}, but where now a different Hamiltonian $H_{ij}$ is chosen for each pair $(V_i,V_j)$, and correspondingly the inhomogeneous terms vary over the punctured discs in a more complicated way. Other approaches use singular cochains \cite{fooo}, or Morse-Bott techniques similar to the one mentioned in Remark \ref{th:generalized-directed} \cite{pozniak, bourgeois02, biran-cornea09, biran-cornea09c, johns08, seidel08b, sheridan11}. In each of these approaches, one can ensure (based on our exactness assumption) that the unit in $e_j\B e_j$ is represented by a generator $u_j$ satisfying
\begin{equation} \label{eq:uuu}
\mu^d_\B(u_j,\dots,u_j) = \begin{cases} u_j & d = 2, \\ 0 & \text{otherwise.} \end{cases}
\end{equation}
As a consequence, $\A = \bigoplus_j R u_j \oplus \bigoplus_{i<j} e_j\B e_i$ is an $A_\infty$-subalgebra.

\begin{remark}
This $\A$ is quasi-isomorphic to the $A_\infty$-algebra obtained by starting from $\bar{\A} = \bigoplus_{i<j} e_j\B e_i$ and adjoining strict units as in \eqref{eq:adjoin-units}. The proof works as follows. First (for purely algebraic reasons) one can find a strictly unital $\A'$ and a quasi-isomorphism
\begin{equation} \label{eq:unity-trick}
\A \longrightarrow \A'.
\end{equation}
Restrict that to $\bar{\A} \rightarrow \A'$, and then extend it uniquely to a strictly unital $A_\infty$-homomorphism $R \oplus \bar{\A} \rightarrow \A'$, which is necessarily a quasi-isomorphism. Finally, compose that with the inverse of \eqref{eq:unity-trick}.
\end{remark}

We then have a short exact sequence of $\A$-bimodules
\begin{equation} \label{eq:quotient-sequence}
0 \longrightarrow \A \longrightarrow \B \longrightarrow \B/\A \longrightarrow 0.
\end{equation}
Take $\Delta$ to be the associated boundary map. The $A_\infty$-algebra (rather than bimodule) structure on $\B$ yields additional information about $\Delta$. Namely, it is {\em ambidextrous}, which means that the following diagram of $\A$-bimodules commutes up to homotopy \cite[Lemma 3.2]{seidel08}:
\begin{equation}
\xymatrix{
 & (\B/\A) \otimes_{\A} (\B/\A)
 \ar[dl]_{\mathrm{id} \otimes \Delta} \ar[dr]^{\Delta \otimes \mathrm{id}} & \\
 (\B/\A) \otimes_{\A} \A
 \ar[dr]_{\htp} &&
 \A \otimes_{\A} (\B/\A) \ar[dl]^{\htp} \\
 & \B/\A.  &
}
\end{equation}
One can prove that, with the definitions just given, $\B/\A$ is quasi-isomorphic to $\A^\vee$ (see \cite[Proposition 5.1]{seidel09} for an outline of the argument), so that the boundary map takes on the same form as in \eqref{eq:boundary-2}. Indeed, one can prove that the two constructions yield the same result up to chain homotopy.

\subsection{Boundary $A_\infty$-algebras\label{sec:boundary-alg}}
There is yet another and more ambitious approach, which builds Poincar{\'e} duality into the chain level structure itself. At least in the form explained here, this only applies over a coefficient field with $\mathrm{char}(\K) = 0$, which in turn requires at least $\Z/2$-gradings; and to simplify one technical point later on, we will actually use $\Z$-gradings. This unfortunately means that the conventions here differ a bit from those elsewhere in the paper.

Following the indication in Remark \ref{th:generalized-directed}, we take the cohomology of manifolds with boundary as a guide. Define an $(n+1)$-dimensional {\em algebra with boundary} to be a finite-dimensional graded unital associative algebra $A$, together with a distinguished element
\begin{equation}
D = \sum_j D^{2,j} \otimes D^{1,j} \in A \otimes A
\end{equation}
of degree $n+1$, satisfying the following two conditions. First, it should be graded symmetric up to a factor of $(-1)^{n+1}$, so
\begin{equation} \label{eq:gr-sym}
\sum_j (-1)^{\mathrm{deg}(D^{j,2})\mathrm{deg}(D^{j,1})} D^{j,1} \otimes D^{j,2} = (-1)^{n+1}
 \sum_j D^{j,2} \otimes D^{j,1}.
\end{equation}
Secondly, for all $a \in A$ we want to have
\begin{equation} \label{eq:both-sides}
\begin{aligned}
& \sum_j a D^{j,2} \otimes D^{j,1} = (-1)^{(n+1) \mathrm{deg}(a)} \sum_j D^{j,2} \otimes D^{j,1} a, \\
& \sum_j (-1)^{\mathrm{deg}(a) \mathrm{deg}(D^{j,2})} D^{j,2} a \otimes D^{j,1} = \sum_j (-1)^{\mathrm{deg}(a) \mathrm{deg}(D^{j,2})} D^{j,2} \otimes a D^{j,1}
\end{aligned}
\end{equation}
(in view of \eqref{eq:gr-sym}, each of the two parts of \eqref{eq:both-sides} implies the other). One interpretation of \eqref{eq:both-sides} is that the two ways of combining the product and $D$ to obtain a coproduct coincide.

As for any finite-dimensional algebra, $A \oplus A^\vee[-n]$ (the bracket means shifting the grading up by $n$) is canonically an $n$-dimensional Frobenius algebra. The product is obtained from the bimodule structure of $A^\vee$ (what's known as a trivial extension algebra). The Frobenius map $\int: A \oplus A^\vee[-n] \rightarrow \K$ is $(a,a^\vee) \mapsto a^\vee(e)$, where $e$ is the unit. The element $D$ gives rise to a compatible differential, namely
\begin{equation}
d(a,a^\vee) = \big( (-1)^{(n+1)\mathrm{deg}(a^\vee)} \sum_j a^\vee(D^{j,2}) D^{j,1}, 0 \big).
\end{equation}
We call the resulting Frobenius dga $(A \oplus A^\vee[-n], d, \int)$ the boundary of $A$, and denote it by $\partial A$.

\begin{example} \label{th:cohomology-boundary}
Take $A = H^*(L;\K)$, where $L$ is an $(n+1)$-dimensional compact oriented manifold with boundary. We can make this into an algebra with boundary, by taking $D$ to be the Poincar{\'e} dual of the diagonal (or equivalently, the dual of the intersection pairing on homology). From the long exact sequence
\begin{equation}
\cdots \rightarrow H^{n+1-*}(L;\K)^\vee \longrightarrow H^*(L;\K) \longrightarrow H^*(\partial L;\K) \rightarrow \cdots
\end{equation}
we see that $H(\partial A) = H^*(\partial L;\K)$ reproduces the cohomology of the boundary as a graded vector space. However, the example $L = [-1,1]$ shows that the induced product structure on $H(\partial A)$ is not necessarily the correct one. This should be blamed on working on the cohomological level, and ignoring higher order operations.
\end{example}

The subspace $A \subset \partial A$ is a dg subalgebra (with vanishing differential), as well as an isotropic subspace for the pairing induced by $\int$. When including higher order operations into our framework, we reverse directions and take that observation as the starting point:

\begin{definition}
Take a finite-dimensional graded vector space $\A$, together with the structure of an $A_\infty$-algebra on $\A \oplus \A^\vee[-n]$, cyclic with respect to the standard pairing. Assume that $\A \subset \A \oplus \A^\vee[-n]$ itself is an $A_\infty$-subalgebra, and that $H(\A)$ contains an element which is a unit for the entire algebra $H(\A \oplus \A^\vee[-n])$. We then call $\A$ an $(n+1)$-dimensional {\em $A_\infty$-algebra with boundary}, and $\partial \A = \A \oplus \A^\vee[-n]$ its boundary.
\end{definition}

In terms of noncommutative geometry \cite[Section 10]{kontsevich-soibelman06}, the data describe a noncommutative formal symplectic vector space, equipped with a cyclic function $S$ (with a critical point at the origin) satisfying $\{S,S\} = 0$, as well as a Lagrangian subspace on which $S$ vanishes identically (this follows a general philosophical guideline of having a symplectic vector space associated to a manifold, and a Lagrangian subspace to a filling of that manifold; compare for instance \cite{salamon96}). On a more pedestrian level, one can give a description of the structure in terms of operations on $\A$ itself. Namely, take any cyclically ordered finite set $I$, partitioned into subsets $I_{\mathrm{in}}$, $I_{\mathrm{out}}$, and subject to the conditions that
\begin{equation} \label{eq:stability-i}
I_{\mathrm{out}} \neq \emptyset, \quad |I_{\mathrm{in}}| + 2|I_{\mathrm{out}}| \geq 3.
\end{equation}
For each such $I$ we then have a multilinear map
\begin{equation} \label{eq:mu-i}
\mu^I: \A^{\otimes I_{\mathrm{in}}} \longrightarrow \A^{\otimes I_{\mathrm{out}}}[2+(n-1)(|I_{\mathrm{out}}|-1)-|I_{\mathrm{in}}|].
\end{equation}
Up to appropriate sign, these operations are preserved under isomorphism of cyclically ordered sets $I$, and moreover they satisfy quadratic relations, which for $|I_{\mathrm{out}}| = 1$ reduce to the ordinary $A_\infty$-relations. Let's look at the first level beyond that. The operations with $|I_{\mathrm{out}}| = 2$ give rise to an $\A$-bimodule homomorphism $\Delta: \A^\vee \rightarrow \A$ of degree $n+1$. On the cohomology $A = H(\A)$, we have a product as well as a class $D \in (A \otimes A)^{n+1}$ (the latter inherited from the operation with $|I_{\mathrm{in}}| = 0$, $|I_{\mathrm{out}}| = 2$), and this makes $A$ into an algebra with boundary in the sense introduced previously.

To apply this to Floer theory, we assume that our symplectic manifold $M$ comes with a trivialization of its canonical bundle, and that the homology spheres $V_i$ are {\em Spin} and have zero Maslov index, hence can be equipped with gradings. Take the coefficient field $\K = \R$. Techniques from \cite{fukaya09} (see also \cite[Section 10]{kontsevich-soibelman06} for the general algebraic framework) allow us to define a cyclic (and graded) $A_\infty$-algebra over $R = \K^m$ whose cohomology is $\bigoplus_{i,j} \mathit{HF}(V_i,V_j)$. In fact, we can assume that this algebra is minimal, meaning that the differential vanishes. Denote it by $\B$, and let $R \subset \B$ be the subspace generated by the identity elements in each $e_i\B e_i \iso \mathit{HF}(V_i,V_i)$. Because of the grading, we can use the same trick as in \eqref{eq:uuu} to show that
\begin{equation}
\A = R \oplus \bigoplus_{i<j} \mathit{HF}(V_i,V_j) \subset \B
\end{equation}
is an $A_\infty$-subalgebra. Moreover, it is a maximal isotropic subspace, with respect to the pairing which is part of the cyclic $A_\infty$-structure of $\B$. One can use the pairing to identify $\B = \A \oplus \A^\vee[-n]$ as vector spaces. Hence, $\A$ is an $A_\infty$-algebra with boundary, of dimension $n+1$.


\section{Geometry of the base\label{sec:analysis}}

This section introduces a toy model for our main theory, where the target space is the complex plane (this will later be thought of as the base of a Lefschetz fibration). We use the simplified situation to introduce basic features as well as some useful technical tricks.

\subsection{Geometric data}
Take the complex plane $\C = \R^2$, with coordinates $w = p+i q$. We equip it with its standard complex structure and symplectic form $\omega_\C = \frac{i}{2} \mathit{dw} \wedge \mathit{d\bar{w}} = \mathit{dp} \wedge \mathit{dq}$. We will occasionally use a primitive $\theta_\C$ of this form (the specific choice of primitive is irrelevant).

Informally speaking, we consider $\C$ as having a base point at infinity in positive real direction. In connection with this, subsets bounded by vertical lines will play a prominent role. For any $a<b$ in $\R$, define
\begin{equation} \label{eq:uw-notation}
\begin{aligned}
& U_{a,b} = \{w = p+iq \in \C \;:\; a \leq p \leq b\}, \\
& W_a = \{w = p+iq \in \C\;:\; a \leq p\}.
\end{aligned}
\end{equation}

\begin{definition} \label{def:vanishing-1}
A {\em vanishing path} $l \subset \C$ is a properly embedded semi-infinite path $l \iso \R^+$, with endpoint on $\C \setminus W_{-2}$, such that $l \cap W_{-2}$ is the graph of a function $q(p)$, which is constant for $p \geq 1$. The constant $o = q(p)$, $p \geq 1$, is called the {\em ordinate} of $l$.

A {\em simply intersecting basis} of vanishing paths is an ordered collection $(l_1,\dots,l_m)$ with the following properties. First of all, the $l_i$ do not intersect each other in $\C \setminus W_{-2}$. Write $l_i \cap W_{-2} = \{q = q_i(p)\}$. Then, the remaining conditions say that for any $i<j$, the difference $q_i(p) - q_j(p)$ is nondecreasing for all $p \geq -2$, is negative at $p = 0$ and positive at $p = 1$, and has positive derivative for $p \in (0,1)$.
\end{definition}

This of course implies that the ordinates $o_i$ of paths $l_i$ in a simply intersecting basis are strictly decreasing. Moreover, any two different such paths $l_i$ and $l_j$ intersect exactly once, and that intersection point is transverse and lies in the interior of $U_{0,1}$.

\begin{definition} \label{def:displacement-1}
Fix $c \in \R$. A {\em $c$-displacement function} is an $h \in \smooth(\C,\R)$ with the following properties. $h$ depends only on $p = \mathrm{re}(w)$. Moreover, $h(p) = 0$ for $p \leq -1$, whereas $h'(p) = c$ is constant for $p \geq 2$.
\end{definition}

We denote by $X_h = h'(p)\,\partial_q$ the associated Hamiltonian vector field, and by $\phi_h(p+iq) = p+iq+ih'(p)$ its time-one map. Take vanishing paths $l_0,l_1$, with ordinates $o_0,o_1$. If we choose $c \in \R$ satisfying
\begin{equation} \label{eq:nondeg}
c \neq o_1 - o_0,
\end{equation}
then $\phi_h(l_0) \cap l_1$ is compact. Different signs of $c - o_1 + o_0$ correspond to the two possible ways of dealing with the intersection at infinity between our paths. Figure \ref{fig:strips} summarizes some of the conditions from Definitions \ref{def:vanishing-1} and \ref{def:displacement-1} (plus one more notion that will be introduced a little later; also, note that for space reasons the figure is rotated by $\pi/2$).
\begin{figure}
\begin{center}
\begin{picture}(0,0)%
\includegraphics{strips-new.pstex}%
\end{picture}%
\setlength{\unitlength}{3947sp}%
\begingroup\makeatletter\ifx\SetFigFont\undefined%
\gdef\SetFigFont#1#2#3#4#5{%
  \reset@font\fontsize{#1}{#2pt}%
  \fontfamily{#3}\fontseries{#4}\fontshape{#5}%
  \selectfont}%
\fi\endgroup%
\begin{picture}(4617,6551)(4489,4635)
\put(6041,7364){\makebox(0,0)[lb]{\smash{{\SetFigFont{10}{12.0}{\rmdefault}{\mddefault}{\updefault}{\color[rgb]{0,0,0}$l_2$}%
}}}}
\put(4651,6014){\makebox(0,0)[lb]{\smash{{\SetFigFont{10}{12.0}{\rmdefault}{\mddefault}{\updefault}{\color[rgb]{0,0,0}vector field $X_h$}%
}}}}
\put(7501,6614){\makebox(0,0)[lb]{\smash{{\SetFigFont{10}{12.0}{\rmdefault}{\mddefault}{\updefault}{\color[rgb]{0,0,0}$p=-1$}%
}}}}
\put(4651,4739){\makebox(0,0)[lb]{\smash{{\SetFigFont{10}{12.0}{\rmdefault}{\mddefault}{\updefault}{\color[rgb]{0,0,0}endpoints of paths}%
}}}}
\put(7501,5414){\makebox(0,0)[lb]{\smash{{\SetFigFont{10}{12.0}{\rmdefault}{\mddefault}{\updefault}{\color[rgb]{0,0,0}$p=-2$}%
}}}}
\put(4651,5864){\makebox(0,0)[lb]{\smash{{\SetFigFont{10}{12.0}{\rmdefault}{\mddefault}{\updefault}{\color[rgb]{0,0,0}is zero}%
}}}}
\put(7501,7814){\makebox(0,0)[lb]{\smash{{\SetFigFont{10}{12.0}{\rmdefault}{\mddefault}{\updefault}{\color[rgb]{0,0,0}$p = 0$}%
}}}}
\put(7501,9014){\makebox(0,0)[lb]{\smash{{\SetFigFont{10}{12.0}{\rmdefault}{\mddefault}{\updefault}{\color[rgb]{0,0,0}$p = 1$}%
}}}}
\put(7501,10214){\makebox(0,0)[lb]{\smash{{\SetFigFont{10}{12.0}{\rmdefault}{\mddefault}{\updefault}{\color[rgb]{0,0,0}$p = 2$}%
}}}}
\put(5176,8339){\makebox(0,0)[lb]{\smash{{\SetFigFont{10}{12.0}{\rmdefault}{\mddefault}{\updefault}{\color[rgb]{0,0,0}$X_h$}%
}}}}
\put(4651,10889){\makebox(0,0)[lb]{\smash{{\SetFigFont{10}{12.0}{\rmdefault}{\mddefault}{\updefault}{\color[rgb]{0,0,0}times $\partial_q$}%
}}}}
\put(4651,11039){\makebox(0,0)[lb]{\smash{{\SetFigFont{10}{12.0}{\rmdefault}{\mddefault}{\updefault}{\color[rgb]{0,0,0}$X_h$ is a constant}%
}}}}
\put(6601,7514){\makebox(0,0)[lb]{\smash{{\SetFigFont{10}{12.0}{\rmdefault}{\mddefault}{\updefault}{\color[rgb]{0,0,0}$l_1$ }%
}}}}
\put(7276,6014){\makebox(0,0)[lb]{\smash{{\SetFigFont{10}{12.0}{\rmdefault}{\mddefault}{\updefault}{\color[rgb]{0,0,0}inhomogeneous}%
}}}}
\put(7276,5864){\makebox(0,0)[lb]{\smash{{\SetFigFont{10}{12.0}{\rmdefault}{\mddefault}{\updefault}{\color[rgb]{0,0,0}term $X_k$ is zero}%
}}}}
\put(7276,8414){\makebox(0,0)[lb]{\smash{{\SetFigFont{10}{12.0}{\rmdefault}{\mddefault}{\updefault}{\color[rgb]{0,0,0}intersection points of paths}%
}}}}
\put(7276,8264){\makebox(0,0)[lb]{\smash{{\SetFigFont{10}{12.0}{\rmdefault}{\mddefault}{\updefault}{\color[rgb]{0,0,0}in a simply intersecting basis}%
}}}}
\put(7276,9689){\makebox(0,0)[lb]{\smash{{\SetFigFont{10}{12.0}{\rmdefault}{\mddefault}{\updefault}{\color[rgb]{0,0,0}vanishing paths are}%
}}}}
\put(7276,9539){\makebox(0,0)[lb]{\smash{{\SetFigFont{10}{12.0}{\rmdefault}{\mddefault}{\updefault}{\color[rgb]{0,0,0}horizontal}%
}}}}
\put(7276,10589){\makebox(0,0)[lb]{\smash{{\SetFigFont{10}{12.0}{\rmdefault}{\mddefault}{\updefault}{\color[rgb]{0,0,0}$X_k = \gamma \otimes \partial_q$}%
}}}}
\end{picture}%
\caption{\label{fig:strips}}
\end{center}
\end{figure}%

Our next task is to define the appropriate class of Riemann surfaces. The difference, with respect to the standard setup in Section \ref{sec:fibre}, is that the one-forms $\gamma$ now take on a more fundamental role.

\begin{definition} \label{def:riemann-2}
A {\em worldsheet} $(S,\gamma,l)$ is a triple of the following form. $S = \bar{S} \setminus \Sigma$ is a Riemann surface as in Definition \ref{def:riemann-1}. $\gamma \in \Omega^1(S)$ is a real one-form with $\gamma|\partial S = 0 \in \Omega^1(\partial S)$. Moreover, it should satisfy
 \begin{equation} \label{eq:sub-closed}
d\gamma \leq 0,
\end{equation}
where the sign is with respect to the complex orientation of $S$, and $d\gamma = 0$ outside a compact subset. Consider a puncture $\zeta \in \Sigma$. For concreteness, take a local holomorphic coordinate in which $\bar{S}$ is the upper half plane and $\zeta$ the origin. Because of the previous assumptions,
\begin{equation} \label{eq:conv}
\mathrm{res}_\zeta(\gamma) \stackrel{\mathrm{def}}{=} \textstyle \int_{|z| = \epsilon} \gamma
\end{equation}
is independent of the choice of $\epsilon$, provided that is sufficiently small, and of the local coordinate. The boundary conditions $l = (l_z)$ are a locally constant family of vanishing paths parametrized by points $z \in \partial S$. In local coordinates around some $\zeta$ as before, we have two vanishing paths associated to points on the positive and negative real half-axes. Denoting their ordinates by $o_{\zeta,\pm}$, we additionally require that
\begin{equation} \label{eq:rzeta}
o_{\zeta,+} - o_{\zeta,-} + \mathrm{res}_\zeta(\gamma) \neq 0.
\end{equation}
\end{definition}

\begin{remark}
A consequence of the condition \eqref{eq:sub-closed} is that punctures are divided into two classes depending on the sign of the left-hand side of \eqref{eq:rzeta}. For instance, if $S$ is the upper half plane with one end, then necessarily $\mathrm{res}_\zeta(\gamma) > 0$. Similar asymmetries appear in other theories with non-compact target spaces, such as wrapped Floer cohomology \cite{abouzaid-seidel07} or string topology.
\end{remark}

As in our previous discussion of Fukaya $A_\infty$-structures, the definition above is the conceptually correct one, but it is technically convenient to add strip-like ends, and to restrict the behaviour of $\gamma$ over those ends. Namely, we want
\begin{equation} \label{eq:gamma-epsilon}
\epsilon_\zeta^* \gamma = g_\zeta(t) \mathit{dt},
\end{equation}
where $g_\zeta$ is some function. Denote by $l_{\zeta,k}$ the vanishing paths associated to boundary points $\epsilon_\zeta(\cdot,k)$, and by $o_{\zeta,k}$ their ordinates. In view of \eqref{eq:conv}, $\int_0^1 g_\zeta(t) \mathit{dt} = \pm \mathrm{res}_\zeta(\gamma)$, where the sign is $-$ for $\zeta \in \Sigma_{\mathrm{in}}$, and $+$ for $\zeta \in \Sigma_{\mathrm{out}}$. We write $c_\zeta$ for this integral. Then, \eqref{eq:rzeta} reduces to saying that \eqref{eq:nondeg} holds for every $\zeta$:
\begin{equation} \label{eq:rzeta-2}
c_\zeta \neq o_{\zeta,1} - o_{\zeta,0}.
\end{equation}
We call $(S,\gamma,l,\epsilon)$ with these properties a {\em worldsheet with strip-like ends}.

\subsection{The perturbed Cauchy-Riemann equation\label{subsec:perturbed-1}}
Let $(S,\gamma,l,\epsilon)$ be a worldsheet with strip-like ends. Given a $1$-displacement function $h$, one can consider the analogue of \eqref{eq:gromov-2}, where the target space $\C$ carries the standard complex structure:
\begin{equation} \label{eq:v0-equation}
\left\{
\begin{aligned}
& v: S \longrightarrow \C, \\
& v(z) \in l_z \text{ for all $z \in \partial S$}, \\
& \bar\partial v = (X_h \otimes \gamma)^{0,1}.
\end{aligned}
\right.
\end{equation}

Occasionally, we will need perturbations in a slightly more general class. Let $k$ be a section of the pullback bundle $T^*\!S = \mathit{Hom}({TS},\R) \rightarrow S \times \C$. One can associate to it a section $X_k$ of $\mathit{Hom}_\R({TS},\C) \rightarrow S \times \C$, so that for each $\sigma \in {TS}$, $X_k(\sigma)$ is the Hamiltonian vector field of the function $k(\sigma)$ on $\C$. For $v$ with boundary conditions as in \eqref{eq:v0-equation}, one can then consider the equation
\begin{equation} \label{eq:v-equation}
\bar\partial v = X_k^{0,1}.
\end{equation}
To preserve important features of the simpler case $k = h \gamma$, we need to impose a number of conditions:

\begin{definition} \label{def:base-perturbation}
An {\em inhomogeneous term} is a $k$ such that for each $\sigma \in {TS}$, $k(\sigma)$ is a $\gamma(\sigma)$-displacement function. There are also restrictions on the behaviour of $k$ on the boundary and strip-like ends of $S$. Namely, for any $z \in \partial S$ and any vector $\sigma$ tangent to $\partial S$ at that point, the restriction of $k(\sigma)$ to $l_z$ should vanish. Finally, over each strip-like end, we ask that $\epsilon_\zeta^* k = k_{\zeta}(t,w)\,\mathit{dt}$ for some $k_\zeta \in \smooth([0,1] \times \C,\R)$ (each $k_{\zeta}(t,\cdot)$ is then automatically a $g_\zeta(t)$-displacement function).
\end{definition}

Since any displacement function depends only on $p$, one can alternatively write \eqref{eq:v-equation} as follows. Consider $\kappa = dk/dp$, which is again a section of $T^*\!S \rightarrow S \times \C$, depending only on the real part of the $\C$ variable. Then
\begin{equation} \label{eq:k-written-out}
X_k^{0,1} = (i\kappa)^{0,1} = \half(i\kappa - *\kappa)
\end{equation}
where $*$ is the Hodge operator on $\Omega^1(S)$ (which is composition with $j_S$, hence conformally invariant), and $\kappa$ is evaluated at point $(z,v(z))$. The special case \eqref{eq:v0-equation} then corresponds to taking $\kappa = (dh/dp) \gamma$. The point of view given by \eqref{eq:k-written-out} is maybe more elementary, but also more remote from standard Floer-theoretic considerations, and we'll therefore usually avoid it.

In parallel with the discussion in \eqref{eq:gromov-2-limits}, the natural asymptotic condition for solutions of \eqref{eq:v-equation} is $\lim_{s \rightarrow \pm\infty} v(\epsilon_\zeta(s,t)) = v_\zeta(t)$ for
\begin{equation} \label{eq:v-limits}
\left\{
\begin{aligned}
& v_\zeta: [0,1] \longrightarrow \C, \\
& dv_\zeta/dt = X_{k_\zeta}(t), \\
& v_\zeta(k) \in l_{\zeta,k} \text{ for $k = 0,1$.}
\end{aligned}
\right.
\end{equation}
Generalizing our previous notation slightly, let $\phi_{k_\zeta}$ be the time-one map of the Hamiltonian isotopy associated to the non-autonomous Hamiltonian $k_\zeta$. Then solutions $v_\zeta$ of \eqref{eq:v-limits} correspond bijectively to points $w_\zeta = v_\zeta(1) \in \phi_{k_\zeta}(l_{\zeta,0}) \cap l_{\zeta,1}$. Since $X_{k_\zeta} = g_\zeta(t) \partial_q$ on $W_2$, the situation on that subset is as in our original discussion of \eqref{eq:nondeg}, meaning that there are no intersection points. This, together with the boundedness of $X_{k_\zeta}$ and the given nature of the vanishing paths, implies that all the $v_\zeta$ lie in a compact subset.

\begin{remark}
The reader may have noticed that in our discussion of Floer cohomology in the fibre, we have limited ourselves to inhomogeneous terms of a relatively simple form \eqref{eq:gromov-2}, while here (and more importantly, in the subsequent discussion in Section \ref{subsec:m}) we allow more general terms \eqref{eq:v-equation}, \eqref{eq:u-equation}. In part, this is due to the difference in expository context. In the case of Floer cohomology for compact exact Lagrangian submanifolds, we were building on an existing body of literature (for instance, it was clear that Floer cohomology groups are independent of the choice of Hamiltonian perturbation), while here we are building a variation on that theory from scratch, hence have to leave enough freedom to accommodate such arguments, see for instance \eqref{eq:cont-pm}. Of course, the added flexibility also makes transversality arguments easier. However, the decisive reason is that we will ultimately end up using the more general form of inhomogeneous terms in the proof of the main result (see Remark \ref{th:general-k} and Section \ref{subsec:ccc}).
\end{remark}

\subsection{A priori bounds\label{subsec:a-priori-1}}
Fix a worldsheet with strip-like ends, as well as an inhomogeneous term $k$, and consider solutions $v$ of \eqref{eq:v-equation}. From now on, it is implicitly assumed that all such solutions have limits \eqref{eq:v-limits}.

\begin{lemma} \label{th:v-semicircle}
Take a constant $R>0$, such that the part of all vanishing paths $l_z$ lying in $\C \setminus W_{-1}$ is contained in the disc of radius $R$ around the origin. Then, if for any point $z$ we have $\mathrm{re}(v(z)) < -1$, then also $|v(z)| \leq R$.
\end{lemma}

\proof Since the boundary of $v$ lies on the vanishing paths $l_z$, and all its limits $v_\zeta$ either lie in the interior of $W_{-1}$ or else are intersection points of $l_z$, the degree of $v$ is a locally constant function on $\C \setminus (W_{-1} \cup \bigcup_z l_z)$. Since the image of $v$ is bounded, that function must necessarily vanish at all points $w \in \C \setminus W_{-1}$ such that $|w| > R$. But $v$ is holomorphic on the subset $v^{-1}(\C \setminus W_{-1})$, hence the vanishing degree means that $v^{-1}(w) = \emptyset$. \qed

\begin{lemma} \label{th:v-maximum}
For any $v$, we have $\mathrm{re}(v) \leq 2$ everywhere.
\end{lemma}

\proof We've already seen that the limits lie in the region $\mathrm{re}(v) < 2$. Suppose that $\mathrm{re}(v(z)) > 2$ at some point $z \in S$. Choose holomorphic local coordinates $z = s+it$ near that point. Then, \eqref{eq:v-equation} reduces to $\partial_t v - i \partial_s v = \gamma(\partial_s) + i \gamma(\partial_t)$, which in turn implies the Laplace equation
\begin{equation} \label{eq:laplace}
\begin{aligned}
\Delta \mathrm{re}(v) = - d\gamma(\partial_s,\partial_t) \geq 0.
\end{aligned}
\end{equation}
If $z$ additionally happens to be a boundary point of $S$, and the local coordinates around it are chosen so that $\partial S = \{t = 0\}$, we have
\begin{equation} \label{eq:normal-d}
\partial_t \mathrm{re}(v) = -\partial_s \mathrm{im}(v) + \gamma(\partial_s) = 0,
\end{equation}
because $\partial_s v$ is tangent to $l_z$ where that path is horizontal, and $\gamma$ vanishes in directions tangent to the boundary. Applying the maximum principle (with Neumann boundary conditions) leads to a contradiction. \qed

There are two versions of energy for a solution of \eqref{eq:v-equation}, the geometric and topological energies
\begin{align} \label{eq:energy}
& E^{\geom}(v) = \int_S {\textstyle\frac{1}{2}} |dv - X_k|^2, \\
& E^{\top}(v) = \int_S v^*\omega_\C - d(v^*k).
\end{align}
In the second line, we are thinking of $k$ as a one-form on $S \times \C$, and pull that back by the graph of $v$. The relation between the two energies is that
\begin{equation} \label{eq:e-e}
E^{\geom}(v) = E^{\top}(v) + \int_S v^*R.
\end{equation}
The curvature $R$ is a section of $\Lambda^2(T^*S) \rightarrow S \times \C$, obtained by taking the exterior derivative of $k$ in $S$-direction.
%
As before, we pull it back by the graph of $v$ to get the integrand in \eqref{eq:e-e}.
%
In the simpler situation of \eqref{eq:v0-equation} these expressions reduce to
\begin{align} \label{eq:simple-energy}
& E^{\geom}(v) = \int_S {\textstyle{\frac{1}{2}}} |dv - X_h \otimes \gamma|^2 = \int_S v^*\omega_\C - d(v^*h) \wedge \gamma, \\
& E^{\top}(v) = \int_S v^*\omega_\C - d(v^*h \cdot \gamma), \\
& R = h \, d\gamma. \label{eq:simplified-small-curvature}
\end{align}

It follows from Definition \ref{def:base-perturbation} that the curvature is bounded on $S \times \overline{\C \setminus W_2}$, and vanishes on the strip-like ends. We therefore get an inequality
\begin{equation} \label{eq:geom-top-constant}
E^{\geom}(v) \leq E^{\top}(v) + C
\end{equation}
for some constant $C$ which is independent of the particular $v$. On the other hand, the topological energy can be explicitly determined by integrating out. On each vanishing path $l_z$ which occurs as boundary condition, choose a function $f_z$ (locally constant in $z$) such that $df_z = \theta_\C|l_z$. If $w_\zeta$ is the point corresponding to a flow line $v_\zeta$ as in \eqref{eq:v-limits}, define its action to be
\begin{equation} \label{eq:little-action}
A(w_\zeta) = \Big( \int_{[0,1]} -v_\zeta^*\theta_\C + k_\zeta(t,v_\zeta(t)) \mathit{dt}\Big)
+ f_{\zeta,1}(v_\zeta(1)) - f_{\zeta,0}(v_\zeta(0)),
\end{equation}
where $f_{\zeta,t}$ are the functions associated to $l_{\zeta,t}$. Then
\begin{equation} \label{eq:energy-action}
E^{\top}(v) =
\sum_{\zeta \in \Sigma_{\mathrm{out}}} A(w_\zeta) - \sum_{\zeta \in \Sigma_{\mathrm{in}}} A(w_\zeta).
\end{equation}
By comparing this with \eqref{eq:geom-top-constant}, and using the compactness of the set of possible limits, one sees that there is a uniform bound on $E^{\geom}(v)$ for all $v$. Given that, we will use a baby version of Gromov compactness to prove the following:

\begin{proposition} \label{th:v-bounded}
For a given $(S,\gamma,l,\epsilon)$ and $k$, all solutions $v$ of \eqref{eq:v-equation} remain within a bounded subset of $\C$.
\end{proposition}

\proof From Lemmas \ref{th:v-semicircle} and \ref{th:v-maximum}, we know a priori that any $v$ remains in the region
\begin{equation} \label{eq:b-region}
U_{-1,2} \cup \{|w| \leq R \text{ and } p = \mathrm{re}(w) \leq -1\} \subset \C
\end{equation}
for a suitable constant $R$. Suppose that we have a sequence of solutions $v_j$ such that $\|v_j\|_\infty$, hence also $\|\mathrm{im}(v_j)\|_{\infty}$, goes to infinity. Then, the point where $|\mathrm{im}(v_j)|$ achieves its maximum must be increasingly far away from the boundary conditions. On the other hand, on $S$ itself every point is within a bounded distance of the boundary, with respect to any metric which is standard on the strip-like ends. By comparing these two facts, one sees that $\|dv_j\|_\infty$ necessarily goes to infinity as well.

Take the points $z_j$ where $|dv_j|$ reaches its maximum, for $j \gg 0$. If $v_j(z_j)$ has a bounded subsequence, we can apply the standard rescaling argument to obtain a limit which is either a non-constant holomorphic sphere $\C P^1 \rightarrow \C$, or else a non-constant holomorphic disc $(D,\partial D) \rightarrow (\C,l_z)$ for some $z \in \partial S$, which is a contradiction.

Now consider the other situation, where $|\mathrm{im}(v_j(z_j))|$ goes to infinity. If $z_j$ itself is convergent (possibly after passing to a subsequence), we work in local holomorphic coordinates near its limit point on $S$. In the other case, (possibly after passing to a subsequence) we can assume that all $z_j$ lie inside the strip-like ends, and can use the standard coordinates on those ends. In either case, we proceed as follows: fix some $\rho$, and take a closed disc of radius $\rho/\|dv_j\|_\infty$ around $z_j$. Assuming that $j$ is large, this disc lies in the interior of $S$, because the distance of $v_j(z_j)$ to $v_j(\partial S)$ is greater than $\rho$. Moreover, the image of this disc under $v_j$ lies in $U_{-1,2}$, because the distance of $v_j(z_j)$ to the other part $\{|w| \leq R\}$ of \eqref{eq:b-region} is again greater than $\rho$. We restrict $v_j$ to this disc, rescale the domain by $\|dv_j\|_{\infty}$, and subtract a constant which is the imaginary part of $v_j(z_j)$. The limit of a subsequence is a holomorphic map $\tilde{v}$ defined on a slightly smaller disc, satisfying $\tilde{v}(0) \in [-1,2] \subset \R$ and $|d\tilde{v}(0)| = \|d\tilde{v}\|_{\infty} = 1$, and whose energy is still bounded by the original upper bound on geometric energies. Repeat this process with subsequently larger values of $\rho$, passing to further subsequences. In the end, one again gets a non-constant holomorphic function $\C P^1 \rightarrow \C$, hence a contradiction. \qed

\subsection{Barriers\label{subsec:barrier}}
Let $(S,\gamma,l,\epsilon)$ be a worldsheet with strip-like ends, and $k$ an inhomogeneous term.

\begin{definition} \label{th:def-barrier}
A {\em barrier} is a vertical line $B = \{p = b\} \subset \C$, where $b \geq -2$, such that
\begin{equation} \label{eq:barrier-def}
X_k | S \times B = \partial_q \otimes \kappa
\end{equation}
for some one-form $\kappa \in \Omega^1(S)$ (note that unlike the general observation in \eqref{eq:k-written-out}, $\kappa$ may not additionally depend on points in $\C$). We say that the barrier is closed if $d\kappa = 0$ and $\kappa|\partial S = 0$.
\end{definition}

Let $v$ be a solution of \eqref{eq:v-equation}. Assume that the limits $w_\zeta$ of $v$ are disjoint from $B$, and that $v$ intersects $B$ transversally (strictly speaking, by this we mean that both $v$ and $v|\partial S$ are transverse to $B$). These assumptions imply that $v^{-1}(B)$ is a compact one-dimensional submanifold of $S$. We orient that submanifold so that it is the boundary of $v^{-1}(\overline{\C \setminus W_b})$. Let $\nu$ be an outwards pointing normal vector field, by which we mean a section of ${TS}|v^{-1}(B)$ satisfying $\mathrm{re}(Dv(\nu)) = 1$. The same argument as in \eqref{eq:simple-energy} shows that
\begin{equation} \label{eq:barrier}
0 \leq \int_{v^{-1}(B)} \iota_\nu\big(\half|dv - X_k|^2\big) = \int_{v^{-1}(B)} v^*\mathit{dq} - \kappa.
\end{equation}
This is a local computation, meaning that the integrands are pointwise equal; in particular, it only depends on the fact that $v$ satisfies \eqref{eq:v-equation} along $v^{-1}(B)$.

The additional closedness condition on $\kappa$ introduced in Definition \ref{th:def-barrier} has not entered into our discussion so far. However, it is useful when it comes to applications, where one wants to use \eqref{eq:barrier} to restrict the global behaviour of solutions. We will explain this through examples (which should be considered as toy models for our later discussion in Section \ref{sec:comparison}).

\begin{example} \label{th:one-sided}
Let $v$ be a solution of \eqref{eq:v-equation} all of whose limits $w_\zeta$ lie on the same side of a closed barrier $B \subset \C$, and which intersects $B$ transversally. 
We claim that then, {\em $v$ never meets $B$}. To see this, take the closed half-plane bounded by $B$ which does not contain any $w_\zeta$, and let $A$ be its preimage under $v$. By assumption, $A$ is a compact subset of $S$, whose boundary is the union of $v^{-1}(B)$ and $A \cap \partial S$. The latter part consists of disjoint circles and intervals. Each circle is a compact connected component of $\partial S$, and the integral of $v^*\mathit{dq} - \kappa$ over it is necessarily zero. Similarly, each interval is mapped to some vanishing path $l_z$, and its endpoints are therefore preimages of the unique point of $l_z \cap B$, which means that the integral over it of $v^*\mathit{dq} - \kappa$ is again zero. Applying Stokes yields
\begin{equation} \label{eq:stokes}
\int_{v^{-1}(B)} v^*\mathit{dq} - \kappa = \pm \int_{\partial A} v^*\mathit{dq} - \kappa = 0,
\end{equation}
which in view of \eqref{eq:barrier} implies that $dv = X_k = \partial_q \otimes \kappa$ along $v^{-1}(B)$. This contradicts the fact that $v$ intersects $B$ transversally, unless $v^{-1}(B) = \emptyset$.
\end{example}

\begin{example} \label{th:1-cross}
Consider $S$ with a choice of distinguished output end $\zeta$, whose adjacent vanishing paths are $l_{\zeta,j}$. Recall that by definition, $l_{\zeta,j} \cap W_{-2} = \{q = q_{\zeta,j}(p)\}$ for some function $q_{\zeta,j}$. Suppose that we have a closed barrier, with the additional condition that $\epsilon_\zeta^*\kappa = k_\zeta(t) \mathit{dt}$ for some function $k_\zeta$ satisfying
\begin{equation} \label{eq:sign-of-translation}
\int_0^1 k_\zeta(t) \mathit{dt} < q_{\zeta,1}(b) - q_{\zeta,0}(b).
\end{equation}
Let $v$ be a solution of \eqref{eq:v-equation} such that $w_\zeta$ lies in the interior of $W_b$, all other limits lie in $\C \setminus W_b$, and which intersects $B$ transversally. Set $A = v^{-1}(W_b) \setminus \epsilon_\zeta((-\infty,s) \times [0,1])$, for $s \ll 0$. Using the same reasoning as in the previous example, we have
\begin{equation} \label{eq:stokes2}
\begin{aligned}
0 & = -\Big( \int_{\partial A} v^*\mathit{dq} - \kappa \Big) \\
& = \Big(\int_{v^{-1}(B)} v^*\mathit{dq} - \kappa\Big)
+ \Big( q_{\zeta,1}(b) - q_{\zeta,0}(b) \Big) - \Big( \int_0^1 k_\zeta(t) dt \Big).
\end{aligned}
\end{equation}
Here, the $q_{\zeta_j}(b)$ terms come from integrating $\mathit{dq}$ over $\{s\} \times [0,1]$ and the two adjacent intervals of $v^{-1}(W_b) \cap \partial S$, which together yield a path going from $l_{\zeta,0} \cap B$ to $l_{\zeta,1} \cap B$. This leads to a contradiction between \eqref{eq:sign-of-translation} and \eqref{eq:barrier}. Hence, {\em no such solutions $v$ can exist.}
\end{example}

\begin{remark} \label{th:sub-closed-barrier}
The closedness assumption can sometimes be weakened. For instance, in Example \ref{th:1-cross} it is sufficient if $d\kappa \leq 0$ (still assuming that $\kappa|\partial S = 0$). The same happens in Example \ref{th:one-sided}, provided that all limits $w_\zeta$ lie to the left of $B$ (correspondingly, if all $w_\zeta$ lie to the right of $B$, we want $d\kappa \geq 0$).
\end{remark}

\subsection{Cauchy-Riemann operators}
We begin with some facts from Sturm-Liouville theory (this is understood in the general sense of \cite{bott56, levitan-sargsjan, weidmann}, but our discussion covers only the most straightforwardly elementary case). Suppose that we are given the following data:
\begin{equation} \label{eq:selfadjoint-data}
\parbox{28em}{
One-dimensional linear real subspaces $\lambda_0,\lambda_1 \subset \C$, and a smooth family $a_t$, $t \in [0,1]$, of symmetric real 2x2 matrices, considered as elements of $\mathit{Hom}_\R(\C,\C)$.
}
\end{equation}
To these, we can associate the following operator $A$:
\begin{equation} \label{eq:selfadjoint}
\left\{
\begin{aligned}
& A Y = i (dY/dt) + a_t Y, \\
& Y: [0,1] \longrightarrow \C, \\
& Y(0) \in \lambda_0, \; Y(1) \in \lambda_1.
\end{aligned}
\right.
\end{equation}
This is a selfadjoint elliptic operator, hence has discrete spectrum and an $L^2$ basis of eigenvectors. Next, because the eigenvalue equation $(A-\mu)Y = 0$ is a first order ODE with one-dimensional space of possible boundary values, the eigenvalues are all simple. Finally, for the same reason, each eigenvector is nowhere vanishing, hence defines a path in $\RP{1}$ from the point $\lambda_0$ to $\lambda_1$.

\begin{lemma}
Each homotopy class of paths in $\RP{1}$ going from $\lambda_0$ to $\lambda_1$ is represented by exactly one eigenvector of $A$. Moreover, passing from one eigenvalue to a bigger one corresponds to composing the given homotopy class with a negative (with respect to the standard identification $\pi_1(\RP{1}) \iso \Z$) loop in $\RP{1}$.
\end{lemma}

\proof The statement is straightforward for the case $\lambda_0 = \lambda_1 = \R$ and $a = 0$. The general case can be reduced to this by a deformation. Since that deformation remains within the same class of operators, all eigenvalues remain simple during it, and can never cross each other; this implies that all eigenvectors vary continuously, that the associated homotopy classes are preserved, and that the order-reversing nature of the correspondence between eigenvalues and homotopy classes is retained. \qed

We will now impose the additional condition that $AY = 0$ should have no nonzero solutions. Equivalently, if $g_1$ is the path in $\mathit{Sp}(2,\R)$ defined by $g_0 = \mathrm{Id}$, $(dg_t/dt) g_t^{-1} = i a_t$, then
\begin{equation} \label{eq:transverse-linearized}
g_1(\lambda_0) \neq \lambda_1.
\end{equation}
Given any path $[0,1] \rightarrow \RP{1}$, one can define its (total signed) angle, which is a real number invariant under homotopies rel endpoints. Our normalization condition is that the loop going once (anticlockwise) around $\RP{1}$ has angle $\pi$. Denote by $\measuredangle$ the angle of the path $t \mapsto g_t(\lambda_0)$. By \eqref{eq:transverse-linearized}, any path from $\lambda_0$ to $\lambda_1$ must have angle in $\R \setminus (\measuredangle + \pi\Z)$.

\begin{lemma} \label{th:signed-eigenvalue}
Among the eigenvectors of $A$, those with negative eigenvalues are precisely those whose associated paths in $\RP{1}$ have total angle greater than $\measuredangle$.
\end{lemma}

\proof The proof is again by a deformation argument. We have to be a little more careful to preserve the additional condition \eqref{eq:transverse-linearized}, which ensures both that eigenvalues never cross zero, and that the angles associated to eigenvectors never cross $\measuredangle$. Nevertheless, within that space, one can always deform to a situation where $a_t = \measuredangle\, \mathrm{Id}$ is a constant multiple of the identity. There, the eigenvector associated to an eigenvalue $\mu$ is $Y(t) = \exp(i (\measuredangle-\mu) t)$. The path in $\RP{1}$ associated to $Y$ has angle $\measuredangle-\mu$, which clearly shows the desired property. \qed

The class of differential operators considered above appears when analyzing the asymptotic behaviour of perturbed Cauchy-Riemann operators on Riemann surfaces with strip-like ends. More precisely, we consider the following context:
\begin{equation}
\parbox{28em}{
$S$ is a punctured-boundary Riemann surface equipped with strip-like ends $\{\epsilon_\zeta\}$. Additionally, it should come with a family $\lambda: \partial S \rightarrow \RP{1}$ of one-dimensional subspaces of $\C$, which is locally constant over the part of $\partial S$ which belongs to the strip-like ends. Finally, we want to have a one-form $a$ on $S$ with values in $\mathit{Hom}_\R(\C,\C)$, whose restriction to each strip-like end is of the form $a_{\zeta,t} \mathit{dt}$, where the $a_{\zeta,t}$ are symmetric matrices.
}
\end{equation}
We can associate to this the Cauchy-Riemann operator
\begin{equation} \label{eq:simple-linearization}
\left\{
\begin{aligned}
& DY = \bar\partial Y - (ia Y)^{0,1}, \\
& Y: S \longrightarrow \C, \\
& Y(z) \in \lambda_z \text{ for all $z \in \partial S$.}
\end{aligned}
\right.
\end{equation}
Here, $a Y$ is the complex-valued one-form obtained by multiplying the matrix (of one-forms) $a$ and vector $Y$. Note that at each end, we have the same data as in \eqref{eq:selfadjoint-data}, namely: two linear subspaces $\lambda_{\zeta,k}$, corresponding to the points $\epsilon_\zeta(s,k)$ for any $s$; and the family of matrices $a_\zeta = (a_{\zeta,t})$. If $A_\zeta$ is the associated differential operator \eqref{eq:selfadjoint}, then the restriction of \eqref{eq:simple-linearization} to any end, where $Y_\zeta(s,t) = Y(\epsilon_\zeta(s,t))$, is given by
\begin{equation} \label{eq:simple-on-the-end}
D_\zeta Y_\zeta = \partial_s Y_\zeta + i (\partial_t Y_\zeta - i a_{\zeta,t} Y_\zeta) = \partial_s Y_\zeta + A_\zeta Y_\zeta.
\end{equation}

Let $(g_{\zeta,t})$ be the path in $\mathit{Sp}(2,\R)$ obtained as before by integrating $(a_{\zeta,t})$. We impose the additional assumption that $g_{\zeta,1}(\lambda_{\zeta,0}) \neq \lambda_{\zeta,1}$, which is just saying that \eqref{eq:transverse-linearized} holds at each end. This implies that the Sobolev completion of \eqref{eq:simple-linearization} (taking $W^{1,p}$ on the domain, and $L^p$ on the range), which we denote by $D: \E^1 \rightarrow \E^0$, is a Fredholm operator.
Let $\hat{S}$ be the compactification of $S$ obtained by adding an interval $\epsilon_\zeta(\{\pm \infty\} \times [0,1])$ to each end. We extend $\lambda$ to a map
\begin{equation} \label{eq:hat-lambda}
\hat\lambda: \partial \hat{S} \longrightarrow \RP{1},
\end{equation}
in a way which is unique up to homotopy. Doing that just means singling out, for each $\zeta$, a homotopy class of paths from $\lambda_{\zeta,0}$ to $\lambda_{\zeta,1}$. As before, let $\measuredangle_\zeta$ be the angle of the path $g_{\zeta,t}(\lambda_{\zeta,0})$. If $\zeta$ is an input, we take the homotopy class whose angle is smaller than $\measuredangle_\zeta$, but as big as possible within that restriction. If $\zeta$ is an output, we similarly take the homotopy class whose angle is bigger than $\measuredangle_\zeta$, but as small as possible within that restriction. We can now state the index formula for $D$, which in its most basic form goes back to \cite{floer88e}:
\begin{equation} \label{eq:index-formula}
\mathrm{index}(D) = \chi(\hat{S}) + \mathrm{deg}(\hat\lambda),
\end{equation}
where the second term is the degree of \eqref{eq:hat-lambda} with respect to the standard orientations  of $\partial\hat{S}$ and $\RP{1}$.

\begin{prop} \label{th:negative-degree}
If $\mathrm{deg}(\hat\lambda) < 0$, $D$ is injective.
\end{prop}

This generalizes the classical fact that a negative degree line bundle over a closed Riemann surface can't have nonzero holomorphic sections (compare also \cite[Lemma 11.5]{seidel04}). To simplify the notation, let's assume that all ends of our Riemann surface are inputs; it is not difficult to see that the definition of $\hat\lambda$ above, in spite of (in fact, because of) its apparent asymmetry, is compatible with relabeling outputs as inputs.

\proof Suppose that $Y$ is a nontrivial solution of $DY = 0$. This has isolated zeros, each of which has finite multiplicity (by \cite[Theorem 2.2]{floer-hofer-salamon94} and its boundary counterpart \cite[Theorem 2.1]{oh96c}). Denote these multiplicities by $m_z$, with the convention that the multiplicity for boundary points is half the order of vanishing, hence usually only a half-integer. Now consider some end, with its associated selfadjoint operator $A_\zeta$. Denote by $(Y_{\zeta,k})$ a basis of eigenvectors for that operator, labeled in increasing order of their eigenvalues, and such that $Y_{\zeta,0}$ corresponds to the smallest positive eigenvalue. In view of Lemma \ref{th:signed-eigenvalue}, this means that $Y_{\zeta,0}$ lies in the same homotopy class as the path used to define $\hat\lambda(\epsilon_\zeta(\infty,\cdot))$. All the other $Y_{\zeta,k}$ differ from it by adding $-k$ times the generators of $\pi_1(\RP{1})$.

Expand $Y_\zeta(0,\cdot) = Y(\epsilon_\zeta(0,\cdot))$ with respect to this basis of eigenvectors. Because $Y_\zeta(s,t)$ decays as $s \rightarrow \infty$, only the $Y_{\zeta,k}$ with $k \geq 0$ can contribute. Let $k_\zeta$ be the index corresponding to the smallest eigenvector that appears with nonzero coefficient $r_\zeta \neq 0$, and $\mu_\zeta > 0$ its eigenvalue. Then
\begin{equation} \label{eq:rescaled-decay}
\textstyle \mathrm{lim}_{s \rightarrow \infty} e^{\mu_\zeta s} Y_\zeta(s,\cdot) = r_\zeta Y_{\zeta,k_\zeta}.
\end{equation}
In particular, $Y \neq 0$ outside a compact subset of $S$. An elementary winding number argument then shows that
\begin{equation}
0 \leq 2 \!\!\sum_{\text{zeros $z$ of $Y$}}\!\! m_z = \mathrm{deg}(\hat\lambda) - \!\! \sum_{\text{ends $\zeta$}} k_\zeta,
\end{equation}
which is a contradiction if $\mathrm{deg}(\hat\lambda) < 0$.
\qed

\subsection{Linearization}
Let $v$ be a solution of an equation \eqref{eq:v-equation}, with limits $v_\zeta$ as in \eqref{eq:v-limits}. The linearization of \eqref{eq:v-equation} at $v$, written in analogy with \eqref{eq:k-written-out} for simplicity, is the operator
\begin{equation} \label{eq:explicit-linearization}
\left\{
\begin{aligned}
& D_v(Y) = \bar\partial Y - (i \mathrm{re}(Y) \otimes \delta)^{0,1}, \\
& Y: S \longrightarrow \C, \\
& Y(z) \in (Tl_z)_{v(z)} \text{ for all $z \in \partial S$.}
\end{aligned}
\right.
\end{equation}
where $\delta = (d^2 k/dp^2)(z,v(z)) \in \Omega^1(S)$. Setting $\lambda_z = (Tl_z)_{v(z)}$ and
\begin{equation}
a = \begin{pmatrix} \delta & 0 \\ 0 & 0 \end{pmatrix}
\end{equation}
(here we're identifying $\mathit{Hom}_\R(\C,\C) = \mathit{Hom}_\R(\R^2,\R^2) = \mathit{Mat}_2(\R)$, whence the matrix notation) turns this into the form \eqref{eq:simple-linearization}. We can still associate to each end a translation-invariant operator on $\R \times [0,1]$, namely
\begin{equation} \label{eq:asymptotic-operator}
\left\{
\begin{aligned}
& D_{v,\zeta}(Y_\zeta) = \partial_s Y_\zeta + i\partial_t Y_\zeta + \delta_{\zeta}(t) \mathrm{re}(Y_\zeta), \\
& Y_\zeta: \R^\pm \times [0,1] \longrightarrow \C, \\
& Y_\zeta(s,k) \in (Tl_{\zeta,k})_{v_\zeta(k)}.
\end{aligned}
\right.
\end{equation}
where $\delta_\zeta = d^2 k_{\zeta}/dp^2$ evaluated at $(t,v_\zeta(t))$. This is again of the same form as \eqref{eq:simple-on-the-end}, with $\lambda_{\zeta,k} = (Tl_{\zeta,k})_{v_\zeta(k)}$ and $a_\zeta = \mathrm{diag}(\delta_\zeta,0)$. However, this time \eqref{eq:explicit-linearization} is only asymptotically equal to \eqref{eq:asymptotic-operator}. In spite of that, we'll see that results from the existing literature allow us to extend the previous argument to cover this situation.

Assume from now on that the $v_\zeta$ are nondegenerate, which means that they correspond to transverse intersection points $w_\zeta \in \phi_{k_\zeta}(l_{\zeta,0}) \cap l_{\zeta,1}$. Then the Sobolev completion of \eqref{eq:explicit-linearization} is an elliptic operator $D_v: \E^1_v \rightarrow \E^0_v$, whose index can be computed as in \eqref{eq:index-formula} (from which it follows by a deformation argument). Results of \cite{agmon-nirenberg67}, as applied to this problem in \cite{robbin-salamon02}, show that the limiting behaviour of solutions of $D_v Y = 0$ remains the same as in \eqref{eq:rescaled-decay} \cite[Lemma 3.5]{robbin-salamon02}. Hence, the proof of Proposition \ref{th:negative-degree} goes through as before, which in particular implies the following {\em automatic regularity} result:

\begin{corollary} \label{th:automatic-regularity}
Suppose that $S$ is a disc with some boundary punctures. Suppose also that the intersections $\phi_{k_\zeta}(l_{\zeta,0}) \cap l_{\zeta,1}$ are all transverse. Then, every index zero solution of \eqref{eq:v-equation} has the property that its linearized operator is invertible. \qed
\end{corollary}

\section{Geometry of the total space}

\subsection{Lefschetz fibrations}
Here is a preliminary notion of Lefschetz fibration, in an essentially topological context. Define an {\em almost complex Lefschetz fibration} $\pi: E \rightarrow \C$ as follows: $E$ is a $(2n+2)$-dimensional manifold, possibly with boundary, equipped with an almost complex structure $I_E$. The map $\pi: E \rightarrow \C$ is proper and pseudo-holomorphic (with respect to $I_E$), and $\pi|\partial E$ should be a submersion. Moreover, we assume that $\pi$ has finitely many critical points, at most one of which lies in each fibre. Near each critical point $x$, $I_E$ must be integrable, and the complex Hessian $D^2_x\pi: TE_x \rightarrow \C$ at that point must be nondegenerate (as a complex quadratic form). Properness of $\pi$, together with the additional condition along the boundary, ensures that the nonsingular part of any almost complex Lefschetz fibration is a differentiable fibre bundle. Near each critical point, the integrability assumption allows one to apply the complex Morse Lemma \cite{atiyah58}, which gives the standard normal form $\pi(x) = \sum_j x_j^2 + (\mathit{constant})$ in local holomorphic coordinates.

To make such a $\pi: E \rightarrow \C$ into a {\em symplectic Lefschetz fibration}, we assume that additionally, $E$ comes with a symplectic structure $\omega_E$ compatible with $I_E$. We also require that the horizontal tangent spaces $TE^h_x$ (defined as the orthogonal complements of $TE^v_x = \mathit{ker}(D\pi_x)$ for any regular point $x$) are tangent to the boundary. This implies that $TE^h$, which is a Hamiltonian connection away from the singular fibres, has well-defined parallel transport maps. Finally, for our purposes we want to add some more conditions:

\begin{definition} \label{def:exact-lefschetz}
An {\em exact symplectic Lefschetz fibration} is a symplectic Lefschetz fibration with the following additional properties. $\omega_E = d\theta_E$ is exact, and $\partial E$ is weakly convex (compare Definition \ref{def:exact-symplectic}).  Moreover, the fibration should be trivial near infinity in the following sense (see \eqref{eq:uw-notation} for notation): there are no critical points in $\pi^{-1}(W_{-2})$; the derivative $D\pi_x: TE^h_x \rightarrow \C$ is a symplectic linear isomorphism for all $x \in \pi^{-1}(W_{-2})$; and parallel transport along any path in $W_{-2}$ preserves almost complex structures on the fibres.
\end{definition}

Triviality near infinity can be equivalently expressed as follows. Pick any point $z \in W_{-2}$, and denote by $M$ the fibre over that point, with associated structures $I_M$ and $\omega_M$. Then there is a partial trivialization
\begin{equation} \label{eq:trivialization-at-infty}
\pi^{-1}(W_{-2}) \stackrel{\iso}{\longrightarrow} W_{-2} \times M
\end{equation}
which takes $I_E$ to the product structure $i \times I_M$, $\omega_E$ to $\omega_\C + \omega_M$, and whose restriction to $\pi^{-1}(z) = M$ is the identity (it is clear that the existence of such a partial trivialization implies the properties stated in Definition \ref{def:exact-lefschetz}; to go in the converse direction, use parallel transport). From now on, we'll work with a fixed exact symplectic Lefschetz fibration $\pi: E \rightarrow \C$.

\begin{definition} \label{def:vanishing-2}
A {\em vanishing path} is an $l$ as in Definition \ref{def:vanishing-1}, whose endpoint is a critical value of $\pi$, and such that all other points on $l$ are regular values. We also want to adapt the notion of {\em simply intersecting basis} $(l_1,\dots,l_m)$ accordingly. Namely, each critical value of $\pi$ should be the endpoint of exactly one $l_i$.
\end{definition}

For any vanishing path in the revised sense defined above, there is a unique properly embedded Lagrangian submanifold $L \subset E$, called the Lefschetz thimble of $l$, such that $\pi(L) = l$ (see e.\ g.\ \cite[Section 1.3]{seidelexact}). It is always diffeomorphic to $\R^{n+1}$. Moreover, with respect to \eqref{eq:trivialization-at-infty}, we have
\begin{equation}
L \cap \pi^{-1}(W_{-2}) = (l \cap W_{-2}) \times V
\end{equation}
where $V \subset M$ is a Lagrangian sphere, called the vanishing cycle. If $(l_1,\dots,l_m)$ is a simply intersecting basis, we also refer to the associated collections $(L_1,\dots,L_m)$ and $(V_1,\dots,V_m)$ as bases.

Fix some constant $c$, and consider functions with the following properties.
\begin{equation} \label{eq:h-prime}
\parbox{28em}{
$\tilde{H} \in \smooth(E,\R)$ is supported in $\pi^{-1}(W_{-1})$ and invariant under translation in $\R$-direction, hence can be thought of a function of $(p,x) \in [-1, \infty) \times M$. Moreover, in a neighbourhood of $[-1,\infty) \times \R \times \partial M$, $\tilde{H}$ is a function of $p$ alone. Finally, on $[2,\infty) \times M$ we have $d\tilde{H}/dp = c$.
}
\end{equation}
Obviously, one way to satisfy these requirements is to take a $c$-displacement function $h$ on $\C$, in the sense of Definition \ref{def:displacement-1}, and pull that back to $E$. Somewhat more generally, one can take $\tilde{H}(p,x) = h(p) + \psi(p)F(x)$, where $\psi$ is a cutoff function satisfying $\psi(p) = 0$ for $p \leq -1$ and $\psi(p) = 1$ for $p \geq 2$, and $F$ is any function on $M$ which vanishes near the boundary.

\begin{definition} \label{def:displacement-2}
A {\em $c$-displacement function} is an $H$ which agrees with some $\tilde{H}$ as in \eqref{eq:h-prime} outside a compact subset of $E \setminus \partial E$.
\end{definition}

As usual, we write $X_H$ for its Hamiltonian vector field and $\phi_H$ for the time-one map. The motivation for allowing extra freedom, compared to \eqref{eq:h-prime}, is transversality of intersections. Namely, suppose that $l_0,l_1$ are vanishing paths, and $L_0,L_1$ the associated Lefschetz thimbles. For any $c$ satisfying \eqref{eq:nondeg}, we can find a $c$-displacement function $H$ such that $\phi_H(L_0) \cap L_1$ consists of finitely many transverse intersection points.

\subsection{Maps into the total space\label{subsec:m}}
Our first task is to define a suitable class of almost complex structures, which we will do in two steps following the model of Definition \ref{def:displacement-2}. First, consider $\tilde{J}$ as follows:
\begin{equation} \label{eq:j-prime}
\parbox{28em}{
$\tilde{J}$ is an $\omega_E$-compatible almost complex structure, which agrees with $I_E$ in a neighbourhood of $\partial E$, and which has the following properties. First, $\pi$ is $\tilde{J}$-holomorphic. Second, the restriction of $\tilde{J}$ to $\pi^{-1}(W_{-1}) \iso W_{-1} \times M = [-1,\infty) \times \R \times M$ is invariant under translation in $\R$-direction.
}
\end{equation}

\begin{definition}
$\JJ(E)$ is the space of all $\omega_E$-compatible almost complex structures $J$ for which there exists a $\tilde{J}$ as in \eqref{eq:j-prime}, such that $\tilde{J} = J$ outside a compact subset of $E \setminus \partial E$.
\end{definition}

Take a worldsheet with strip-like ends, $(S,\gamma,l,\epsilon)$, and a $1$-displacement function $H$. To the boundary condition $l = (l_z)$ which is part of its structure, one can associate the corresponding family $L = (L_z)$ of Lefschetz thimbles. Choose a family $J = (J_z)$ in $\JJ(E)$ depending on $z \in S$, such that on the strip-like ends $J_{\epsilon_\zeta(s,t)} = J_{\zeta,t}$ for some $J_\zeta = (J_{\zeta,t})$ parametrized by $t \in [0,1]$. Given such data we can write down an analogue of \eqref{eq:v0-equation}:
\begin{equation} \label{eq:u0-equation}
\left\{
\begin{aligned}
& u: S \longrightarrow E, \\
& u(z) \in L_z \text{ for all $z \in \partial S$}, \\
& (du - X_H \otimes \gamma)^{0,1} \\ & \quad = \half(du(z) - X_H \otimes \gamma + J_z \circ du(z) \circ j_S - J_z X_H \otimes \gamma \circ j_S) = 0.
\end{aligned}
\right.
\end{equation}

While this is sufficient for basic purposes, we will also need a more flexible version, which follows the model of \eqref{eq:v-equation} by allowing inhomogeneous terms of the form
\begin{equation} \label{eq:u-equation}
(du - X_K)^{0,1} = 0.
\end{equation}
Here, $K$ is a section of the bundle $T^*\!S \rightarrow S \times E$, and we are using the associated section $X_K$ of $\mathit{Hom}_\R(TS,TE) \rightarrow S \times E$. The set of allowed $K$ is defined as follows. First, consider sections $\tilde{K}$ as follows:
\begin{equation} \label{eq:k-prime}
\parbox{28em}{
For any $\sigma \in TS$, the function $\tilde{H} = \tilde{K}(\sigma)$ is as in \eqref{eq:h-prime}, with constant $c = \gamma(\sigma)$.
}
\end{equation}

\begin{definition} \label{def:k}
An {\em inhomogeneous term} is a $K$ which agrees with a $\tilde{K}$ as in \eqref{eq:k-prime} outside $S \times (\text{\it a compact subset of $E \setminus \partial E$})$. Additionally, there are restrictions on the behaviour of $K$ on $\partial S$ and over the strip-like ends. Namely, if $\sigma \in T(\partial S)_z$, then $K(\sigma)$ must vanish on $L_z$. Finally, as in Definition \ref{def:base-perturbation} we require that $\epsilon_\zeta^*K = K_\zeta(t)\, \mathit{dt}$, where $K_\zeta$ is a function on $[0,1] \times E$. Note that each $K_\zeta(t)$ is then automatically a $g_\zeta(t)$-displacement function.
\end{definition}

If we are given such a $K$ as well as a $J$ as before, the pair $(J,K)$ is called a {\em perturbation datum}, since this is what's needed to write down the perturbed pseudo-holomorphic map equation \eqref{eq:u-equation}.

\begin{remark}
Readers having worked their way through the definitions above may be tempted to take their frustration out on the author in person. Before that idea takes root, let me give a quick preview of how these technical properties are used, starting with the more restrictively defined $(\tilde{J},\tilde{K})$ from \eqref{eq:j-prime}, \eqref{eq:k-prime}. On the region $\pi^{-1}(\C \setminus W_{-1})$, we have $D\pi \circ \tilde{J} = i \circ D\pi$ and $\tilde{K} = 0$. On $\pi^{-1}(W_2)$, we similarly have $D\pi \circ \tilde{J} = i \circ D\pi$ and $D\pi(X_{\tilde{K}}) = \partial_q \otimes \gamma$. In both cases, this allows one to argue by projection to the base, where methods from Section \ref{subsec:a-priori-1} can be applied. Finally, on $\pi^{-1}(W_{-1}) \iso [-1,\infty) \times \R \times M$, the main additional property of $(\tilde{J},\tilde{K})$ was independence of the $\R$ variable. This allows us to translate solutions in that direction, and thereby to restrict considerations to a bounded subset. The actual $(J,K)$ satisfy these properties only outside a compact subset, but that is sufficient for our main a priori estimate (Proposition \ref{th:u-bounded}).
\end{remark}

The natural limit condition for solutions of \eqref{eq:u-equation} is $\lim_{s \rightarrow \pm\infty} u(\epsilon_\zeta(s,t)) = u_\zeta(t)$, where
\begin{equation} \label{eq:u-limits}
\left\{
\begin{aligned}
& u_\zeta: [0,1] \longrightarrow E, \\
& du_\zeta/dt = X_{K_\zeta}(t), \\
& u_\zeta(k) \in L_{\zeta,k} \text{ for $k = 0,1$.}
\end{aligned}
\right.
\end{equation}
Here, $L_{\zeta,k}$ are the Lefschetz thimbles associated the boundary components adjacent to $\zeta$.
If $\phi_{K_\zeta}$ is the time-one map of the Hamiltonian isotopy induced by $K_\zeta$, then solutions of \eqref{eq:u-limits} correspond bijectively to points $y_\zeta \in \phi_{K_\zeta}(L_{\zeta,0}) \cap L_{\zeta,1}$.

\begin{lemma}
All limits $y_\zeta$ are contained inside a compact subset of $E \setminus \partial E$.
\end{lemma}

\proof It follows from Definition \eqref{def:k} that $D\pi(X_{K_{\zeta}}(t)) = g_\zeta(t) \partial_q$ on $\pi^{-1}(W_C)$ for $C \gg 0$. In view of \eqref{eq:rzeta-2}, this implies that the limits $y_\zeta$ must lie outside $\pi^{-1}(W_C)$. Next, we also know that $D\pi(X_{K_{\zeta}})$ is bounded. This, together with the position of the Lefschetz thimbles, ensures that all limits $y_\zeta$ must stay within a compact subset.

Since those limits themselves form a closed subset, all we need to show is that no limit can touch $\partial E$. But that is obvious, since $X_{K_\zeta}$ is always tangent to the boundary, while the Lefschetz thimbles are contained in the interior. \qed

\subsection{A priori bounds}
Consider solutions of \eqref{eq:u-equation}, which are always assumed to have limits \eqref{eq:u-limits}.
We now carry over the discussion from Section \ref{subsec:a-priori-1} to this situation.

\begin{lemma} \label{th:no-boundary}
$u^{-1}(\partial E) = \emptyset$.
\end{lemma}

\proof
To see why that is the case, suppose on the contrary that $x = u(z) \in \partial E$, which necessarily implies that $z \in S \setminus \partial S$. If $\mathrm{re}(\pi(x)) < -1$, we have $J = I_E$ and $K = 0$ near $(z,x)$. Then \eqref{eq:u-equation} reduces to the standard $I_E$-holomorphic map equation, so by the weak convexity of $\partial E$, an entire neighbourhood of $z$ must be mapped to $\partial E$. On the other hand, if $\mathrm{re}(\pi(x)) > -2$, one can use the isomorphism $\pi^{-1}(W_{-2}) \iso W_{-2} \times M = [-2,\infty) \times \R \times M$ to write locally $u(z) = (u_1(z),u_2(z))$. We still have $J = I_E = i \times I_M$. Moreover, $K(p+iq,x) = k_1(p)$, and hence $X_K = (X_{k_1},0)$. It follows that $u_2$ is an $I_M$-holomorphic map, and by applying weak convexity to $M$, we reach the same conclusion as before. We have now shown the set $u^{-1}(\partial E)$ is open and closed, which contradicts the fact that $\partial S \neq \emptyset$.
\qed

\begin{lemma} \label{th:u-semicircle}
There is an $R \gg 0$ such that for all solutions $u$, the following holds: if $\mathrm{re}(\pi(u(z))) < -1$, then $|\pi(u(z))| \leq R$.
\end{lemma}

\proof For a suitable $R$, we know that if $z$ is such that $\mathrm{re}(\pi(u(z))) < -1$ and $|\pi(u(z))| > R$, then $J = \tilde{J}$ and $K = \tilde{K} = 0$ locally, which implies that $v(z) = \pi(u(z))$ is a holomorphic function. Now apply the same degree argument as in Lemma \ref{th:v-semicircle}. \qed

\begin{lemma} \label{th:u-maximum}
There is a $C \gg 0$ such that for all solutions $u$, we have $\mathrm{re}(\pi(u)) \leq C$ everywhere.
\end{lemma}

\proof For a suitable $C$, we know that if $\mathrm{re}(\pi(u(z))) > C$, then $J = \tilde{J}$ and $K = \tilde{K}$ locally, which implies that $D\pi$ is $J_z$-holomorphic and $D\pi(X_K) = \partial_q \otimes \gamma$. Hence, $v = \pi(u(z))$ is locally a solution of $\bar\partial v = (\partial_q \otimes \gamma)^{0,1}$. Now apply the same maximum principle argument as in Lemma \ref{th:v-maximum}. \qed

Choose, for each boundary condition $L_z$, a function $F_z$ (locally constant in $z$) such that $\theta_E|L_z = dF_z$. With that, the energy-action calculus looks as follows. Define
\begin{align} \label{eq:u-energy}
& E^{\geom}(u) = \int_S {\textstyle\frac{1}{2}} \|du - X_K\|^2, \\
& E^{\top}(u) = \int_S u^*\omega_E - u^*dK, \\
& A(y_\zeta) = \Big(\int_{[0,1]} -u_\zeta^*\theta_E + K_\zeta(t,u_\zeta(t)) \,\mathit{dt}\Big)
+ F_{\zeta,1}(u_\zeta(1)) - F_{\zeta,0}(u_\zeta(0)),
\end{align}
where $F_{\zeta,k} = F_{\epsilon_{\zeta}(s,k)}$, and $u_\zeta$ is the flow line corresponding to $y_\zeta$. Then
\begin{align}
\label{eq:e-geom} & E^{\geom}(u) = E^{\top}(u) + \int_S u^*R, \\
\label{eq:e-a} & E^{\top}(u) = \sum_{\zeta \in \Sigma_{\mathrm{out}}} A(y_\zeta) - \sum_{\zeta \in \Sigma_{\mathrm{in}}} A(y_\zeta),
\end{align}
where the curvature, in local coordinates $z = s+it$ on $S$, is given by
\begin{equation} \label{eq:big-curvature}
R = \big( \partial_s K(\partial_t) - \partial_t K(\partial_s) + \{K(\partial_s),K(\partial_t)\} \big) \, \mathit{ds} \wedge \mathit{dt},
\end{equation}
with $\{\cdot,\cdot\}$ the Poisson bracket (the absence of a similar term in \eqref{eq:e-e} is explained by the fact that in that situation, both $k(\partial_s)$ and $k(\partial_t)$ would be functions of $p$, hence Poisson-commute). The curvature integrand in \eqref{eq:e-geom} is actually uniformly bounded, independently of the particular $u$. To see this, temporarily replace $K$ by $\tilde{K}$ as in \eqref{eq:k-prime}. Then, its curvature $\tilde{R}$ vanishes on $S \times \pi^{-1}(\C \setminus W_{-1})$, and is invariant under translation in $q$-direction on $S \times \pi^{-1}(W_{-1})$, hence is bounded on any subset of the form $S \times \pi^{-1}(\C \setminus W_C)$. $R$ differs from $\tilde{R}$ only on a compact subset, which shows that it has the same boundedness property. Finally, the integrand vanishes on the strip-like ends, which concludes the argument. The consequence of this and \eqref{eq:e-a} is this:

\begin{lemma} \label{th:energy-bound}
There is a uniform bound on $E^{\geom}(u)$ over all solutions $u$. \qed
\end{lemma}

\begin{proposition} \label{th:u-bounded}
All solutions $u$ remain within a compact subset of $E$.
\end{proposition}

\proof From Lemmas \ref{th:u-semicircle} and \ref{th:u-maximum}, we know a priori that any $u$ remains in the region
\begin{equation} \label{eq:e-region}
\{|\pi(x)| \leq R\} \cup \{-1 \leq \mathrm{re}(\pi(x)) \leq C\} \subset E
\end{equation}
for suitable constants $R,C$. Due to the conditions on $E$ and $J$, that region has bounded geometry (see for instance \cite{audin-lafontaine} for a consideration of that notion in the context of pseudo-holomorphic curve theory), and moreover, $K$ and all its derivatives are bounded on it. In principle, knowing this is already sufficient to apply a suitable form of Gromov compactness, which implies the desired result. Nevertheless, it is maybe appropriate to explain more explicitly how the given problem can be reduced to a more standard form.

Suppose that we have a sequence of solutions $u_j$ such that $\|\pi(u_j)\|_\infty$, hence also $\|\mathrm{im}(\pi(u_j))\|_{\infty}$, goes to infinity. Then, $\|du_j\|_\infty$ necessarily goes to infinity as well. Take the points $z_j$ where $\|du_j\|$ reaches its maximum, for $j \gg 0$. If $u_j(z_j)$ has a bounded subsequence, we can apply the standard rescaling argument (taking into account what was said above concerning bounded geometry) to obtain a nonconstant $J_z$-holomorphic map $\C P^1 \rightarrow E$ for some $z \in S$, or else a nonconstant $J_z$-holomorphic disc $(D,\partial D) \rightarrow (E,L_z)$ for some $z \in \partial S$, contradicting exactness.

Now consider the other situation, where $|\mathrm{im}(\pi(u_j(z_j)))|$ goes to infinity. We work in local coordinates in $S$ exactly as in the corresponding part of Proposition \ref{th:v-bounded}, and proceed similarly: fix some $\rho$, and take a closed disc of radius $\rho/\|du_j\|_\infty$ around $z_j$. Assuming that $j$ is large, this disc lies in the interior of $S$, because the distance of $u_j(z_j)$ to $u_j(\partial S)$ is greater than $\rho$. Moreover, the image of the disc under $u_j$ is contained in $\pi^{-1}(W_{-1})$, because the distance of $u_j(z_j)$ to the other part $\{|\pi(x)| \leq R\}$ of \eqref{eq:e-region} is again greater than $\rho$. We restrict $u_j$ to this disc, rescale the domain by $\|du_j\|_{\infty}$, and using the trivialization $\pi^{-1}(W_{-1}) \iso W_{-1} \times M$, subtract a constant which is the imaginary part of $\pi(u_j(z_j))$. The last-mentioned step works because the image of $u_j$ will lie in a region where $J = \tilde{J}$ is invariant under translation in imaginary direction. The limit of a subsequence of these translated maps is a pseudo-holomorphic map $\tilde{u}$ defined on a slightly smaller disc, satisfying $\tilde{u}(0) \in \pi^{-1}([-1,C])$ and $\|d\tilde{u}(0)\| = \|d\tilde{u}\|_{\infty} = 1$, and whose energy is still bounded by the original upper bound on geometric energies. Repeat this process with subsequently larger values of $\rho$, passing to further subsequences. In the end, one again gets a non-constant pseudo-holomorphic map $\C P^1 \rightarrow W_{-1} \subset E$, hence a contradiction. \qed

\section{Algebraic structures from Lefschetz thimbles\label{sec:lefschetz-2}}

\subsection{Floer cohomology\label{subsec:pm-floer}}
We begin by considering the cohomology level theory. This discussion is kept brief, since the overall framework is  familiar, and anyway most of it will be superseded by the cochain level constructions to be introduced in more detail later on.

Fix an exact symplectic Lefschetz fibration $\pi: E \rightarrow \C$, a pair of vanishing paths $(l_0,l_1)$ with ordinates $(o_0,o_1)$, and a real number $c$ satisfying \eqref{eq:nondeg}. Given that, one can define Floer cohomology $\mathit{HF}(L_0,L_1)^c$ for the associated Lefschetz thimbles, as follows. Choose a $c$-displacement function $H$, with the additional property that $\phi_H(L_0)$ intersects $L_1$ transversally. We then define the Floer cochain space to be the $\K$-vector space freely generated by points of $\phi_H(L_0) \cap L_1$. Fix a function $g \in \smooth([0,1],\R)$ such that $\int_0^1 g(t)\, \mathit{dt} = 1$. Additionally, choose a family $J = (J_t)_{0 \leq t \leq 1}$ in $\JJ(E)$. The Floer equation is
\begin{equation}
\label{eq:floer-lefschetz}
\left\{
\begin{aligned}
& u: \R \times [0,1] \longrightarrow E, \\
& u(s,0) \in L_0, \;\; u(s,1) \in L_1, \\
& \partial_s u + J_t (\partial_t u - g(t) X_H) = 0.
\end{aligned}
\right.
\end{equation}
For $c \neq 0$, this is a special case of \eqref{eq:u0-equation}, with $S = \R \times [0,1]$, worldsheet structure $\gamma = c g(t) \, \mathit{dt}$, and Hamiltonian $c^{-1} H$. However, it is more natural to think of it as an instance of \eqref{eq:u-equation} with $K = H g(t)\,\mathit{dt}$, and that covers the case $c = 0$ as well. Given two intersection points, a generic choice of $J$ ensures that we have a well-defined count of solutions of our equation with suitable asymptotics. The ingredients that enter into the proof of this fact are: transversality, which is not difficult since we have complete freedom to perturb the almost complex structures on any compact subset of $E \setminus \partial E$; Proposition \ref{th:u-bounded}, which gives a $C^0$ bound on solutions; and standard compactness arguments. The same techniques, together with Lemma \ref{th:no-boundary} which excludes families of solutions escaping to the boundary, shows that the resulting Floer differential squares to zero. Its cohomology is $\mathit{HF}(L_0,L_1)^c$.

\begin{example} \label{th:hf-basis}
Let $(l_1,\dots,l_m)$ be a simply intersecting basis of vanishing paths, with associated thimbles $(L_1,\dots,L_m)$. Take a pair $(L_i,L_j)$ with $i<j$, and suppose that $c < o_j - o_i$ is less than the difference between their ordinates. Because of the special position of the paths, the Hamiltonian can be chosen in such a way that $\phi_H(L_i) \cap L_j = \emptyset$, so $\mathit{HF}(L_i,L_j)^c$ vanishes. The same holds if $i>j$ and $c > o_j - o_i$.
\end{example}

For all $c^+ \geq c^-$ such that both sides are defined, we have canonical maps
\begin{equation} \label{eq:increase-nu}
\mathit{HF}(L_0,L_1)^{c^-} \longrightarrow \mathit{HF}(L_0,L_1)^{c^+}.
\end{equation}
These compose in the obvious way, and reduce to the identity when $c^+ = c^-$. The definition is in terms of continuation map equations, which are of the form \eqref{eq:u-equation} on the following worldsheet:
\begin{equation} \label{eq:cont-pm}
\parbox{28em}{
$S = \R \times [0,1]$, with a one-form $\gamma \in \Omega^1(S)$ which satisfies $\gamma = c^\pm g^\pm(t) \mathit{dt}$ for $\pm s \ll 0$ (this can exist only for $c^+ \geq c^-$, since otherwise \eqref{eq:sub-closed} would be violated). One similarly chooses the inhomogeneous term $K$ such that $K = H^\pm g^\pm(t) \mathit{dt}$ for $\pm s \ll 0$.}
\end{equation}

Suppose that we have $c^\pm$ which both lie on the same side of $o_1 - o_0$ (both smaller or both larger). By carefully choosing $H^{\pm}$ and $K$, one can achieve that the associated Floer cochain spaces coincide, and that the cochain level map underlying \eqref{eq:increase-nu} is an isomorphism. We omit the details, which are fairly standard (one arranges that the map is compatible with the action filtration, and that the only terms which strictly preserve the action are constant solutions). The outcome is that there are really only two essentially different groups
\begin{equation} \label{eq:plus-cohomology}
\begin{aligned}
&
\mathit{HF}(L_0,L_1)^- = \mathit{HF}(L_0,L_1)^c \quad \text{for $c < o_1 - o_0$,} \\
&
\mathit{HF}(L_0,L_1)^+ = \mathit{HF}(L_0,L_1)^c \quad \text{for $c > o_1 - o_0$.}
\end{aligned}
\end{equation}

Finally, there are canonical Poincar{\'e} duality isomorphisms
\begin{equation} \label{eq:c-duality}
\mathit{HF}(L_0,L_1)^c \iso (\mathit{HF}(L_1,L_0)^{-c})^\vee
\end{equation}
which allow one to reduce the framework further by using only $\mathit{HF}(L_0,L_1)^+$. These groups fit into an asymmetric TQFT framework whose structure is dictated by the notion of worldsheet. As a standard part of this, they carry associative products, as well as unit elements in $\mathit{HF}(L,L)^+$. The latter are defined in terms of a surface $S$ which is a once-punctured disc with a single output end, carrying a one-form $\gamma$ such that $\gamma = c g(t) \mathit{dt}$ over the end, with $c>0$. Unlike the more familiar case of Floer cohomology for compact manifolds, there is no canonical linear map from $\mathit{HF}(L,L)^+$ to scalars, and no pairing between Floer cohomology groups $\mathit{HF}(L_0,L_1)^+$ and $\mathit{HF}(L_1,L_0)^+$. However, the nontrivial map \eqref{eq:increase-nu} translates into a diagonal class
\begin{equation} \label{eq:diagonal}
D_{L_0,L_1} \in \mathit{HF}(L_0,L_1)^+ \otimes \mathit{HF}(L_1,L_0)^+,
\end{equation}
satisfying equations which parallel \eqref{eq:gr-sym}, \eqref{eq:both-sides}. In fact, for any finite family of vanishing cycles $\{L_1,\dots,L_m\}$, the direct sum $A = \bigoplus_{i,j} \mathit{HF}(L_i,L_j)^+$ is naturally a boundary algebra in the sense of Section \ref{sec:boundary-alg}, except for the missing grading.

\begin{example} \label{th:single-l}
Consider a single $L$. Then $\mathit{HF}(L,L)^+ \iso H^*(L;\K)$, $\mathit{HF}(L,L)^- \iso H^*_{\mathrm{cpt}}(L;\K)$, where the latter group is compactly supported cohomology. The homomorphism \eqref{eq:increase-nu} is the standard map from compactly supported cohomology to the ordinary one (or equivalently, $D_{L,L}$ is the Poincar{\'e} dual of the diagonal as in Example \ref{th:cohomology-boundary}). Of course, since $L$ is contractible, this map vanishes.
\end{example}

\subsection{An $A_\infty$-algebra\label{subsec:t-algebra}}
From now on, we will work with a fixed simply intersecting basis of vanishing paths $(l_1,\dots,l_m)$. We will assume that the associated vanishing cycles $(V_1,\dots,V_m)$ in $M$ are in general position in the sense of \eqref{eq:generic-1}. For fixed paths, this assumption is generic within a suitable class of perturbations of the symplectic form on $E$, but here we assume that it holds for the original $\omega_E$.

\begin{remark} \label{th:transverse-thimbles}
The assumption \eqref{eq:generic-1} implies the corresponding general position statement for the Lefschetz thimbles $(L_1,\dots,L_m)$, without requiring any further condition on the vanishing paths (other than those imposed by the notion of simply intersecting basis, see Definition \ref{def:vanishing-1}). Namely, any two distinct $l_i$ intersect transversally, and all intersection points lie in $\pi^{-1}(U_{0,1})$, which means that the intersections of the vanishing thimbles over such points are governed by the intersections of the vanishing cycles, hence are transverse. Similarly, while it is possible for three distinct vanishing paths to meet at the same point, the associated three vanishing cycles will not have a common point in that fibre, hence $L_i \cap L_j \cap L_k = \emptyset$ for $i<j<k$.
\end{remark}

Let $\mathit{CF}(L_i,L_j)$, for $i<j$, be the $\K$-vector space freely generated by $L_i \cap L_j$. We turn
\begin{equation} \label{eq:l-bar-a}
\bar\TT = \bigoplus_{i<j} \mathit{CF}(L_i,L_j)
\end{equation}
into an $A_\infty$-algebra over $R = \K^m$, as follows. For each $(i,j)$, choose a generic family of almost complex structures $J_{i,j}^1 = (J_{i,j,t}^1)_{0 \leq t \leq 1}$ in $\JJ(E)$, and consider the associated pseudo-holomorphic strip equation, which gives rise to a differential $\mu^1$. This is just the special case of \eqref{eq:floer-lefschetz} where $c$ and $H$ vanish. To define the higher order operations
\begin{equation} \label{eq:lefschetz-a}
\mu^d: \mathit{CF}(L_{i_{d-1}},L_{i_d}) \otimes \cdots \otimes \mathit{CF}(L_{i_0},L_{i_1}) \longrightarrow \mathit{CF}(L_{i_0},L_{i_d})
\end{equation}
we consider $(d+1)$-punctured discs, with boundary conditions $(l_{i_0},\dots,l_{i_d})$ for $i_0 < \cdots < i_d$. On the universal family of such discs, we choose a family $J_{i_0,\dots,i_d}^d$ of almost complex structures taken from $\JJ(E)$, and subject to standard consistency conditions. The one-forms on these surfaces are supposed to be zero, and consequently we consider the pseudo-holomorphic map equation \eqref{eq:u0-equation} with zero inhomogeneous term. Counting points in the resulting moduli spaces $\scrM^d(y_0,\dots,y_d)$ leads to numbers $n^d(y_0,\dots,y_d) \in \K$ which are the coefficients of \eqref{eq:lefschetz-a}. We write
\begin{center}
\begin{tabular}{|l|l|}
\hline
$CF_{\bar\TT}(L_i,L_j)$ & Floer cochain spaces, defined for all $i < j$ \\
$\mu^d_{\bar\TT}$ & $A_\infty$-operations \\
$\scrS^d_{\bar{\TT}} \rightarrow \scrR^d_{\bar{\TT}}$ & moduli space of $(d+1)$-punctured discs \\
$J^d_{\bar{\TT},i_0,\dots,i_d}$ & families of almost complex structures \\
$\scrM_{\bar{\TT}}^d(y_0,\dots,y_d)$ & moduli spaces of pseudo-holomorphic maps \\
$n_{\bar{\TT}}^d(y_0,\dots,y_d)$ & numbers counting points in the spaces above \\
\hline
\end{tabular}
\end{center}

\subsection{A family of bimodules\label{subsec:uu-modules}}
If $(o_1,\dots,o_m)$ are the ordinates of our vanishing paths, fix any number $c$ such that $c \neq o_i - o_j$ for all $i,j$ (in particular, $c \neq 0$). For each such $c$ we will define an $A_\infty$-bimodule $\UU^c$ over $\bar\TT$. Choose a $c$-displacement function $H$. We require transversality of the intersections $\phi_H(L_i) \cap L_j$, and genericity of triple intersections, in the same sense as in \eqref{eq:b-general-position}. Write $\mathit{CF}(L_i,L_j)$ for the vector space freely generated by points in $\phi_H(L_i) \cap L_j$, and set
\begin{equation}
\UU^c = \bigoplus_{i,j} \mathit{CF}(L_i,L_j).
\end{equation}
Fix a function $g \in \smooth([0,1],\R)$ vanishing near $t = 0$ and $t = 1$, and such that $\int_0^1 g(t)\, \mathit{dt} = 1$. We consider surfaces $S$ as in \eqref{eq:punctured-strip}, with the worldsheet structure $\gamma = c g(t) \, \mathit{dt}$, and carrying the inhomogeneous term $K = H g(t) \,\mathit{dt}$. In the simplest case, which is $S = \R \times [0,1]$ equipped with boundary conditions $(l_i,l_j)$, we choose a family $J^{0|1|0}_{i,j} = (J_{i,j,t}^{0|1|0})$ in $\JJ(E)$. The moduli space $\scrM^{0|1|0}(y_0,y_1)$ of solutions of the associated equation \eqref{eq:floer-lefschetz} yields numbers $n^{0|1|0}(y_0,y_1)$ which define the differential $\mu^{0|1|0}$ on $\UU^c$. This is just one case of the Floer differential from Section \ref{subsec:pm-floer}. In particular, in terms of \eqref{eq:plus-cohomology} we have
\begin{equation} \label{eq:ppc-cohomology}
H(\UU^c) = \bigoplus_{i,j} \mathit{HF}(L_i,L_j)^c \iso
\bigoplus_{\substack{i \leq j \\ o_j - o_i < c}} \mathit{HF}(L_i,L_j)^+ \;\oplus
\bigoplus_{\substack{i \geq j \\ o_j - o_i > c}} \mathit{HF}(L_i,L_j)^-.
\end{equation}
The restrictions to $i \leq j$ in the first sum, and to $i \geq j$ in the second sum, follow from the fact that we are considering a simply intersecting basis, see Example \ref{th:hf-basis}.

More generally, consider surfaces $S$ with additional $p+q$ punctures as in \eqref{eq:punctured-strip}, where the boundary components are labeled by $(l_{i_0},\dots,l_{i_{p+q+1}})$ for some $i_0 < \cdots < i_p$ and $i_{p+1} < \cdots < i_{p+q+1}$. We equip each such $S$ with the same $\gamma$ and $K$ as before, and with a family $J = (J_z)$ of almost complex structures in $\JJ(E)$. For points $z = (s,t)$ with $\pm s \ll 0$, this should agree with the previously defined $J^{0|1|0}_{i_0,i_{p+q}}$ and $J^{0|1|0}_{i_p,i_{p+1}}$, respectively, while near all the other punctures it is determined by $J^1_{\bar\TT,i_k,i_{k+1}}$. More precisely, we want to choose a family $J^{q|1|p}_{i_0,\dots,i_{p+q+1}}$ on the universal family of surfaces $\scrS^{q|1|p} \rightarrow \scrR^{q|1|p}$, and this is subject to suitable consistency conditions. By counting points in the associated moduli spaces $\scrM^{q|1|p}(y_0,\dots,y_{p+q+1})$, we get numbers $n^{q|1|p}(y_0,\dots,y_{p+q+1})$ which form the coefficients of the structure maps
\begin{multline} \label{eq:b-lefschetz-structure}
\qquad \mu^{q|1|p}: \mathit{CF}_{\bar\TT}(L_{i_{p+q}},L_{i_{p+q+1}}) \otimes \cdots \otimes \mathit{CF}(L_{i_p},L_{i_{p+1}}) \otimes \cdots \\ \cdots \otimes \mathit{CF}_{\bar\TT}(L_{i_0},L_{i_1}) \longrightarrow \mathit{CF}(L_{i_0},L_{i_{p+q+1}}). \qquad
\end{multline}
As usual, we retrospectively change the relevant notation:
\begin{center}
\begin{tabular}{|l|l|}
\hline
$H_{\UU^c}$ & Hamiltonian \\
$g_{\UU^c}$ & cut-off function \\
$\mathit{CF}_{\UU^c}(L_i,L_j)$ & $H_{\UU^c}$-perturbed Floer cochain spaces, \\
& \hspace{2em} defined for all $i,j$ \\
$\mu^{q|1|p}_{\UU^c}$ & $A_\infty$-bimodule operations on $\UU^c$ \\
$\scrS^{q|1|p}_{\UU^c} \rightarrow \scrR^{q|1|p}_{\UU^c}$ & moduli space of $(p+q+2)$-punctured discs \\
$J^{q|1|p}_{\UU^c,i_0,\dots,i_{p+q+1}}$ & families of almost complex structures \\
$\scrM_{\UU^c}^{q|1|p}(y_0,\dots,y_{p+q+1})$ & moduli spaces of perturbed \\
& \hspace{2em} pseudo-holomorphic maps \\
$n_{\UU^c}^{q|1|p}(y_0,\dots,y_{p+q+1})$ & numbers counting points in the spaces above \\
\hline
\end{tabular}
\end{center}

\subsection{Bimodule homomorphisms}
For any $c^- \leq c^+$ such that $\UU^{c^\pm}$ are defined, one has a natural bimodule map
\begin{equation} \label{eq:gamma-map}
\Gamma: \UU^{c^-} \longrightarrow \UU^{c^+}.
\end{equation}
To construct that, start with $\R \times [0,1]$ equipped with a one-form $\gamma$ and inhomogeneous term $K$ as in \eqref{eq:cont-pm}, whose behaviour as $s \rightarrow \pm\infty$ is determined by the previously chosen data for $\UU^{c^{\pm}}$, and with the additional assumption that both should vanish near $\R \times \{0,1\}$. Consider surfaces $S$ as in \eqref{eq:punctured-strip}, but which now come with a fixed isomorphism of their partial compactification with $\R \times [0,1]$. This means that we no longer identify surfaces differing by a common translation of the marked points $\zeta_i$ in $s$-direction, and allows us to restrict $\gamma$ and $K$ to $S$ uniquely. On each such surface, choose a family $J = (J_z)$ of almost complex structures in $\JJ(E)$ as follows. For points $z = (s,t)$ with $\pm s \ll 0$, it should agree with the previously defined almost complex structures for $\UU^{c^\pm}$, while near all the other punctures it is determined by those for $\bar\TT$. Overall, this leads to a family $J^{q|1|p}$ parametrized by a suitably defined universal surface $\scrS^{q|1|p} \rightarrow \scrR^{q|1|p}$. By counting points in the associated moduli spaces $\scrM^{q|1|p}(y_0,\dots,y_{p+q+1})$, we get numbers $n^{q|1|p}(y_0,\dots,y_{p+q+1})$ which form the coefficients of $\Gamma^{q|1|p}$. At first sight this may look identical to the previous definition of the bimodules themselves. However, because we are no longer dividing by translation, the moduli spaces themselves are different: $\scrR^{q|1|p}$ is of dimension $p+q$, and its compactification $\bar\scrR^{q|1|p}$ has more strata. Figure \ref{fig:bimodule-map} shows the four basic kinds of codimension one strata in the compactification (the vertical dashed lines symbolically indicate the breaking of translation invariance, which means that these are the principal components in the compactification process, the others being thought of as bubbles). These are responsible for the terms in the bimodule homomorphism formula
\begin{equation} \label{eq:bimodule-hom}
\begin{aligned}
& \sum_{r,s} \Gamma^{q-s|1|p-r}(c_{p+q},\dots,c_{p+s+1},\mu^{s|1|r}_{\UU^{c^-}}(t_{p+s},\dots,t_{p+1},u, \\[-1em]
& \hspace{18em} t_p,\dots,t_{p-r+1}),t_{p-r}, \dots,c_1) \\
& + \sum_{r,s} \mu^{q-s|1|p-r}_{\UU^{c^+}}(c_{p+q},\dots,c_{p+s+1},\Gamma^{s|1|r}(t_{p+s},\dots,t_{p+1},u, \\[-1em]
& \hspace{18em} t_p,\dots,t_{p-r+1}),t_{p-r}, \dots,c_1)
\\ = & \sum_{j,s} \Gamma^{q-s+1|1|p}(c_{p+q},\dots,\mu_{\bar\TT}^s(t_{p+j+s},\dots,t_{p+j+1}),t_{p+j},\dots,
t_{p+1},u,t_p,\dots,t_1) \\ + & \sum_{i,r} \Gamma^{q|1|p-r+1}(t_{p+q},\dots,t_{p+1},u,t_p,\dots,t_{i+r+1},
\mu_{\bar\TT}^r(t_{i+r},\dots,t_{i+1}),t_i,\dots,t_1).
\end{aligned}
\end{equation}
\begin{figure}
\begin{center}
\begin{picture}(0,0)%
\includegraphics{bimodule-map.pstex}%
\end{picture}%
\setlength{\unitlength}{3947sp}%
\begingroup\makeatletter\ifx\SetFigFont\undefined%
\gdef\SetFigFont#1#2#3#4#5{%
  \reset@font\fontsize{#1}{#2pt}%
  \fontfamily{#3}\fontseries{#4}\fontshape{#5}%
  \selectfont}%
\fi\endgroup%
\begin{picture}(5724,5644)(-2186,-3408)
\end{picture}%
\caption{\label{fig:bimodule-map}}
\end{center}
\end{figure}

\begin{remark}
More abstractly, one can think of $\scrR^{q|1|p}$ as a moduli space of discs with $p+q+2$ boundary points, and with an additional marked path in its interior (a geodesic for the hyperbolic metric on the open disc). To get to the picture above, one first removes the points marked $0$ and $p+1$ and identifies their complement with $\R \times [0,1]$, in such a way that the marked path corresponds to $\R \times (0,1)$.

This may remind the reader of the discussion of the spaces $\scrR^d_\Phi$ in Remark \ref{th:moduli1} (the same spaces will come up again in Section \ref{subsec:hh2} below). Obviously, those spaces as well as the ones considered here are topologically open balls, but the differences become clear once one passes to the appropriate compactifications; this reflects the structural difference between \eqref{eq:hochschild-d} and \eqref{eq:bimodule-hom}. For instance, looking at the two-dimensional moduli spaces, $\bar\scrR^2_\Phi$ is a hexagon (drawn for instance in \cite[p.\ 255]{seidel04}), whereas here we have two different two-dimensional moduli spaces $\bar\scrR^{1|1|1}$ and $\bar\scrR^{2|1|0}$ (which is isomorphic to $\bar\scrR^{0|1|2}$), of which the first is a hexagon while the second is a pentagon. In one dimension higher, $\bar\scrR^3_\Phi$ is a polytope whose boundary faces are four hexagons, four pentagons, and four squares; whereas $\partial\bar\scrR^{2|1|1}$ (as well as $\partial\bar\scrR^{1|1|2}$) consists of three hexagons, four pentagons, and four squares; and $\partial\bar\scrR^{3|1|0}$ (as well as $\partial\bar\scrR^{0|1|3}$) of six pentagons and three squares.
\end{remark}

Here's the summary of notation for future reference (where the $c^\pm$ have been suppressed, as in \eqref{eq:gamma-map}, to keep things reasonably simple):
\begin{center}
\begin{tabular}{|l|l|}
\hline
$\gamma_{\Gamma}$ & one-form on $\R \times [0,1]$ \\
$K_{\Gamma}$ & inhomogeneous term \\
$\scrS^{q|1|p}_{\Gamma} \rightarrow \scrR^{q|1|p}_{\Gamma}$ & configuration space of $p+q$ points on $\R \times \{0,1\}$ \\
$J^{q|1|p}_{\Gamma,i_0,\dots,i_{p+q+1}}$ & families of almost complex structures \\
$\scrM_{\Gamma}^{q|1|p}(y_0,\dots,y_{p+q+1})$ & moduli spaces of perturbed \\
& \hspace{2em} pseudo-holomorphic maps \\
$n_{\Gamma}^{q|1|p}(y_0,\dots,y_{p+q+1})$ & numbers counting points in the spaces above \\
\hline
\end{tabular}
\end{center}

\begin{remark} \label{th:general-k}
In principle, one could simplify the setup a little by assuming that $H^\pm = c^\pm H$ are both multiples of the same $1$-displacement function, and then choosing $K = H \gamma$. However, the freedom to choose the two functions independently will be crucial for our applications, which is why we're working with the wider class of $K$ (in contrast to the stylistic choices made at other places in the paper). This has another advantage, this time on the technical side, namely that transversality is quite straightforward if we allow ourselves to vary $K$ as well as the almost complex structures.
\end{remark}

As should be obvious from the definition, the map between cohomology groups \eqref{eq:ppc-cohomology} induced by $\Gamma$ is the direct sum of \eqref{eq:increase-nu}. As a consequence, if the interval $[c^-,c^+]$ contains none of the $o_j - o_i$, $\Gamma$ is a quasi-isomorphism, which shows that the $\UU^c$ are divided into finitely many quasi-isomorphism classes. On the other hand, if we take $c^- \ll 0$ and $c^+ \gg 0$, it follows from \eqref{eq:ppc-cohomology} and Examples \ref{th:hf-basis}, \ref{th:single-l} that the cohomology level map $H(\Gamma)$ vanishes (but this of course does not mean that $\Gamma$ itself is trivial, not even up to homotopy). We also want to mention briefly some additional properties. The maps \eqref{eq:gamma-map} are independent of all choices up to homotopy. In the same sense, they are compatible with composition. Finally, if $c^- = c^+$, then $\Gamma$ is homotopic to the identity. The proofs of these statements involve parametrized versions of the moduli spaces in our construction, together with appropriate gluing theorems, much as in the classical definition of continuation maps \cite{salamon-zehnder92}.

\subsection{Hochschild cocycles\label{subsec:hh2}}
Let $\TT = R \oplus \bar{\TT}$ be the result of adjoining units to $\bar{\TT}$, as in \eqref{eq:adjoin-units}. Fix a $c > 0$ such that $\UU^c$ is well-defined. By looking at \eqref{eq:ppc-cohomology} and Example \ref{th:single-l}, one sees that there is a natural inclusion
\begin{equation} \label{eq:include-tt}
H(\TT) = \bigoplus_{i \leq j} \mathit{HF}(L_i,L_j)^+ \subset H(\UU^c).
\end{equation}
To raise this basic insight to the cochain level, we will construct a Hochschild cocycle $\bar\Pi \in \mathit{CC}(\bar\TT,\UU^c)$, in parallel with the discussion in Section \ref{subsec:hh}. Namely, on the upper half-plane $\R \times \R^+$, choose a worldsheet structure $\gamma$ and an inhomogeneous term $K$, with the following property. Both should vanish near the boundary, and their pullback to the strip-like end \eqref{eq:half-plane-end} should be equal to $c g_{\UU^c} \mathit{dt}$ and $H_{\UU^c} \mathit{dt}$, respectively (because of \eqref{eq:sub-closed}, this would be impossible for negative $c$).

The general class of surfaces which we will consider are $(d+1)$-punctured discs with an additional interior marked point, realized as in \eqref{eq:punctured-plane} as version of the upper half-plane with additional punctures, and boundary conditions $(l_{i_0},\dots,l_{i_d})$ for some $i_0 < \cdots < i_d$. We equip them with the restrictions of $\gamma$ and $K$, as well as with almost complex structures $J = (J_z)$ whose restriction to the end $\mathrm{im}(z) \rightarrow \infty$ is given by $J_{\UU^c,i_0,i_d}^{0|1|0}$, and which near the other punctures are given by $J_{\bar\TT,i_k,i_{k+1}}^1$. As usual, we actually need one family $J_{i_0,\dots,i_d}^d$ on the universal family of surfaces $\scrS^d \rightarrow \scrR^d$, for each $d \geq 0$. Having that, we build moduli spaces $\scrM^d(y_0,\dots,y_d)$ and get numbers $n^d(y_0,\dots,y_d)$, which are the coefficients of $\bar\Pi^d$. The adjusted notation is:
\begin{center}
\begin{tabular}{|l|l|}
\hline
$\gamma_{\bar{\Pi}}$ & one-form on the upper half plane \\
$K_{\bar{\Pi}}$ & inhomogeneous term \\
$\scrS^d_{\bar{\Pi}} \rightarrow \scrR^d_{\bar{\Pi}}$ & moduli space of $d+1$-punctured discs \\
& \hspace{1em} with an additional interior marked point \\
$J^d_{\bar{\Pi},i_0,\dots,i_{p+q+1}}$ & family of almost complex structures \\
$\scrM_{\bar{\Pi}}^d(y_0,\dots,y_d)$ & moduli spaces of pseudo-holomorphic maps \\
$n_{\bar{\Pi}}^d(y_0,\dots,y_d)$ & numbers counting points in the spaces above \\
\hline
\end{tabular}
\end{center}

The Hochschild cohomology class $[\bar\Pi]$ is independent of all choices. $[\bar\Pi^0] \in H(\UU^c)$ is the sum of unit elements in $\mathit{HF}(L_i,L_i)^+ \iso H^*(L_i;\K)$ for all $i$. Consider the bimodule map $\Pi = X(\bar\Pi): \TT \rightarrow \UU^c$ constructed from $\bar\Pi$ using \eqref{eq:functorial-hh}. On the cohomological level, is given by multiplying with $[\bar\Pi^0]$. Hence, its components are the identity maps on $\mathit{HF}(L_i,L_j)^+$ for all $i \leq j$ such that $c > o_j - o_i$. In particular:

\begin{cor} \label{th:plus-u-quasi-iso}
If $c > o_j - o_i$ for all $i,j$, then $\Pi: \TT \rightarrow \UU^c$ is a $\bar\TT$-bimodule quasi-isomorphism. \qed
\end{cor}

There is also a dual picture, which is precisely parallel to that in Section \ref{subsec:dual-diagonal}. One uses the same Riemann surfaces, but where the end $\mathrm{im}(z) \rightarrow \infty$ is considered as an input, and the condition \eqref{eq:sub-closed} correspondingly requires that $c \leq 0$. For each such $c$ such that $\UU^c$ is defined (this of course excludes $c = 0$), we then get a cocycle $\bar\Xi \in \mathit{CC}(\bar\TT,(\UU^c)^\vee)$. The leading term $[\bar\Xi^0] \in H(\UU^c)^-$ is the direct sum of integration maps $\mathit{HF}(L_i,L_i)^- \iso H^*_{\mathrm{cpt}}(L_i;\K) \rightarrow \K$ (compare Example \ref{th:single-l}). When composed with multiplication, these integration maps yield nondegenerate pairings
\begin{equation}
\mathit{HF}(L_j,L_i)^- \otimes \mathit{HF}(L_i,L_j)^+ \longrightarrow \mathit{HF}(L_i,L_i)^- \longrightarrow \K.
\end{equation}
The bimodule homomorphism $\Xi = Y(\bar\Xi)$, on the cohomology level, consists of the associated isomorphisms $\mathit{HF}(L_i,L_j)^- \rightarrow (\mathit{HF}(L_j,L_i)^+)^\vee$ for all $i \geq j$ such that $c < o_j - o_i$ (these isomorphisms were already mentioned in \eqref{eq:c-duality} as part of the general Floer cohomology setup for Lefschetz thimbles). Therefore:

\begin{cor} \label{th:minus-u-quasi-iso}
If $c < o_j - o_i$ for all $i,j$, then $\Xi: \UU^c \rightarrow \TT^\vee$ is a $\bar\TT$-bimodule quasi-isomorphism. \qed
\end{cor}

In conclusion, we have arrived at the following picture. Fix constants $c^\pm$ such that $c^- < o_j - o_i < c^+$ for all $i,j$. Then we have two $\bar\TT$-bimodules $\UU^{c^+}$ and $\UU^{c^-}$, which are quasi-isomorphic to the diagonal bimodule $\TT$ and its dual $\TT^\vee$, respectively. Using these quasi-isomorphisms, one can turn $\Gamma$ into a bimodule homomorphism
\begin{equation}
\Pi^{-1} \circ \Gamma \circ \Xi^{-1}: \TT^\vee \longrightarrow \TT,
\end{equation}
which is well-defined up to homotopy. This should be considered as a natural addition to the $A_\infty$-algebra structure of $\TT$, arising from the specific geometry of the Lefschetz thimbles.

\section{Comparison results\label{sec:comparison}}

\subsection{Outline}
Fix an exact symplectic Lefschetz fibration $\pi: E \rightarrow \C$, as well as a simply intersecting basis of vanishing paths $(l_1,\dots,l_m)$. Consider on one hand, the associated vanishing cycles $(V_1,\dots,V_m)$ in the fibre $M$, assumed to be in general position in the sense of \eqref{eq:generic-1}; and on the other hand, the Lefschetz thimbles $(L_1,\dots,L_m)$ in $E$. The algebraic structures obtained from either datum (in Sections \ref{sec:fibre} and \ref{sec:lefschetz-2}, respectively) run partially in parallel:
\begin{center}
\begin{tabular}{l|l}
fibre & total space \\ \hline \\[-1em]
$A_\infty$-algebra $\bar\A$ & $A_\infty$-algebra $\bar\TT$ \\
bimodules $\B^\pm$ & bimodules $\UU^{\pm}$ \\
morphism $\Delta: \B^- \rightarrow \B^+$ & morphism $\Gamma: \UU^- \rightarrow \UU^+$
\end{tabular}
\end{center}
where in the right column, $\UU^- = \UU^{c^-}$ for $c^- \ll 0$, and similarly $\UU^+ = \UU^{c^+}$ for $c^+ \gg 0$. It is not difficult to see that
\begin{equation}
\mathit{HF}(L_i,L_j)^+ \iso \mathit{HF}(V_i,V_j) \quad\text{for $i<j$.}
\end{equation}
Here, the left hand side is Floer cohomology in $E$ in the sense of \eqref{eq:plus-cohomology}, and the right hand side is the standard Floer cohomology in $M$. This means that $H(\bar\TT) \iso H(\bar\A)$ as $R$-bimodules. Further inspection shows that $\bar\A$ and $\bar\TT$ are quasi-isomorphic as $A_\infty$-algebras, and actually identical if the choices in both definitions are correlated appropriately. Concerning the bimodules, we know from Propositions \ref{th:plus-quasi-iso}, \ref{th:minus-quasi-iso} that there are quasi-isomorphisms $\B^+ \htp \A$, $\B^- \htp \A^\vee$. Similarly, we know that $\UU^+ \htp \TT$, $\UU^- \htp \TT^\vee$ by Corollaries \ref{th:plus-u-quasi-iso}, \ref{th:minus-u-quasi-iso}.

This leaves the bimodule homomorphisms, which are not covered by the results mentioned so far. Assume that things have been set up in such a way that $\bar\A$ and $\bar\TT$ can be identified directly. By further careful choices, we will single out a sub-bimodule $\VV \subset \UU^+$. This has the following properties with respect to the map $\Gamma$, as well as the Hochschild cocycles $\bar\Phi$, $\bar\Pi$ from Sections \ref{subsec:hh} and \ref{subsec:hh2}, respectively.

\begin{theorem} \label{th:comparison}
(i) The map $\Gamma$ lands in $\VV$, and defines a quasi-isomorphism $\UU^- \rightarrow \VV$.
(ii) $\UU^+/\VV$ is isomorphic to $\B$ as a bimodule over $\bar\A = \bar\TT$.
(iii) The image of the Hochschild cocycle $\bar\Pi$ in $\mathit{CC}(\bar\TT,\UU^+/\VV) = \mathit{CC}(\bar\A,\B)$ is equal to $\bar\Phi$.
\end{theorem}

As consequence of this and Proposition \ref{th:plus-quasi-iso}, Corollary \ref{th:plus-u-quasi-iso}, we have the following commutative diagram in the homotopy category $H(\frakC)$ of $A_\infty$-bimodules, whose rows are exact triangles, and whose vertical maps are isomorphisms in $H(\frakC)$ (which means quasi-isomorphisms of bimodules):
\begin{equation} \label{eq:big-picture}
\xymatrix{
\Cone(\Gamma) \ar@{-->}[d] \ar[r] & \UU^- \ar[r]^-{\Gamma} \ar[d]^-{\Gamma} & \UU^+ \ar@{=}[d] \ar[r] & \Cone(\Gamma) \ar@{-->}[d] \\
\UU^+/\VV \ar@{=}[d] \ar[r] & \VV \ar@{^{(}->}[r] & \UU^+ \ar@{->>}[r] & \UU^+/\VV \ar@{=}[d] \\
\B \ar@{=}[d] \ar[r] & \Cone(\Phi) \ar@{-->}[d] \ar@{-->}[u] \ar[r] & \A \ar[d]^-{\Phi} \ar[u]_-{\Pi} \ar[r]^-{\Phi} & \B \ar@{=}[d] \\
\B \ar@{->>}[r] & \B^- \ar[r]^-{\Delta} & \B^+ \ar@{^{(}->}[r] & \B
}
\end{equation}
The maps marked $\hookrightarrow$ and $\twoheadrightarrow$ are inclusions and projections, respectively. Those marked $\dashrightarrow$ are obtained by completing a commutative square to a map of exact triangles.

\begin{remark} \label{th:cone}
In the general framework of triangulated categories, the morphisms completing a map of exact triangles are not unique (this is the well-known issue that cones are unique only in a weak sense). In the specific situation above, all squares are commutative on the cochain level, and hence the completions are unique up to homotopy. However, we will not really make use of this additional degree of precision.
\end{remark}

Now suppose that we add quasi-inverses of the quasi-isomorphisms $\Psi: \B^- \rightarrow \A^\vee$ and $\Xi: \UU^- \rightarrow \TT^\vee$ to the picture. From the middle two columns of \eqref{eq:big-picture} we get a commutative diagram in $H(\frakC)$, where the dotted arrow is again an isomorphism:
\begin{equation} \label{eq:big-picture-2}
\xymatrix{
\A^\vee \ar[drr]^-{\;\;\;\;\;\Pi^{-1} \circ \Gamma \circ \Xi^{-1}} \ar@{<-->}[dd] && \\
 && \A \\
\A^\vee \ar[urr]_-{\;\;\;\;\;\Phi^{-1} \circ \Delta \circ \Psi^{-1}} &&
}
\end{equation}
It seems plausible to think that the isomorphism is in fact the identity (this question makes sense thanks to Remark \ref{th:cone}), but we will not attempt to prove that.
In the absence of this piece of the overall puzzle, \eqref{eq:big-picture-2} shows the following, as originally stated in the Introduction:

\begin{corollary} \label{th:main}
There is a bimodule quasi-isomorphism $\Theta$ from $\A^\vee$ to itself, such that the bimodule homomorphisms $\Pi^{-1} \circ \Gamma \circ \Xi^{-1}$ and $\Phi^{-1} \circ \Delta \circ \Psi^{-1} \circ \Theta$ are homotopic. \qed
\end{corollary}

The proof of Theorem \ref{th:comparison} will take up the whole of this section. After initially setting up the comparison $\bar\A = \bar\TT$, we will study the structure of the bimodules $\UU^\pm$, leading to a proof of part (ii) of the theorem. Next, we consider the map $\Gamma$ and derive part (i). Finally, we compare the Hochschild cocycles appearing in (iii).

\begin{remark}
Up to homotopy, $\Gamma$ factors through the intermediate bimodules $\UU^c$. On the cohomology level, these factorizations agree with the ones constructed in terms of the fibre in Remark \ref{th:factoring}. It would be nice to show that the correspondence actually holds on the cochain level, but we have not tried to do that.
\end{remark}

\subsection{The $A_\infty$-algebras}
By definition, each intersection $l_i \cap l_j$ ($i<j$) consists of a single point $w_{ij}$, which lies in the interior of $U_{0,1}$ and is transverse. Let $S$ be a $(d+1)$-punctured disc, for some $d \geq 1$. Equip it with boundary conditions $(l_{i_0},\dots,l_{i_d})$ for some $i_0 < \cdots < i_d$. We consider holomorphic maps $v: S \rightarrow \C$ satisfying these boundary conditions, and with limits $w_{i_0,i_d}$ (over the end $\zeta_0$) respectively $w_{i_{k-1},i_k}$ (over the ends $\zeta_k$, $1 \leq k \leq d$).

\begin{lemma} \label{th:base-polygon}
There is precisely one such map $v$ for each $S$. Its image is entirely contained in the interior of $W_0$. Moreover, the associated linearized operator $D_v$ is invertible.
\end{lemma}

\proof We begin with the case $d = 1$, which is somewhat trivial. Since there is only intersection point $w_{i_0,i_1}$, any holomorphic strip $v$ is necessarily constant, and its linearized operator is the standard Dolbeault operator on $\R \times [0,1]$, with locally constant boundary conditions on $\R \times \{k\}$ given by the tangent space of $l_{i_k}$ at the point $w_{i_0,i_1}$. It is straightforward to see that this is invertible (for instance: the index is zero, and the kernel is zero by a maximum principle).

Assume from now on that $d \geq 2$. By assumption, the open subset $\C \setminus (W_0 \cup l_1 \cup \cdots \cup l_m)$ is connected and unbounded. Hence, we can apply the degree argument from Lemma \ref{th:v-semicircle} to show that no holomorphic map $v$ can enter this subset. Because of holomorphicity and the fact that the limits lie in the interior of $W_0$, the whole of $v(S)$ is then necessarily contained in that interior. \eqref{eq:index-formula} shows that $D_v$ has index zero, which by Corollary \ref{th:automatic-regularity} implies invertibility.

After these preliminaries, the rest is a deformation argument. In principle, we could have chosen our basis so that all the $l_i$ intersect in the same point. Then, the constant map would be a solution of our equation, and actually the only one, since any solution necessarily has energy zero by \eqref{eq:energy-action}. Now deform from this very specific situation back to our original choice of vanishing paths, while remaining within the class of simply intersecting bases (this is always possible). Since there are no nontrivial holomorphic strips with boundary on any two paths in such a basis, the parametrized moduli space of maps $v$ appearing in the deformation is compact. Corollary \ref{th:automatic-regularity} shows that the parametrized moduli space fibers smoothly over the deformation parameter. Hence, the number of maps $v$ for each parameter value is constant (this is the actual number, not a signed or modulo $2$ count, since regularity for each parameter value rules out birth-death processes). \qed

As pointed out in Remark \ref{th:transverse-thimbles}, it follows from our assumptions that the Lefschetz thimbles $L_i$ intersect transversally. Moreover, there is an obvious bijection $V_i \cap V_j \rightarrow L_i \cap L_j$ for $i<j$, which takes $x \in M$ to $y = (w_{ij},x) \in U_{0,1} \times M \iso \pi^{-1}(U_{0,1}) \subset E$. This induces an isomorphism of vector spaces
\begin{equation} \label{eq:cf-a-t}
\mathit{CF}_{\bar\A}(V_i,V_j) \iso \mathit{CF}_{\bar\TT}(L_i,L_j).
\end{equation}

Fix families of almost complex structures $J_{\bar{\A},i_0,\dots,i_d}^d$ on the fibre $M$ as needed to define $\bar\A$. When defining the corresponding families $J_{\bar{\TT},i_0,\dots,i_d}^d$ on $E$, we impose the following conditions:
\begin{equation} \label{eq:j-splits}
\parbox{28em}{The restriction of $J_{\bar{\TT},i_0,\dots,i_d}^d$ to $\pi^{-1}(W_0) = W_0 \times M$ is the product of the standard complex structure $i$ on the base and $J_{\bar{\A},i_0,\dots,i_d}^d$ on the fibre.}
\end{equation}
This is not a generic property, of course, so we'll have to explicitly discuss transversality issues later on.

\begin{lemma} \label{th:left-barrier-t}
Let $u: S \rightarrow E$ be a map which corresponds to a point of $\scrM_{\bar{\TT}}^d(y_0,\dots,y_d)$. Then $u(S) \subset \pi^{-1}(W_0) \iso W_0 \times M$.
\end{lemma}

\proof We apply the barrier argument from Section \ref{subsec:barrier} to $v = \pi(u)$, with $b>0$ small and $\kappa = 0$. More precisely, take $b$ so that all $w_{ij}$ lie in the interior of $W_b$, and so that $v$ intersects $B = \{p = b\}$ transversally. Crucially, $v$ is a holomorphic map on the subset $v^{-1}(W_0) \subset S$ of the domain, which contains $v^{-1}(B)$, hence \eqref{eq:barrier} applies. This is in fact a particularly simple special case of Example \ref{th:one-sided}. The argument there shows that $v(S)$ lies in the interior of $W_b$, which is stronger than what we needed. \qed

\begin{lemma} \label{th:bijection-of-a-t-spaces}
Projection to the fibre yields a bijection
\begin{equation} \label{eq:project-to-fibre}
\scrM_{\bar{\TT}}^d(y_0,\dots,y_d) \stackrel{\iso}{\longrightarrow} \scrM_{\bar{\A}}^d(x_0,\dots,x_d),
\end{equation}
where the $x_i$ and $y_i$ correspond as in \eqref{eq:cf-a-t}.
\end{lemma}

\proof The projection map is well-defined thanks to Lemma \ref{th:left-barrier-t}. In converse direction, let $u: S \rightarrow M$ be a map which contributes to $\scrM_{\bar\A}^d(x_0,\dots,x_d)$, where the boundary conditions are $(V_{i_0},\dots,V_{i_d})$. By Lemma \ref{th:base-polygon}, there is a unique holomorphic map $v: S \rightarrow \C$ with boundary conditions $(l_{i_0},\dots,l_{i_d})$ and appropriate limits. Moreover, the image of that map remains inside $W_0$. Hence, $(v,u)$ yields an element of $\scrM_{\bar\TT}^d(y_0,\dots,y_d)$, which is the unique preimage of $u$ under \eqref{eq:project-to-fibre}. \qed

Suppose, to make the notation slightly more homogeneous, that $d \geq 2$. To any map $u: S \rightarrow M$ in $\scrM_{\bar\A}(x_0,\dots,x_d)$ belongs an extended linearized operator $\tilde{D}_u: \tilde\E^1_u \rightarrow \E^0_u$ (called extended because it includes deformations of the conformal structure on $S$; the restriction to a fixed conformal structure would be given by the restriction to a subspace $\E^1_u \subset \tilde{\E}^1_u$ of codimension $d-2$, which would be denoted by $D_u$). We also have the ordinary linearized operator associated to $v: S \rightarrow \C$, denoted by $D_v: \E^1_v \rightarrow \E^0_v$. Finally, if we consider $(v,u)$ as a point of $\scrM_{\bar\TT}^d(y_0,\dots,y_d)$, its associated extended linearized operator $\tilde{D}_{(v,u)}: \tilde\E^1_{(v,u)} \rightarrow \E^0_{(v,u)}$ sits in a commutative diagram with short exact rows, as follows:
\begin{equation} \label{eq:compare-deformations}
\xymatrix{
0 \ar[r] & \E^1_v \ar[r] \ar[d]^-{D_v} & \tilde\E^1_{(v,u)} \ar[r] \ar[d]^-{\tilde{D}_{(v,u)}} & \tilde\E^1_u \ar[r] \ar[d]^-{\tilde{D}_u} & 0 \\
0 \ar[r] & \E^0_v \ar[r] & \E^0_{(v,u)} \ar[r] & \E^0_u \ar[r] & 0.
}
\end{equation}
The intuitive meaning is as follows: a first order deformation of $(S,v,u)$ induces, by forgetting about $v$, a first order deformation of $(S,u)$. The kernel of that forgetful map consists of first order deformations of $v$, for fixed $S$. As part of our assumptions when defining $\bar{\A}$, we required the choices of $J_{\bar{\A},i_0,\dots,i_d}^d$ to be such that the associated moduli spaces are regular, which by definition means that $\tilde{D}_u$ is onto. On the other hand, $D_v$ is invertible, see Lemma \ref{th:base-polygon}. In view of \eqref{eq:compare-deformations}, this implies that $\tilde{D}_{(v,u)}$ is onto as well. Therefore, choosing almost complex structures on $E$ as in \eqref{eq:j-splits} makes the spaces $\scrM_{\bar{\TT}}^d(y_0,\dots,y_d)$ regular. A similar (even a little simpler) argument applies to the case $d = 1$. Using these almost complex structures to define $\bar{\TT}$, we therefore find:

\begin{lemma}
With respect to the isomorphism from \eqref{eq:cf-a-t}, the $A_\infty$-structures $\mu_{\bar{\A}}^d$ and $\mu_{\bar{\TT}}^d$ coincide. \qed
\end{lemma}

\subsection{The choice of Hamiltonians\label{subsec:ccc}}
Our further argument will depend crucially on specific choices of the Hamiltonian perturbations used to define the bimodules $\UU^\pm$. Recall that $l_i \cap W_{-2}$ is the graph of a function $q_i$, with $q_i(p) = o_i$ equal to the ordinate for $p \geq 1$. Fix $c^\pm$ such that for all $i,j$,
\begin{equation} \label{eq:choose-cplusminus}
\begin{aligned}
& c^- < q_i(p) - q_j(p) \quad \text{for all $p \geq 0$,} \\
& c^+ > o_i - o_j.
\end{aligned}
\end{equation}
Choose $c^\pm$-displacement functions $h^\pm$ on $\C$ (see Definition \ref{def:displacement-1}) with the following additional properties:
\begin{equation} \label{eq:little-h}
\parbox{28em}{
$dh^-/dp = c^-$ for all $p \geq 0$. On the other hand, $h^+(p) = h^-(p)$ for $p \leq 1$, and $d^2 h^+/dp^2 > 0$ for $p \in (1,2)$.
}
\end{equation}
Figure \ref{fig:strips2} shows the effect of the resulting time-one maps $\phi_{h^\pm}$ on our basis of vanishing paths. Note in particular that for any $(i,j)$, the intersection $\phi_{h^+}(l_i) \cap l_j$ has a unique point lying in $W_1$. We denote it by $w_{ij}^+$. It is automatically transverse.

Take a function $H$ on the fibre $M$ satisfying \eqref{eq:b-general-position}. Starting from this and $h^\pm$, we choose $c^\pm$-displacement functions $H^\pm$ on $E$ (see Definition \ref{def:displacement-2}) as follows:
\begin{equation}
\parbox{28em}{
First, take a cutoff function $\psi$, satisfying $\psi(p) = 0$ for $p \leq -1$, and $\psi(p) = 1$ for $p \geq 0$. Then, the function $\tilde{H}^-(p,x) = h^-(p) + \psi(p) H(x)$ lies in the class from \eqref{eq:h-prime}. The actual $H^-$ should differ from this $\tilde{H}^-$ only on a compact subset of $E \setminus \partial E$, which is disjoint from $\pi^{-1}(W_0)$. Finally, we set $H^+(p,x) = h^+(p) + H(x)$ on $\pi^{-1}(W_1)$, and $H^+ = H^-$ on the rest of $E$.
}
\end{equation}
With respect to the trivialization $\pi^{-1}(W_0) \iso W_0 \times M$, we then have
\begin{equation} \label{eq:xh-split}
\begin{aligned}
& X_{H^-}|\pi^{-1}(W_0) = (c^- \partial_q,X_H), \\
& X_{H^+}|\pi^{-1}(W_0) = ((dh^+/dp) \partial_q,X_H).
\end{aligned}
\end{equation}
\begin{figure}
\begin{center}
\begin{picture}(10,0)%
\includegraphics{strips2.pstex}%
\end{picture}%
\setlength{\unitlength}{2723sp}%
\begingroup\makeatletter\ifx\SetFigFont\undefined%
\gdef\SetFigFont#1#2#3#4#5{%
  \reset@font\fontsize{#1}{#2pt}%
  \fontfamily{#3}\fontseries{#4}\fontshape{#5}%
  \selectfont}%
\fi\endgroup%
\begin{picture}(6038,3690)(3528,11) 
\put(9651,364){\makebox(0,0)[lb]{\smash{{\SetFigFont{10}{12}{\rmdefault}{\mddefault}{\updefault}{\color[rgb]{0,0,0}$\phi_{h^-}(l_i)$}%
}}}}
\put(4001, 89){\makebox(0,0)[lb]{\smash{{\SetFigFont{10}{12}{\rmdefault}{\mddefault}{\updefault}{\color[rgb]{0,0,0}$p = -2$}%
}}}}
\put(5401, 89){\makebox(0,0)[lb]{\smash{{\SetFigFont{10}{12}{\rmdefault}{\mddefault}{\updefault}{\color[rgb]{0,0,0}$p = -1$}%
}}}}
\put(6651, 89){\makebox(0,0)[lb]{\smash{{\SetFigFont{10}{12}{\rmdefault}{\mddefault}{\updefault}{\color[rgb]{0,0,0}$p = 0$}%
}}}}
\put(7901, 89){\makebox(0,0)[lb]{\smash{{\SetFigFont{10}{12}{\rmdefault}{\mddefault}{\updefault}{\color[rgb]{0,0,0}$p = 1$}%
}}}}
\put(9201, 89){\makebox(0,0)[lb]{\smash{{\SetFigFont{10}{12}{\rmdefault}{\mddefault}{\updefault}{\color[rgb]{0,0,0}$p = 2$}%
}}}}
\put(9651,2564){\makebox(0,0)[lb]{\smash{{\SetFigFont{10}{12}{\rmdefault}{\mddefault}{\updefault}{\color[rgb]{0,0,0}$\phi_{h^+}(l_i)$}%
}}}}
\put(9651,2114){\makebox(0,0)[lb]{\smash{{\SetFigFont{10}{12}{\rmdefault}{\mddefault}{\updefault}{\color[rgb]{0,0,0}$l_j$ ($j<i$)}%
}}}}
\put(9651,1439){\makebox(0,0)[lb]{\smash{{\SetFigFont{10}{12}{\rmdefault}{\mddefault}{\updefault}{\color[rgb]{0,0,0}$l_i$}%
}}}}
\put(9651,1064){\makebox(0,0)[lb]{\smash{{\SetFigFont{10}{12}{\rmdefault}{\mddefault}{\updefault}{\color[rgb]{0,0,0}$l_k$ ($k>i$)}%
}}}}
\end{picture}%
\caption{\label{fig:strips2}}
\end{center}
\end{figure}

\begin{lemma} \label{th:hminus}
All points of $\phi_{H^-}(L_i) \cap L_j$, for any $(i,j)$, lie in $\pi^{-1}(\C \setminus W_0)$.
\end{lemma}

\proof The subset $\pi^{-1}(W_0) \subset E$ is invariant under $X_{H^-}$. Hence, if $y \in \phi_{H^-}(L_i) \cap L_j$ lies in that subset, so does the corresponding flow line which joins $L_i$ and $L_j$. In view of \eqref{eq:xh-split}, this shows that $w = \pi(y)$ satisfies $w \in l_j$ and $w - ic^- \in l_i$. But this means that $q_j(p) - q_i(p) = c^-$ for $p = \mathrm{re}(w)$, contradicting \eqref{eq:choose-cplusminus}. \qed

\begin{lemma} \label{th:hplus}
Each point of $\phi_{H^+}(L_i) \cap L_j$ is of one of the following two kinds. (i) The first kind lie in $\pi^{-1}(\C \setminus W_0)$; these agree with the points of $\phi_{H^-}(L_i) \cap L_j$. (ii) The second kind project to the interior of $U_{1,2}$; these are in bijection with points of $\phi_H(V_i) \cap V_j$.
\end{lemma}

\proof As before, $\pi^{-1}(W_1)$ is invariant under $X_{H^+}$. Outside that subset, $H^+ = H^-$, so we immediately get part (i), and the absence of any intersection points in $\pi^{-1}(U_{0,1})$. For part (ii), we know that $\phi_{h^+}(l_i) \cap l_j \cap W_1 = \{w_{ij}^+\}$. In the fibre over this point, we have $\phi_{H^+}(L_i) \cap L_j \cap \pi^{-1}(w_{ij}^+) = \phi_H(V_i) \cap V_j$, whence the bijection. \qed

Since the choice of $H^-$ on any compact subset of $\pi^{-1}(\C \setminus W_0) \setminus \partial E$ is free, Lemma \ref{th:hminus} also implies that by a generic choice, we can achieve that the intersections $\phi_{H^-}(L_i) \cap L_j$ are transverse. Suppose that this has been done, and consider the intersections $\phi_{H^+}(L_i) \cap L_j$. The points of type (i) in Lemma \ref{th:hplus} are automatically again transverse. For type (ii), we know that the intersection point $w_{ij}^+ \in \phi_{h^+}(l_i) \cap l_j$ is transverse, and so is $\phi_H(V_i) \cap V_j$ by assumption, hence the same follows for the resulting points of $\phi_{H^+}(L_i) \cap L_j$. Versions of the same arguments show that the triple intersections $\phi_{H^\pm}(L_i) \cap \phi_{H^\pm}(L_j) \cap L_k$ and $L_i \cap L_j \cap \phi_{H^\pm}(L_k)$ will be empty for $i<j$ and arbitrary $k$, as long as $H^-$ is chosen generically, which we will again assume to be the case from now on.

\subsection{The bimodules\label{subsec:bimodule-comparison}}
When defining $\B$, we use the Hamiltonian $H$ introduced above, and some function $g$, which is assumed to be nonnegative everywhere. When defining $\UU^\pm$, we use $H^\pm$ and the same $g$. The families of almost complex structures should be chosen so that the following holds:
\begin{equation} \label{eq:j-plus-minus}
\parbox{28em}{
For any member of $J^{q|1|p}_{\UU^-,i_0,\dots,i_{p+q+1}}$, projection $\pi: E \rightarrow \C$ is pseudo-holomorphic over $U_{0,1}$. For $J^{q|1|p}_{\UU^+,i_0,\dots,i_{p+q+1}}$, we impose the more restrictive requirement that its restriction to $\pi^{-1}(W_0) \iso W_0 \times M$ should be the product of $i$ and $J^{q|1|p}_{\B,i_0,\dots,i_{p+q+1}}$.
}
\end{equation}
Again, since these are not generic conditions, some care needs to be exercised concerning transversality. Take a map $u: S \rightarrow M$ representing a point in $\scrM^{q|1|p}_{\UU^-}(y_0,\dots,y_{p+q+1})$, and consider its restriction to the output end, which is a solution of a Floer-type equation:
\begin{equation}
\left\{
\begin{aligned}
& u: (-\infty,s_0] \times [0,1] \longrightarrow E, \\
& \partial_s u + J^{0|1|0}_{\UU^-,i_0,i_{p+q+1},t}(u)(\partial_t u - g(t) X_{H^-}) = 0, \\
& u(s,0) \in L_{i_0}, \;\; u(s,1) \in L_{i_{p+q+1}}, \\
& \textstyle\lim_{s \rightarrow -\infty} u(s,t) = u_0(t),
\end{aligned}
\right.
\end{equation}
where $u_0$ is the $X_{H^-}$-flow line corresponding to the intersection point $y_0$. A generic point $(s,t)$ is {\em injective} in the following sense:
\begin{equation} \label{eq:injective}
\parbox{28em}{
First, $\partial_s u \neq 0$. Second, $u(s,t) \neq u_0(t)$. Finally, $u(s,t) \notin u((-\infty,s_0] \setminus \{s\},t)$.
}
\end{equation}
The existence of such a point follows from strong unique continuation, as in \cite[Section 4]{floer-hofer-salamon94}. Essentially, if one starts by assuming that no such point exists, this leads to the conclusion that $u(s,t) = u_0(t)$, and that equality then extends over the whole of $S$ (not just on the strip-like end; it is crucial that $(s,t) \mapsto u_0(t)$ is indeed a solution of the relevant equation), violating the assumption made when choosing $H$. Now, since the limit $u_0$ necessarily lies in $E \setminus \pi^{-1}(W_0)$, the choice of almost complex structure $J^{0|1|0}_{\UU^-,i_0,i_{p+q+1},t}$ at points $u(s,t)$ with $s \ll 0$ is essentially unconstrained, and this allows one to obtain regularity in the same way as for the ordinary Floer equation. Of course, we could have run the same argument using the opposite end $s \gg 0$. For $\UU^+$, the same arguments yield the regularity of $\scrM_{\UU^+}(y_0,\dots,y_{p+q+1})$ whenever $y_0$ or $y_{p+1}$ is of type (i). The remaining case will have to be dealt with more concretely, which we begin to do now.

\begin{lemma} \label{th:submodule}
Suppose that (in the terminology of Lemma \ref{th:hplus}) $y_0 \in \phi_{H^+}(L_{i_0}) \cap L_{i_{p+q+1}}$ is of type (ii), $y_{p+1} \in \phi_{H^+}(L_{i_p}) \cap L_{i_{p+1}}$ is of type (i), and we have arbitrary auxiliary points $y_k \in L_{i_{k-1}} \cap L_{i_k}$ (for $1 \leq k \leq p$ and $p+2 \leq k \leq p+q+1$). Then $\scrM^{q|1|p}_{\UU^+}(y_0,\dots,y_{p+q+1}) = \emptyset$.
\end{lemma}

\proof Assume the contrary, and let $u: S \rightarrow E$ be a point in that moduli space, where by definition $S = (\R \times [0,1]) \setminus \{\zeta_1,\dots,\zeta_p,\zeta_{p+2},\dots,\zeta_{p+q+1}\}$. Consider the projection $v = \pi(u): S \rightarrow \C$. On $v^{-1}(W_0) \subset S$ this satisfies an equation \eqref{eq:v-equation}, where the boundary values are $l_{i_0},\dots,l_{p+q+1}$, and the inhomogeneous term $X_k = (dh^+/dp) \partial_q \otimes g(t) \mathit{dt}$. This is a consequence of \eqref{eq:j-plus-minus} and \eqref{eq:xh-split}.

We will apply a barrier argument to $v = \pi(u)$, taking $b<1$ so that $1-b$ is small, and $v$ intersects $B = \{p = b\}$ transversally. Since $dh^+/dp = c^-$ on the interval $[0,1]$, we can take $\kappa = c^- g(t) \mathit{dt}$ and \eqref{eq:barrier-def} will be satisfied. For each $k \neq 0,p+1$ we know by definition of a simply intersecting basis that $\pi(y_k)$ lies in the interior of $U_{0,1}$. Hence, we may choose $b$ such that all these points lie to the left of the barrier. By assumption, $\pi(y_{p+1})$ lies to the left of the barrier, while $\pi(y_0)$ lies to its right. This puts us in the situation from Example \ref{th:1-cross}. The requirement \eqref{eq:sign-of-translation} amounts to
$c^- < q_{i_{p+q+1}}(b) - q_{i_0}(b)$, which is satisfied because of \eqref{eq:choose-cplusminus}. \qed

For each $(i,j)$, take $\mathit{CF}_\VV(L_i,L_j) \subset \mathit{CF}_{\UU^+}(L_i,L_j)$ to be the subspace generated by those $y$ of type (i). Lemma \ref{th:submodule} says that the direct sum of these subspaces, denoted by $\VV \subset \UU^+$, is a sub-bimodule over $\bar\TT$. Our next goal is to determine the structure of the quotient $\UU^+/\VV$ explicitly.

\begin{lemma} \label{th:plus-preliminary}
Suppose that both $y_0 \in \phi_{H^+}(L_{i_0}) \cap L_{i_{p+q+1}}$ and $y_{p+1} \in \phi_{H^+}(L_{i_p}) \cap L_{i_{p+1}}$ are of type (ii), and we have arbitrary auxiliary points $y_k \in L_{i_{k-1}} \cap L_{i_k}$ (for $1 \leq k \leq p$ and $p+2 \leq k \leq p+q+1$). Then, for any map $u: S \rightarrow M$ representing a point in $\scrM^{q|1|p}_{\UU^+}(y_0,\dots,y_{p+q+1})$, the whole image $u(S)$ is contained in the interior of $\pi^{-1}(W_0)$.
\end{lemma}

\proof This is another barrier argument, taking the same $\kappa$ as before, but now using some small $b>0$. One can choose this so that all limits $\pi(y_k)$ lie to the right of the barrier $B$, and so that $v = \pi(u)$ intersects $B$ transversally. This puts us in the situation of Example \ref{th:one-sided}, and the result follows immediately. \qed

We temporarily return to considering only the base. Let $S$ be a surface of the form \eqref{eq:punctured-strip}, for some $p,q \geq 0$. Equip it with boundary conditions $(l_{i_0},\dots,l_{i_{p+q+1}})$, where $i_0 < \cdots < i_p$ and $i_{p+1} < \cdots < i_{p+q+1}$. We consider solutions $v: S \rightarrow \C$ of the equation
\begin{equation} \label{eq:hplus-equation}
\bar\partial v = (X_{h^+} \otimes g(t) \mathit{dt})^{0,1},
\end{equation}
with limits corresponding to $w^+_{i_0,i_{p+q+1}}$ and $w^+_{i_p,i_{p+1}}$ as $s \rightarrow \pm \infty$, respectively, and $w_{i_k,i_{k+1}}$ on the other ends. In analogy with Lemma \ref{th:base-polygon}, we have:

\begin{lemma} \label{th:base-polygon-2}
There is precisely one such map $v$ for each $S$. Its image is entirely contained in the interior of $W_0$. Moreover, the associated linearized operator $D_v$ is invertible.
\end{lemma}

\proof The case $p + q = 0$ is again trivial, since a solution of \eqref{eq:hplus-equation} with both limits equal to $w^+_{i_0,i_1}$ is necessarily constant.

Assume from now on that $p + q > 0$. By the same barrier argument as in Lemma \ref{th:plus-preliminary}, any solution $v$ must remain inside $W_0$. Invertibility of $D_v$ follows as before by combining \eqref{eq:index-formula} and Corollary \ref{th:automatic-regularity}. Moreover, the moduli space of maps $v$ is always compact, for the following reasons. Bubbling along the ends $\zeta_k$ for $k = 1,\dots,p,p+2,\dots,p+q+1$ is ruled out since it would lead to a non-constant holomorphic strip with boundary on $(l_{i_k},l_{i_{k+1}})$, which is impossible since these two paths have only a single intersection point. Consider one of the two other ends, say $\zeta_{p+1}$ which corresponds to $s \rightarrow \infty$. Bubbling at that end would lead to a non-stationary solution of \eqref{eq:hplus-equation} defined on the whole of $\R \times [0,1]$, and with boundary on $(l_{i_p},l_{i_{p+1}})$. The intersection $\phi_{h^+}(l_{i_p}) \cap l_{i_{p+1}}$ can consist of more than just one point, so such solutions are not generally ruled out. However, all intersection points except for $w^+_{i_p,i_{p+1}}$ lie outside $W_0$, hence can't occur as limit bubbles of solution $v$, since those stay inside that subset. The same argument applies to the remaining end.

\begin{figure}
\begin{center}
\includegraphics{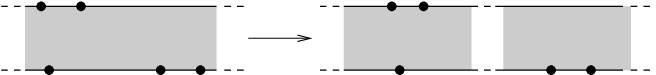}
\end{center}
\caption{\label{fig:reduce}}
\end{figure}

At this point, the standard deformation argument shows that the number of maps $v$ is independent of the particular basis of vanishing paths, as well as of the conformal structure of $S$. Denote this (unsigned) number by $N^{q|1|p}_{i_0,\dots,i_{p+q+1}} \geq 0$. We can get additional information about this by degenerating the conformal structure. For instance,
\begin{equation} \label{eq:reduce-split}
N^{q|1|p}_{i_0,\dots,i_{p+q+1}} =
N^{q'|1|p'}_{i_0,\dots,i_{p'},i_{p + q'' + 1},\dots,i_{p+q+1}}
N^{q''|1|p''}_{i_{p'},\dots,i_{p + q'' + 1}}
\end{equation}
for $p = p' + p''$, $q = q' + q''$. This corresponds to a neck-stretching degeneration, in the limit of which $S$ splits into two pieces $S'$ and $S''$ of the same type, each carrying part of the marked points (Figure \ref{fig:reduce}). Of course, \eqref{eq:reduce-split} depends on suitable compactness and gluing results, but these are relatively straightforward in the present context.

There is one more deformation which is useful for our purpose, depending on an additional parameter $r \in [0,1]$. We keep $S$ fixed, but change $g(t)\mathit{dt}$ to $(1-r)\,g(t) \mathit{dt}$, and simultaneously replace the boundary conditions $l_{i_k}$ by $\phi_{r h^+}(l_{i_k})$ for all $k \leq p$, while leaving all others constant. Note that for any value of $r$, the limits of solutions $v$ at the ends $\zeta_0$, $\zeta_{p+1}$ still correspond to points of $\phi_{h^+}(l_{i_0}) \cap l_{p+q+1}$ and $\phi_{h^+}(l_{i_p}) \cap l_{i_{p+1}}$, respectively. The barrier argument used before applies with minor modifications: setting $\kappa^r = r c^- g(t) \mathit{dt}$, and taking into account the modified position of the vanishing paths, one finds that \eqref{eq:sign-of-translation} is again satisfied. Finally, the same regularity and compactness results continue to hold, so the number of maps $v$ is constant in $r$.

From this, the general result that $N^{q|1|p}_{i_0,\dots,i_{p+q+1}} = 1$ can be derived in the following steps. First, using \eqref{eq:reduce-split} one can reduce the situation inductively to the case where $p+q = 1$. For simplicity of notation, let's consider only one of the two significant cases, namely where $p = 0$ and $q = 1$. There, the deformation we have just described allows one to compute the number $N^{1|1|0}_{i_0,i_1,i_2}$ by an ordinary count of holomorphic maps $v: S \rightarrow \C$, where $S = (\R \times [0,1]) \setminus \{(0,1)\}$ is a three-punctured disc with boundary conditions $\phi_{h^+}(l_{i_0})$, $l_{i_1}$, $l_{i_2}$, and limits given by the points $w_{i_0,i_1}^+$, $w_{i_1,i_2}$, $w_{i_0,i_2}^+$. By a well-known uniformization argument, this comes down to counting immersed triangles, and inspection shows that there is exactly one of those. \qed

Lemma \ref{th:hplus} gives a bijection between points of $\phi_H(V_i) \cap V_j$ and points of type (ii) in $\phi_{H^+}(L_i) \cap L_j$, hence an isomorphism of vector spaces
\begin{equation} \label{eq:cf-b-u}
\mathit{CF}_\B(V_i,V_j) \iso \mathit{CF}_{\UU^+}(L_i,L_j)/\mathit{CF}_{\VV}(L_i,L_j).
\end{equation}
By the same argument as in Lemma \ref{th:bijection-of-a-t-spaces}, this time based on Lemmas \ref{th:plus-preliminary} and \ref{th:base-polygon-2}, one sees that if $y_0$ and $y_{p+1}$ are of type (ii), then projection to the fibre yields a bijection
\begin{equation} \label{eq:project-to-fibre-2}
\scrM_{\UU^+}^{q|1|p}(y_0,\dots,y_{p+q+1}) \stackrel{\iso}{\longrightarrow} \scrM_{\B}^d(x_0,\dots,x_{p+q+1}).
\end{equation}
Moreover, the regularity of the moduli space on the right hand again implies that of the one on the left, which completes our transversality argument for $\UU^+$. The conclusion is the following statement, which is part (ii) of Theorem \ref{th:comparison}:

\begin{lemma}
With respect to the isomorphism from \eqref{eq:cf-b-u}, the $A_\infty$-bimodule structures $\mu_{\UU^+/\VV}^{q|1|p}$ and $\mu_{\B}^{q|1|p}$ coincide. \qed
\end{lemma}

\subsection{The bimodule map}
When defining the map $\Gamma: \UU^- \rightarrow \UU^+$, we choose the worldsheet structures as follows. Let $\eta$ be a cutoff function, satisfying $\eta(s) = 0$ for $s \ll 0$, $\eta(s) = 1$ for $s \gg 0$ and $\eta'(s) \geq 0$ for all $s$. For the inhomogeneous term, set
\begin{equation} \label{eq:k-and-tilde-k}
K = \big( \eta(s) H^- + (1-\eta(s)) H^+ \big) g(t) \mathit{dt}.
\end{equation}
The underlying worldsheet has $\gamma = \big(\eta(s) c^- + (1-\eta(s)) c^+ \big) g(t) \mathit{dt}$, which satisfies $d\gamma \leq 0$ because of the condition on the derivative of $\eta$, together with the nonnegativity of $g$. Finally, for the almost complex structure, we make the following requirement:
\begin{equation} \label{eq:j-gamma}
\parbox{28em}{
For any member of $J^{q|1|p}_{\Gamma,i_0,\dots,i_{p+q+1}}$, projection $\pi: E \rightarrow \C$ is pseudo-holomorphic over $W_0$.
}
\end{equation}
Transversality is established by an argument similar to the one for $\UU^\pm$. Let $u: S \rightarrow E$ be any map representing a point in $\scrM_{\Gamma}^{q|1|p}(y_0,\dots,y_{p+q+1})$. Let $u_{p+1}$ be the $X_{H^-}$-flow line corresponding to $y_{p+1}$. Since that lies in $\pi^{-1}(\C \setminus W_0)$ by Lemma \ref{th:hminus}, the constraint \eqref{eq:j-gamma} does not affect our choices of almost complex structures at points $u(s,t)$ for $s \gg s_0$. Then, either a version of the injectivity property \eqref{eq:injective} holds for such points, or else $u(s,t) = u_{p+1}(t)$ on the whole of $S$. The second alternative actually happens, but only if $p+q = 0$. Note that on $\pi^{-1}(\C \setminus W_1)$ we have $H^- = H^+$, hence the equation for $u$ reduces to an ordinary Floer-type equation, except of course that we don't divide out by translation. It is a standard fact that stationary solutions are regular.

\begin{lemma} \label{th:gamma-into-v}
Take $y_0 \in \phi_{H^+}(L_{i_0}) \cap L_{i_{p+q+1}}$ which is of type (ii), and arbitrary $y_{p+1} \in \phi_{H^-}(L_{i_p}) \cap L_{i_{p+1}}$ as well as $y_k \in L_{i_{k-1}} \cap L_{i_k}$, $k = 1,\dots,p,p+2,\dots,p+q+1$. Then $\scrM_\Gamma^{q|1|p}(y_0,\dots,y_{p+q+1}) = \emptyset$.
\end{lemma}

\proof On $\pi^{-1}(U_{0,1})$ we have $K = H^- g(t) \mathit{dt}$, hence $X_K = (c^- \partial_q, X_H)$. Therefore, the same barrier argument as in Lemma \ref{th:submodule} applies. \qed

This says that $\Gamma$ is in fact a map from $\UU^-$ to the sub-bimodule $\VV \subset \UU^+$. We will now consider the linear terms in this map, which are chain homomorphisms
\begin{equation} \label{eq:linear-gamma}
\Gamma^{0|1|0}: \mathit{CF}_{\UU^-}(L_i,L_j) \longrightarrow \mathit{CF}_\VV(L_i,L_j).
\end{equation}
Both sides have the same generators, corresponding to points of $\phi_{H^-}(L_i) \cap L_j$. By construction, $H^+ \geq H^-$ everywhere, which means that $\partial_s K(\partial_t) = \eta'(s) (H^- - H^+) g(t) \leq 0$. Since all other terms vanish, the curvature \eqref{eq:big-curvature} is nonpositive. Hence, if $u$ is any solution of the equation underlying \eqref{eq:linear-gamma} with limits $y^\pm$ as $\pm s \rightarrow -\infty$, then by \eqref{eq:e-geom} and \eqref{eq:e-a} we have
\begin{equation}
0 \leq E^{\geom}(u) \leq E^{\top}(u) = A(y^+) - A(y^-).
\end{equation}
The only solutions with energy zero are those with $\partial_s u$ vanishing, which we already considered above. The upshot is that \eqref{eq:linear-gamma} is the identity plus another term which is strictly lower triangular with respect to the action filtration. As a consequence, we get part (i) of Theorem \ref{th:comparison}:

\begin{lemma}
$\Gamma: \UU^- \rightarrow \VV$ is a quasi-isomorphism of $A_\infty$-bimodules. \qed
\end{lemma}

\subsection{The Hochschild cocycle}
We first consider the relevant piece of data in the fibre, which is $\bar\Phi \in \mathit{CC}(\bar\A,\B)$. In order to define this, we require that our original $H$ should satisfy \eqref{eq:generic-3}, and choose a suitable one-form on the upper half-plane $\R \times \R^+ \subset \C$. It will be convenient for us to suppose that, in the exponential coordinates $\R \times [0,1] \rightarrow \R \times \R^+$ from \eqref{eq:half-plane-end}, this is of the form $\gamma = \xi(s) c^+  g(t) \mathit{dt}$ where $\xi(s) = 1$ for $s \ll 0$, $\xi(s) = 0$ for $s \gg 0$, and $\xi'(s) \leq 0$ everywhere. For $\bar\Pi \in \mathit{CC}(\bar\TT,\UU^+)$ we use the same $\gamma$, and define the inhomogeneous term in the following two-step process:
\begin{equation} \label{eq:k-tilde-k}
\begin{aligned}
& \tilde{K} = \xi(s) H^+ g(t) \mathit{dt}, \\
& K = \tilde{K} + (\text{\it compactly supported perturbation}),
\end{aligned}
\end{equation}
where the last-mentioned term is supported in a compact subset of $(S \setminus \partial S) \times (E \setminus \partial E)$, which moreover is disjoint from $\pi^{-1}(W_0)$. As for the almost complex structures, we make the following assumptions:
\begin{equation}
\parbox{28em}{
The restriction of $J^d_{\bar\Pi,i_0,\dots,i_d}$ to $\pi^{-1}(W_0) \iso W_0 \times M$ should be the product of $i$ and $J^d_{\bar\Phi,i_0,\dots,i_d}$.
}
\end{equation}
Let's consider the resulting moduli spaces $\scrM_{\bar\Pi}^d(y_0,\dots,y_d)$. If $y_0 \in \phi_{H^+}(L_{i_0}) \cap L_{i_d}$ is of type (i), all maps in such a moduli space must leave $\pi^{-1}(W_0)$, which makes it easy to establish transversality because of the additional freedom provided by \eqref{eq:k-tilde-k}. The type (ii) case will be addressed by explicitly determining the moduli spaces. The starting point for that is the following:

\begin{lemma} \label{th:plus-preliminary-cc}
Suppose that $y_0$ is of type (ii). Then, for any map $u: S \rightarrow M$ representing a point in $\scrM^d_{\bar\Pi}(y_0,\dots,y_d)$, the whole image $u(S)$ is contained in the interior of $\pi^{-1}(W_0)$.
\end{lemma}

\proof We place a barrier $B$ at some small $b>0$, setting $\kappa = c^- \xi(s) g(t) \mathit{dt}$. All limits lie to the right of that barrier. Of course, the barrier is not closed, but $d\kappa = c^- \xi'(s) g(t) \mathit{ds} \wedge \mathit{dt}$ is nonnegative everywhere, since $c^- < 0$. Hence, the argument from Example \ref{th:one-sided} still applies, see also Remark \ref{th:sub-closed-barrier}. \qed

Let's temporarily restrict to taking $S$ to be the whole upper half plane, with our given $\gamma$ and boundary condition given by some $l_{i_0}$, and consider the associated equation $\bar\partial v = (X_{h^+} \otimes \xi(s) g(t) \mathit{dt})^{0,1}$, with limit $w^+_{i_0,i_0}$. For this equation,

\begin{lemma} \label{th:base-polygon-3}
The only solution is the constant map $v \equiv w^+_{i_0,i_0}$. Moreover, that solution is regular.
\end{lemma}

\proof By construction, $X_{h^+} = 0$ at the point $w = w^+_{i_0,i_0} \in l_{i_0}$. Hence, the constant map is indeed a solution of our equation. One easily computes that its index is zero, and then regularity follows as usual from Corollary \ref{th:automatic-regularity}.

Suppose that we have another solution $v$. The same barrier argument as in Lemma \ref{th:plus-preliminary-cc} shows that its image is contained in $W_0$. Suppose that we temporarily modify $h^+$ by a constant, so that $h^+(p) = 0$ at $p = \mathrm{re}(w)$, which of course does not affect $X_{h^+}$. Then, the fact that $h^+$ has vanishing derivative at $p$, and is convex on $[0,\infty)$ by \eqref{eq:little-h}, implies that the curvature term \eqref{eq:simplified-small-curvature} is nonnegative. At the same time, $A(w) = 0$ by definition \eqref{eq:little-action}, so $E^{\geom}(v) \leq E^{\top}(v) = 0$ by \eqref{eq:energy-action}. This implies that $dv = X_{h^+} \otimes \gamma$, which because of the limit $w$ means that $v$ must be constant. \qed

By the same gluing argument as in Lemma \ref{th:base-polygon-2}, one can deduce from Lemma \ref{th:base-polygon-3} that for any $S$ of the form \eqref{eq:punctured-plane}, with boundary conditions $l_{i_0},\dots,l_{i_d}$ for $i_0 < \cdots < i_d$, we have a unique solution of the associated equation with limits $(w_{i_0,i_d}^+,w_{i_0,i_1},\dots,w_{i_{d-1}},w_{i_d})$, and that these solutions are regular. This implies that projection to the fibre yields a bijection
\begin{equation} \label{eq:project-to-fibre-3}
\scrM_{\bar\Pi}^d(y_0,\dots,y_d) \stackrel{\iso}{\longrightarrow} \scrM_{\bar\Phi}^d(x_0,\dots,x_d)
\end{equation}
whenever $y_0$ is of type (ii), and $x_0$ is the corresponding point of $\phi_H(V_{i_0}) \cap V_{i_d}$. Regularity of the moduli space on the right hand side then implies that on the left, and we arrive at the following conclusion, which is part (iii) of Theorem \ref{th:comparison}:

\begin{lemma}
The image of $\bar\Pi^d$ in the quotient
\begin{equation}
\mathit{CC}^d(\bar\TT,\UU^+/\VV) = \mathit{CC}^d(\bar\TT,\UU^+)/\mathit{CC}^d(\bar\TT,\VV) \iso \mathit{CC}^d(\bar\A,\B)
\end{equation}
agrees with $\bar\Phi^d$. \qed
\end{lemma}


\end{document}